\documentclass[11pt,reqno]{amsart}
\usepackage{amsmath,amssymb,latexsym,esint,cite,mathrsfs}
\usepackage{verbatim,wasysym}
\usepackage[left=1.8cm,right=1.8cm,top=2.4cm,bottom=2.4cm]{geometry}
\usepackage{tikz,enumitem,graphicx, subfig, microtype, color}
\usepackage{epic,eepic}

\usepackage[colorlinks=true,urlcolor=blue, citecolor=red,linkcolor=blue,
linktocpage,pdfpagelabels, bookmarksnumbered,bookmarksopen]{hyperref}
\usepackage[english]{babel}
\usepackage[hyperpageref]{backref}
\usepackage{appendix}
\usepackage{bm}
\DeclareMathOperator{\diag}{diag}%
\DeclareMathOperator{\ind}{ind}
\DeclareMathOperator{\nul}{nul}

\usepackage{array,booktabs,tabularx,mathtools,needspace}

\numberwithin{equation}{section}

\newtheorem{thm}{Theorem}[section]
\newtheorem{lem}[thm]{Lemma}
\newtheorem{cor}[thm]{Corollary}
\newtheorem{Prop}[thm]{Proposition}

\newtheorem{Rem}[thm]{Remark}

\newtheorem*{Rem*}{Remark}

\newcommand{\R}{\mathbb R}
\newcommand{\Om}{\Omega}\newcommand{\Id}{\mathrm I}
\newcommand{\cD}{\mathcal D^{1,2}}\newcommand{\cK}{\mathcal K}
\newcommand{\cN}{\mathcal N}\newcommand{\cS}{\mathcal S}
\newcommand{\cP}{\mathcal P}\newcommand{\cQ}{\mathcal Q}
\newcommand{\cB}{\mathcal B}\newcommand{\cG}{\mathcal G}

\newcommand{\one}{\mathbf 1}

\newcommand{\Span}{\operatorname{span}}

\begin{document}

	\title[The four-dimensional Brezis--Nirenberg problem]
	{Fine spectral analysis of concentrating solutions
to Brezis--Nirenberg problem in dimension $4$}
	
\author[S. Cingolani]{Silvia Cingolani}
	\address{\noindent Silvia Cingolani  \newline
		Dipartimento di Matematica, Universit\`{a} degli Studi di Bari Aldo Moro,\newline
		Via Orabona 4, 70125 Bari, Italy.}\email{silvia.cingolani@uniba.it}

	\author[M. Yang]{Minbo Yang}
	\address{\noindent Minbo Yang  \newline
		School of Mathematical Sciences, Zhejiang Normal University,\newline
		Jinhua, China.}
\email{mbyang@zjnu.edu.cn}
	
	\author[S. Zhao]{Shunneng Zhao}
	\address{\noindent Shunneng Zhao  \newline
	School of Mathematical Sciences, Zhejiang Normal University,\newline
		Jinhua, China.}
		\email{snzhao@zjnu.edu.cn}



	\allowdisplaybreaks

\date{}

\begin{abstract}
		{\small
This paper is devoted to the qualitative properties and Morse indices of
concentrating solutions to the Brezis--Nirenberg problem
\[
    -\Delta u=u^3+\epsilon u, \quad u>0, \quad\text{in }\Omega,
    \qquad
    u=0 \quad \text{on }\partial\Omega,
\]
where $\Omega\subset\mathbb R^4$ is a smooth bounded domain, $\epsilon>0$ is small parameter.
For a family of $m$-bubble solutions with distinct interior
concentration points and comparable concentration scales,
we establish precise asymptotic expansions for the first
$6m$ eigenvalues of the associated linearized eigenvalue
problem and determine the corresponding eigenfunction profiles.
}
	\end{abstract}

\vspace{3mm}
\maketitle

\section{Introduction}\label{sec:introduction}
\subsection{Motivation}
Let us consider the well-known Brezis-Nirenberg problem, i.e.
\begin{equation}\label{intro:classicalBN}
\begin{cases}
-\Delta u=|u|^{2^{\ast}-2}u+\epsilon u &\text{in }\Om,\\
u=0 &\text{on }\partial\Om.
\end{cases}
\end{equation}
where $\Omega$ is a bounded regular domain in $\mathbb{R}^N$ with $N\geq3$, $\epsilon$ is a positive parameter and
$2^{\ast}:=\frac{2N}{N-2}$ is the critical exponent for the Sobolev embedding. In a celebrated paper  \cite{Brezis-1983}, Brezis and Nirenberg
proved that, for $N\ge4$, a positive solution exists whenever
$0<\epsilon<\lambda_1(\Om)$, where $\lambda_1(\Om)$ is the first
Dirichlet eigenvalue of $-\Delta$. At $\epsilon=0$, the geometry of
the domain becomes decisive: Poho\v{z}aev's identity
\cite{Pohozaev-1965} rules out positive solutions on star-shaped
domains, whereas Bahri and Coron \cite{Bahri-1988} established
existence under a nontrivial homology condition. Thus the existence
of solutions for $\epsilon>0$ raises a further question: how do they
behave as the perturbation vanishes, especially when the limiting
problem has no positive solution?
This question leads to the study of concentration, whose basic
features are the number, locations and rates of the bubbles.
Following Brezis and Peletier \cite{BP}, Han \cite{HANZCHAO} and
Rey \cite{Rey-1989} studied positive families approaching the best
Sobolev quotient. For the Brezis--Nirenberg problem in dimensions
$N\ge5$, they identified a single interior concentration point,
characterized it as a critical point of the Robin function, and
determined the concentration rate and limiting profile. Conversely,
Rey \cite{Rey-1990} constructed single-bubble solutions near
nondegenerate critical points of the Robin function for $N\ge5$.
These results show that concentration is governed by geometric data
associated with the Dirichlet Green function.
For solutions with several bubbles, the location of each bubble is
also affected by the others. Musso and Pistoia
\cite{Musso-Pistoia-2002} constructed positive multi-bubble solutions
for $N\ge5$ using stable critical points of a reduced function that
contains both the Robin function and the Green interactions. The lower-dimensional cases require a different treatment.
In $n=3$, del Pino, Dolbeault and Musso \cite{DDM} and
Musso and Salazar \cite{MS-2018} obtained single-bubble and
multi-bubble constructions, respectively, as the linear parameter
approaches a suitable positive value. Complementary information
comes from the fine multibubble analysis of K\"onig and Laurain
\cite{K-L-2024,K-L-2026} for $N\ge4$, which relates concentration
rates and limiting configurations to the Green--Robin interaction
matrix. Cao, Luo and Peng \cite{CLP} also established local
uniqueness and counting results for $N\ge6$ under suitable
assumptions on the reduced data.

 In the case $N=4$, Pistoia and Rocci \cite{PR} constructed single-bubble solutions of
\begin{equation}\label{intro:BN}
\begin{cases}
-\Delta u=u^3+\epsilon u &\text{in }\Om,\\
u>0 &\text{in }\Om,\\
u=0 &\text{on }\partial\Om.
\end{cases}
\end{equation}
The existence theory for positive multi-bubble solutions in dimension four has developed more recently.
Pistoia, Rago and Vaira \cite{P-R-V} constructed positive multi-bubble solutions to (1.2) associated with stable critical sets of the least eigenvalue of the Green-Robin interaction matrix, under the assumption that this eigenvalue is positive on these sets. Their result is governed by the smallest eigenvalue of the Green-Robin interaction matrix.
The concentration scales of these solutions are exponentially small in \(1/\epsilon\), whereas those in higher dimensions decay algebraically. This change in the concentration rate leads to a different balance between the perturbation and the bubble interactions, which plays a central role in the refined spectral analysis below.

The preceding results establish the existence of concentrating solutions and describe their asymptotic behavior. A natural next question is how the concentration points and the interactions between bubbles determine the Morse index and nondegeneracy of these solutions. Addressing this question requires a refined analysis of the spectrum of the associated linearized operator.
 For the slightly subcritical problem, Bahri, Li and Rey \cite{B-L-R} related the Morse index of multi-bubble solutions to a finite-dimensional reduced function. Grossi and Pacella \cite{GP05} obtained eigenvalue and eigenfunction asymptotics in the Sobolev-minimizing single-bubble regime, while Choi, Kim and Lee \cite{CKL-1} developed the corresponding spectral analysis for positive multi-bubble solutions. For the Brezis-Nirenberg problem \eqref{intro:classicalBN}, Takahashi \cite{T-1} studied the single-bubble spectrum, obtaining asymptotic estimates for the first eigenpair when $N\ge5$ and refined expansions for the translation and dilation modes when $N\ge6$. In particular, the leading corrections to the translation eigenvalues are determined by the Hessian of the Robin function at the concentration point, establishing a direct connection between the concentration geometry and the linearized spectrum. Other concentration
regimes have recently been studied by Li, Vaira, Wei and Wu \cite{LVWW},
who constructed bubbling solutions in dimensions four and five as the
linear parameter approaches Dirichlet eigenvalues.
He, Liu, Wei and Wu \cite{HLWW} used an inverse reduction argument to obtain refined one-bubble asymptotics for
$N\ge4$ as the linear parameter approaches a positive limit.

More recently, Yang and Zhao \cite{YZ,YZmulti} investigated the spectral properties of single-bubble and multi-bubble solutions with a subcritical perturbation \(\epsilon u^q\). Their general spectral convergence theorem \cite[Theorem~1.2]{YZmulti} covers the case \(N=4\), \(q=1\), whereas their refined multi-bubble analysis for \(m\ge2\) requires \(q>1\) in dimension four. Thus, although the translation and dilation eigenvalues are known to converge to \(1\), their leading corrections in the case \(N=4\), \(q=1\) remain undetermined by these results. Such corrections are needed to distinguish eigenvalues lying above or below \(1\) and hence to investigate the Morse index and nondegeneracy of the concentrating solutions. To the best of our knowledge, a complete description of these refined spectral asymptotics and the corresponding exterior eigenfunction profiles has not yet been obtained for positive separated multi-bubble solutions of the four-dimensional Brezis-Nirenberg problem as \(\epsilon\to0\).

The purpose of this paper is to determine the refined spectral behavior of positive separated single-bubble and multi-bubble families for \eqref{intro:BN}. In particular, we identify the different orders of the translation, relative dilation and
common dilation eigenvalues. The resulting expansions relate
the Morse index to the Hessian of the least interaction
eigenvalue and yield nondegeneracy when that Hessian is
nondegenerate. They therefore connect the four-dimensional
existence theory with the qualitative properties of the
solutions it produces.

\subsection{Main results}\label{subsec:main-results}
We first introduce the notation needed to state our main results.
Let $G$ be the Dirichlet Green function of $-\Delta$ in $\Om$,
and let $H$ be its regular part. We write
\begin{equation}\label{eq:bubble}
U_{\delta,\xi}(x)
=\frac{2\sqrt2\,\delta}{\delta^2+|x-\xi|^2},
\qquad
U=U_{1,0},
\qquad
G(x,y)=\frac{1}{4\pi^2|x-y|^2}-H(x,y),
\end{equation}
where $\delta>0$ and $\xi\in\Om$, and denote by
$\tau(x)=H(x,x)$ the Robin function.
The Dirichlet projection $PU_{\delta,\xi}\in H^1_0(\Om)$
is defined by
\[
-\Delta(PU_{\delta,\xi})=U_{\delta,\xi}^3
\quad\text{in }\Om,
\qquad
PU_{\delta,\xi}=0
\quad\text{on }\partial\Om.
\]
The dilation and translation fields associated with $U$ are
\begin{equation}\label{eq:Jacobi}
Z_0(y)=-U(y)-y\cdot\nabla U(y),
\qquad
Z_k(y)=-\partial_kU(y),
\quad 1\le k\le4.
\end{equation}

Fix an integer $m\ge1$.
For a configuration
$\boldsymbol\xi=(\xi_1,\ldots,\xi_m)$
of pairwise distinct points in $\Om$, define the
Green--Robin interaction matrix
$M(\boldsymbol\xi)
=(M_{ij}(\boldsymbol\xi))_{1\le i,j\le m}$ by
\begin{equation}\label{eq:M}
M_{ii}(\boldsymbol\xi)=\tau(\xi_i),
\qquad
M_{ij}(\boldsymbol\xi)=-G(\xi_i,\xi_j)
\quad(i\ne j).
\end{equation}
We denote its least eigenvalue by $\Lambda(\boldsymbol\xi)$
and the associated positive unit eigenvector by
$e(\boldsymbol\xi)$, so that
\[
M(\boldsymbol\xi)e(\boldsymbol\xi)
=\Lambda(\boldsymbol\xi)e(\boldsymbol\xi),
\qquad
|e(\boldsymbol\xi)|=1,
\qquad
e_i(\boldsymbol\xi)>0
\quad(1\le i\le m).
\]

We consider a family $\{u_\epsilon\}$ of positive solutions
to the four-dimensional Brezis--Nirenberg problem
\eqref{intro:BN}, concentrating at $m$ distinct interior
points as $\epsilon\to0$.
Let $\boldsymbol\xi^0=(\xi_1^0,\ldots,\xi_m^0)$
denote the limiting concentration configuration.
For each sufficiently small $\epsilon>0$, let
$\delta_{i,\epsilon}>0$ and $\xi_{i,\epsilon}\in\Om$
be the bubble widths and centers, respectively, and set
\[
U_i:=U_{\delta_{i,\epsilon},\xi_{i,\epsilon}},
\qquad
P_i:=PU_i,
\qquad
\rho_\epsilon
:=\left(\sum_{i=1}^m\delta_{i,\epsilon}^2\right)^{1/2}.
\]
The projected dilation and translation fields span the space
\begin{equation}\label{eq:K}
\mathcal K_\epsilon
:=\operatorname{span}
\left\{
\delta_{i,\epsilon}\partial_{\delta_i}P_i,\,
\delta_{i,\epsilon}\partial_{\xi_{i,k}}P_i:
1\le i\le m,\ 1\le k\le4
\right\}.
\end{equation}
We assume that the concentration parameters can be chosen
so that $u_\epsilon$ admits the orthogonal decomposition
\begin{equation}\label{eq:branch}
\begin{gathered}
u_\epsilon=\sum_{i=1}^mP_i+\phi_\epsilon,
\qquad
\phi_\epsilon\perp\mathcal K_\epsilon,\\
\|\phi_\epsilon\|_{H^1_0(\Om)}
\le C\gamma_\epsilon,
\qquad
\gamma_\epsilon:=\epsilon\rho_\epsilon+\rho_\epsilon^2,
\end{gathered}
\end{equation}
where orthogonality is taken with respect to the Dirichlet
inner product and $C$ is independent of $\epsilon$.
The widths satisfy
\begin{equation}\label{eq:widths}
\begin{gathered}
\rho_\epsilon\to0,
\qquad
\delta_{i,\epsilon}=\rho_\epsilon t_{i,\epsilon},
\qquad
t_\epsilon=(t_{1,\epsilon},\ldots,t_{m,\epsilon}),\\
|t_\epsilon|=1,
\qquad
t_\epsilon\to t=e(\boldsymbol\xi^0)>0,
\end{gathered}
\end{equation}
and the centers and common scale satisfy
\begin{equation}\label{eq:balance}
\begin{gathered}
\xi_{i,\epsilon}\to\xi_i^0\in\Om,
\qquad
\xi_i^0\ne\xi_j^0\quad(i\ne j),\quad
\Lambda:=\Lambda(\boldsymbol\xi^0)>0,
\qquad
\epsilon\log\rho_\epsilon^{-1}\to8\pi^2\Lambda.
\end{gathered}
\end{equation}

Conditions \eqref{eq:branch}--\eqref{eq:balance}
are satisfied by the modulated solution families constructed
in \cite[Proposition~2.1, (2.7), (2.18),
(3.3)--(3.4)]{P-R-V}, after fixing a limiting configuration
and adopting the above normalization of the bubbles and widths.
We take the existence of these concentrating families as an
input and study their weighted spectral asymptotics.
No nondegeneracy assumption on
$D^2\Lambda(\boldsymbol\xi^0)$ is imposed at this stage.

For each such solution $u_\epsilon$, we consider the weighted
eigenvalue problem
\begin{equation}\label{eq:spectral}
-\Delta w_{\ell,\epsilon}
=\mu_{\ell,\epsilon}(3u_\epsilon^2+\epsilon)w_{\ell,\epsilon}
\quad\text{in }\Omega,
\qquad
w_{\ell,\epsilon}\in H^1_0(\Omega).
\end{equation}
Its eigenvalues are arranged in nondecreasing order and counted
with multiplicity. Define the weighted inner product
\begin{equation}\label{eq:Dform}
D_\epsilon(f,h)
:=\int_\Omega(3u_\epsilon^2+\epsilon)fh\,dx.
\end{equation}
We choose a real eigenbasis orthonormal with respect to
$D_\epsilon$ and introduce the corresponding
$L^\infty$-normalized eigenfunctions:
\begin{equation}\label{tp:Dnormal}
D_\epsilon(w_{\ell,\epsilon},w_{j,\epsilon})
=\delta_{\ell j},
\qquad
v_{\ell,\epsilon}
:=\frac{w_{\ell,\epsilon}}
{\|w_{\ell,\epsilon}\|_{L^\infty(\Omega)}}.
\end{equation}

Our first result describes the eigenvalue expansions about
$1/3$ and the corresponding eigenfunction profiles for
$1\le \ell\le m$.
Write $M=M(\boldsymbol\xi^0)$ and set
\begin{equation}\label{eq:N}
D_t:=\operatorname{diag}(t_1,\ldots,t_m),
\qquad
\mathcal N:=M-\Lambda\mathrm I_m,
\qquad
\mathcal S:=D_t\mathcal N D_t.
\end{equation}
As shown in Lemma~\ref{lem:matrix}, both $\mathcal N$ and
$\mathcal S$ are positive semidefinite, with
\[
\ker\mathcal N=\operatorname{span}\{t\},
\qquad
\ker\mathcal S=\operatorname{span}\{\mathbf1\},
\qquad
\mathbf1=(1,\ldots,1)^T\in\mathbb R^m.
\]
Define
\begin{equation}\label{eq:amplitude}
\mathcal A_{1}
:=D_t\left(M-\frac{\Lambda}{3}\mathrm I_m\right)D_t
=\mathcal S+\frac{2\Lambda}{3}D_t^2.
\end{equation}
Let $\sigma_1\le\cdots\le\sigma_m$ denote its eigenvalues,
counted with multiplicity. Our first result gives the
asymptotic behavior of the first $m$ eigenvalues of
\eqref{eq:spectral} and the corresponding eigenfunction profiles.
\begin{thm}\label{thm:amp}
Under assumptions \eqref{eq:branch}--\eqref{eq:balance},
the matrix $\mathcal A_{1}$ is positive definite.
Moreover, for each $1\le\ell\le m$, the eigenvalue
$\mu_{\ell,\epsilon}$ satisfies
\begin{equation}\label{am:shift}
\mu_{\ell,\epsilon}=\frac13+4\pi^2\sigma_\ell\rho_\epsilon^2
+o(\rho_\epsilon^2).
\end{equation}
The eigenfunction profiles satisfy
\begin{equation}\label{am:corew}
\delta_{i,\epsilon}w_{\ell,\epsilon}(\xi_{i,\epsilon}+\delta_{i,\epsilon}y)
\rightharpoonup b_{\ell,i}U(y)
\quad\hbox{in }\cD(\R^4),
\end{equation}
and strongly in $C^1_{\rm loc}(\R^4)$, with zero extension
for the weak convergence. Their coefficient vectors satisfy
\begin{equation}\label{am:bmatrix}
\mathcal{A}_1 b_\ell=\sigma_\ell b_\ell,\qquad
g_0 \,b_\ell\cdot b_j=\delta_{\ell j}\quad\mbox{where}\quad g_{0}
:=3\int_{\mathbb R^4}U^4\,dy
 =32\pi^2.
\end{equation}
Moreover, there exists a constant $c_3$ defined by
\begin{equation}\label{cc3}
c_3=\int_{\mathbb R^4}U^3\,dy
 =8\sqrt2\,\pi^2
\end{equation}
such that
\begin{equation}\label{am:outerw}
\rho_\epsilon^{-1}w_{\ell,\epsilon}\longrightarrow
c_3\sum_i t_i b_{\ell,i}G(\cdot,\xi_i^0)\quad\text{in }C^1_{\mathrm{loc}}
\bigl(\Omega\setminus\{\xi_1^0,\ldots,\xi_m^0\}\bigr).
\end{equation}
\end{thm}

We next consider the range $m+1\le\ell\le5m$, where the
eigenvalues converge to $1$.
Our first result in this range identifies the limiting local
profiles and shows that they contain only the translation
fields $Z_1,\ldots,Z_4$.
For its limiting local
coefficients $a_{\ell,i,k}$ in the fields $Z_k$, we use the vectors
\begin{equation}\label{eq:d}
d_\ell=\big(t_1 a_{\ell,1,1},\cdots,t_1 a_{\ell,1,4},t_2 a_{\ell,2,1},\cdots,t_m a_{\ell,m,4}\big)\in\R^{4m}.
\end{equation}

\begin{thm}\label{thm:pure}
Under \eqref{eq:branch}--\eqref{eq:balance}, for
$m+1\le\ell\le5m$,
\begin{equation}\label{tp:mainprofile}
v_{\ell,\epsilon}(\xi_{i,\epsilon}+\delta_{i,\epsilon}y)
\rightharpoonup\sum_{k=1}^4a_{\ell,i,k}Z_k(y)
\quad\hbox{in}\hspace{3mm}\cD(\R^4),
\end{equation}
and strongly in $C^1_{\rm loc}(\R^4)$.
At least one vector $a_{\ell,i}=(a_{\ell,i,k})_{k=1}^4$ is nonzero.
The coefficient vectors $d_\ell$ defined in \eqref{eq:d}
are nonzero and mutually orthogonal.
\end{thm}

For the same range $m+1\le\ell\le5m$, we now describe the
eigenfunctions away from the concentration points.
The following theorem determines the exterior expansion of
$v_{\ell,\epsilon}$ at order $\rho_\epsilon^3$ in terms of
the local coefficients obtained in Theorem~\ref{thm:pure}.
Define the matrix
\begin{equation}\label{eq:Ndagger}
\cN^\dagger=(\cN+tt^T)^{-1}-tt^T.
\end{equation}
For every $f\in\mathbb R^m$ satisfying $t^Tf=0$,
the vector $z=\mathcal N^\dagger f$ is the unique solution of
\[
\mathcal N z=f,\qquad t^Tz=0.
\]
Moreover, $\mathcal N^\dagger t=0$.
We denote the entries of $\mathcal N^\dagger$ by
$\eta^{ij}$, so that
$\mathcal N^\dagger=(\eta^{ij})_{1\le i,j\le m}$.

Define the $m\times4m$ matrix $\cP$ and the row
vectors $\cG$, $\widetilde{\cG}$ by
\begin{equation}\label{eq:P}
\cP_{i,(j,k)}=
\begin{cases}
-\dfrac12t_i^2\partial_k\tau(\xi_i^0),&j=i,\\
t_j^2\partial_{y_k}G(\xi_i^0,\xi_j^0),&j\ne i,
\end{cases}
\end{equation}
\begin{equation}\label{eq:G}
\cG(x)=(G(x,\xi_i^0))_i,\qquad
\widetilde{\cG}(x)=(t_i^2\partial_{y_k}G(x,\xi_i^0))_{(i,k)}.
\end{equation}
\begin{thm}\label{thm:exterior}
Under assumptions \eqref{eq:branch}--\eqref{eq:balance},
for $m+1\le\ell\le5m$ and along the subsequences
in Theorem~\ref{thm:pure}, we have
\begin{equation}\label{ep:final}
\begin{aligned}
\rho_\epsilon^{-3}v_{\ell,\epsilon}(x)
\longrightarrow c_3\Bigg[
&\sum_{i=1}^{m}\sum_{j=1}^{m}\sum_{k=1}^{4}
\eta^{ij}
\Bigg(
-\frac12 t_j^3a_{\ell,j,k}
\frac{\partial\tau}{\partial x_k}(\xi_j^0)
\\
&\qquad
+\sum_{\substack{r=1\\r\ne j}}^{m}
t_r^3a_{\ell,r,k}
\frac{\partial G}{\partial y_k}(\xi_j^0,\xi_r^0)
\Bigg)G(x,\xi_i^0)
\\
&+\sum_{i=1}^{m}\sum_{k=1}^{4}
t_i^3a_{\ell,i,k}
\frac{\partial G}{\partial y_k}(x,\xi_i^0)
\Bigg]
\end{aligned}
\end{equation}
in
$C^1_{\mathrm{loc}}
\bigl(\Omega\setminus\{\xi_1^0,\ldots,\xi_m^0\}\bigr)$
as $\epsilon\to0$.
Moreover, the vector of coefficients multiplying
$G(x,\xi_i^0)$ is orthogonal to $t$, namely,
\(t^T\mathcal N^\dagger\mathcal P d_\ell=0.
\)
\end{thm}

Set
\begin{equation}\label{eq:Q4}
(\cQ_4)_{(i,k),(j,h)}=
\begin{cases}
t_i^2t_j^2\partial_{x_ky_h}G(\xi_i^0,\xi_j^0),&i\ne j,\\[1mm]
-\dfrac{t_i^4}{2}\partial_{kh}\tau(\xi_i^0)
+t_i^3\displaystyle\sum_{r\ne i}t_r
\partial_{x_kx_h}G(\xi_i^0,\xi_r^0),&i=j,
\end{cases}
\end{equation}
\begin{equation}\label{eq:B4}
\cB_4=\cP^T\cN^\dagger\cP+\cQ_4,\qquad
\mathsf D=\diag(t_1\Id_4,\ldots,t_m\Id_4).
\end{equation}
The matrix $\cB_4$ is symmetric. Denote its eigenvalues in decreasing
order by $\nu_{m+1}\ge\cdots\ge\nu_{5m}$.

The exterior expansion allows us to determine the leading
corrections to $\mu_{\ell,\epsilon}$ for $m+1\le\ell\le5m$.
The following theorem identifies the matrix governing these
corrections and relates it to the Hessian of $\Lambda$ at
the limiting concentration configuration.
\begin{thm}\label{thm:translation}
Under \eqref{eq:branch}--\eqref{eq:balance},
\begin{equation}\label{tr:mainshift}
\mu_{\ell,\epsilon}=1-20\pi^2\nu_\ell\rho_\epsilon^4+o(\rho_\epsilon^4),
\qquad m+1\le\ell\le5m.
\end{equation}
The vectors in \eqref{eq:d} satisfy
$\cB_4d_\ell=\nu_\ell d_\ell$ and form an orthogonal basis of
$\R^{4m}$. Furthermore,
\begin{equation}\label{tr:Hessianidentity}
\nabla\Lambda(\bm\xi^0)=0,\qquad
\cB_4=-\frac12\mathsf D\,D^2\Lambda(\bm\xi^0)\,\mathsf D.
\end{equation}
\end{thm}

It remains to consider the range $5m+1\le\ell\le6m$,
corresponding to the dilation modes.
The following theorem distinguishes the common dilation
mode, with $\ell=5m+1$, from the relative dilation modes,
with $5m+2\le\ell\le6m$, and gives their corrections at
orders $\epsilon\rho_\epsilon^2$ and $\rho_\epsilon^2$,
respectively. When $m=1$, only the common dilation mode
is present.
\begin{thm}\label{thm:dilation}
Let
$0<\kappa_1\le\cdots\le\kappa_{m-1}$ be the positive eigenvalues
of $\cS$ in \eqref{eq:N}, counted with multiplicity. Under \eqref{eq:branch}--\eqref{eq:balance},
\begin{align}
\mu_{5m+1,\epsilon}
&=1+\frac5{2m}\epsilon\rho_\epsilon^2+o(\epsilon\rho_\epsilon^2),
\label{ds:commonformula}\\
\mu_{5m+1+j,\epsilon}
&=1+20\pi^2\kappa_j\rho_\epsilon^2+o(\rho_\epsilon^2),
\qquad 1\le j\le m-1.\label{ds:relativeformula}
\end{align}
All these eigenvalues exceed $1$ for small $\epsilon$.
Their rescaled $D_\epsilon$-normalized profiles are
\begin{equation}\label{ds:Dprofile}
\delta_{i,\epsilon}w_{\ell,\epsilon}(\xi_{i,\epsilon}+\delta_{i,\epsilon}y)
\rightharpoonup b_{\ell,i}Z_0(y)
\quad\hbox{in }\cD(\R^4)
\quad\hbox{and strongly in }C^1_{\rm loc}(\R^4),
\end{equation}
where
\begin{equation}\label{ds:relativevectors}
b_{5m+1}=\pm\frac{\one}{\sqrt{mg_1}},\qquad
\cS b_{5m+1+j}=\kappa_jb_{5m+1+j},\qquad
g_1\,b_\ell\cdot b_{\ell'}=\delta_{\ell\ell'}
\end{equation}
with
$$g_1:=3\int_{\mathbb R^4}U^2Z_\beta^2\,dy
 =\frac{32\pi^2}{5},\quad\mbox{for}\quad0\le\beta\le4.$$
\end{thm}

For a real symmetric matrix $\mathcal{A}$, let $n_-(\mathcal{A})$, $n_+(\mathcal{A})$ and
$n_0(A)$ denote its negative index, positive index and nullity.
The Morse index and nullity of $u_\epsilon$ are those of the
Dirichlet operator $-\Delta-3u_\epsilon^2-\epsilon$, equivalently,
$\ind(u_\epsilon)=\#\{\ell:\mu_{\ell,\epsilon}<1\}$ and
$\nul(u_\epsilon)=\#\{\ell:\mu_{\ell,\epsilon}=1\}$.
The preceding expansions give the following consequence.
\begin{cor}\label{cor:morse}
Let $\mathcal{H}=D^2\Lambda(\bm\xi^0)$. Under the preceding
hypotheses, for sufficiently small $\epsilon$,
\begin{equation}\label{eq:morse}
m+n_-(\mathcal{H})\le\ind(u_\epsilon)
\le\ind(u_\epsilon)+\nul(u_\epsilon)
\le m+n_-(\mathcal{H})+n_0(\mathcal{H}).
\end{equation}
In particular, $\nul(u_\epsilon)\le n_0(\mathcal{H})$.
If $\mathcal{H}$ is nondegenerate, then
$\ind(u_\epsilon)=m+n_-(\mathcal{H})$ and $\nul(u_\epsilon)=0$.
\end{cor}

\subsection{Main difficulties and proof strategy}\label{subsec:strategy}
The proof follows the scheme of Choi, Kim and Lee \cite{CKL-1}
and Yang and Zhao \cite{YZmulti}: identify the local spectral
profiles, obtain their exterior expansions through Green's
formula, and use local bilinear identities to determine the
finite-dimensional eigenvalue problems. In \cite{YZmulti}, although their spectral convergence theorem \cite[Theorem 1.2]{YZmulti} covers \(N=4\), \(q=1\), their refined multi-bubble analysis requires \(q>1\) in dimension four. The endpoint \(q=1\) presents additional difficulties: the concentration scales are exponentially small, logarithmic volume moments contribute at leading order, the decay iteration loses its strict exponent gain, and the leading scale matrix has a kernel corresponding to a common dilation of all bubbles; see \cite[ Introduction and Table 1]{YZmulti}. These features prevent a direct extension of the refined analysis to this endpoint.  We explain how the lemmas below address these difficulties and fit into the proof.

\medskip
\noindent\textbf{\textbf{Decay estimates and the first $m$ eigenpairs.}}
The comparison iteration in \cite[Appendix~A]{YZmulti}
requires a strict gain in the decay exponent, which is lost at
the present endpoint. Lemma~\ref{tp:endpoint} replaces this
step by a Green-function argument on the physical domain with
large bubble cores removed. The critical potential has a small
weighted Green norm on this region, while the constant linear
potential contributes $O(\epsilon)$. Absorbing these terms gives
the optimal decay with exponent $2$, in the role played by
\cite[Lemma~2.8]{CKL-1}. This estimate will be used throughout the paper to control the contributions to the local integrals from outside
large balls in the rescaled variables.
A second difference already appears in the amplitude cluster.
The bubble is not square integrable in $N=4$, and
$\epsilon\log\rho^{-1}$ has a nonzero limit. Thus the volume integral term containing $\epsilon\int u_\epsilon w_{\ell,\epsilon}$ cannot be
discarded from the local identity. Lemma~\ref{am:logmoment}
computes its logarithmic coefficient by matching the core and
intermediate-annulus expansions. Combined with the boundary integral
calculation in Lemma~\ref{am:flux}, this gives the matrix
$\mathcal{A}_1$. The profile and orthogonality argument in
Lemma~\ref{am:profiles} then identifies all $m$ amplitude
eigenpairs, including multiplicities, and proves
Theorem~\ref{thm:amp}.

\medskip
\noindent\textbf{Separation of translations and dilations.}
For eigenvalues converging to $1$, the local limits may
contain both $Z_0$ and $Z_1,\ldots,Z_4$.
Proposition~\ref{tp:upper}, using the test functions
of \cite[Proposition~4.1]{CKL-1}, gives
$\mu_{\ell,\epsilon}\le1+C\rho_\epsilon^4$
for $m+1\le\ell\le5m$.
Lemma~\ref{tp:moments} computes the logarithmic term
in the local scaling identity, and
Proposition~\ref{tp:leadingconstraint} gives the
equation for the dilation coefficients.
The corresponding matrix in \cite{YZmulti} is
positive definite, whereas here
$\ker\cS=\Span\{\one\}$.
Thus, for these indices, the equation only shows
that the coefficients of $Z_0$ are equal at all
bubbles. It does not show that they are zero,
even when $m=1$.
To prove that this remaining coefficient is zero,
we keep the centers fixed and vary the widths.
Set $s=\log\rho^{-1}$, let $h$ be local coordinates
for the width ratios, and denote the reduced energy
by $F_\epsilon(s,h)$.
In addition to the existence and $H^1_0$ estimates
in \cite{P-R-V}, we need estimates for derivatives
with respect to $s$ and $h$.
Lemma~\ref{tp:scaledifferentiation} provides these
estimates for the projected bubbles, while
Lemma~\ref{tp:scaleLS} gives the first and second
derivatives of $\phi_\epsilon$.
Proposition~\ref{tp:energy} then gives an expansion
of $F_\epsilon$ that includes the terms of order
$\epsilon\rho^2$ and controls the first two
derivatives of the error.
At a solution, $F_s=F_h=0$.
Using these identities, Lemma~\ref{tp:curvature} gives
\[
F_{ss}=16\pi^2\epsilon\rho^2+o(\epsilon\rho^2),
\qquad
F_{sh}=o(\epsilon\rho^2),
\qquad
F_{hh}\ge c\rho^2 I.
\]
The variables $h$ are absent when $m=1$.
Lemma~\ref{tp:translations} constructs translation
functions that are orthogonal to $\mathcal K_\epsilon^\perp$
with respect to the bilinear form \eqref{eq:forms}.
Together with the estimates above, the decomposition
in Lemma~\ref{tp:decomposition} and the eigenvalue
bounds in Lemma~\ref{tp:twosided} give the coefficient
estimates in Lemma~\ref{tp:coefficients}.
All coefficients of $Z_0$ then tend to zero for
$m+1\le\ell\le5m$.
This proves Theorem~\ref{thm:pure} without assuming
that $D^2\Lambda(\bm\xi^0)$ is invertible.

\medskip
\noindent\textbf{Expansion away from the concentration points.}
For $m+1\le\ell\le5m$, the local limits obtained above
do not yet give the expansion away from the centers
at order $\rho^3$.
The decay estimate only gives a bound of order $\rho^2$.
Taking absolute values in the linear term of the Green
representation gives a bound of order $\epsilon\rho^2$,
which is too large for this expansion.
Following \cite[Section~7]{YZmulti}, we keep
$q=(q_i)_{i=1}^m$ in the error estimates, where
\[
q_i
=3\int_{B(\xi_{i,\epsilon},r)}
u_\epsilon^2v_{\ell,\epsilon}\,dx
\]
and $r>0$ is fixed and small.
Lemma~\ref{ep:weightedlemma} estimates the difference
between the eigenfunction and its local translation
terms.
Integration by parts in Lemma~\ref{ep:linearpotential}
then improves the estimate for the linear term.
These estimates give the Green function expansion
in Lemma~\ref{ep:multipole}, with an error that
still depends on $q$.
Since $\ker\cN=\Span\{t\}$, the equations for the
limiting coefficients do not determine the component
along $t$.
We therefore need a separate estimate for $t^Tq$.
The estimate for $F_{ss}$ alone is not enough here.
Lemma~\ref{ep:commonlemma} proves
\(
(\phi_s+\phi)/(\epsilon\rho)\longrightarrow0
\quad\hbox{in }H^1_0(\Om),
\)
where $\phi_s=\partial_s\phi$.
For the bilinear form \eqref{eq:forms}, this gives
a bound of $o(\epsilon\rho^3)$ for the pairing of
$\partial_s(\sum_iP_i+\phi)$ with the translation
functions in Lemma~\ref{tp:translations}.
Using this bound, the integral estimates in
Lemma~\ref{ep:sourcemoments} and the coefficient
estimates above, Lemma~\ref{ep:massclosure} shows that
\[
q=O(\rho^3),
\qquad
t^Tq=o(\rho^3).
\]
The error in the Green function expansion is then
$o(\rho^3)$.
Lemmas~\ref{ep:monopole} and \ref{ep:dipolelemma}
determine the coefficients of $G$ and its derivatives,
respectively, and prove Theorem~\ref{thm:exterior}.
The condition $t^Tq=o(\rho^3)$ places the limit of
$\rho^{-3}q$ in $t^\perp$, where $\cN$ is invertible.
This explains the use of $\cN^\dagger$ in the theorem.

\medskip
\noindent\textbf{Eigenvalue expansions and Morse index.}
Using the limits away from the centers, we apply
Green's identity as in \cite{CKL-1,YZmulti}.
Lemma~\ref{tr:flux} computes the boundary integral
at order $\rho^4$, and Lemma~\ref{tr:pairing} gives
the limit of the volume integral.
These results yield the matrix equation for $\cB_4$.
The coefficient vectors form an orthogonal basis,
so all eigenvalues of $\cB_4$ occur, counted with
multiplicity.
Their ordering gives the eigenvalue expansions in
Theorem~\ref{thm:translation}.
Differentiating $\Lambda$ gives the relation between
$\cB_4$ and $D^2\Lambda(\bm\xi^0)$ in the same theorem.
For $5m+1\le\ell\le6m$, Lemma~\ref{ds:pure} shows
that the local limits contain only $Z_0$.
Lemma~\ref{ds:leading} determines the terms of order
$\rho^2$ for $5m+2\le\ell\le6m$.
To study $\mu_{5m+1,\epsilon}$,
Lemma~\ref{ds:common} tests the eigenvalue equation
with $\partial_s(\sum_iP_i+\phi)$ and uses the
estimate for $F_{ss}$.
This gives a positive leading term of order
$\epsilon\rho^2$ in $\mu_{5m+1,\epsilon}-1$.
Here the bilinear estimate $O(\rho^3)$ is enough,
since $\rho^3=o(\epsilon\rho^2)$.
The stronger estimate in Lemma~\ref{ep:commonlemma}
is needed for the expansion away from the centers.
This proves Theorem~\ref{thm:dilation}.
The first $m$ eigenvalues are below $1$, while those
with $5m+1\le\ell\le6m$ are above $1$.
Since $\mu_{6m+1,\epsilon}\to2$, all larger eigenvalues
are also above $1$ for small $\epsilon$.
The expansions for $m+1\le\ell\le5m$ and the relation
with $D^2\Lambda(\bm\xi^0)$ then give the bounds in
Corollary~\ref{cor:morse}.
When this Hessian is invertible, they give the exact
Morse index and show that $u_\epsilon$ is nondegenerate.

\medskip
\noindent
The rest of the paper is organized as follows:
Section~\ref{sec:prelim} introduces the notation
and gives the limiting spectrum and decay estimates.
Section~\ref{sec:amplitude} proves
Theorem~\ref{thm:amp}.
Sections~\ref{sec:pure} and \ref{sec:exterior}
prove Theorems~\ref{thm:pure} and \ref{thm:exterior}.
Sections~\ref{sec:translation} and \ref{sec:dilation}
give the translation and dilation eigenvalue
expansions and prove Corollary~\ref{cor:morse}.

\section{Notation and preliminaries}\label{sec:prelim}
In this section, we introduce some notations and important results that will be used throughout this paper.
We suppress $\epsilon$ in the widths and centers when no confusion
is possible and write
\[
\rho=\rho_\epsilon,\quad s=\log\rho^{-1},\quad
\gamma=\epsilon\rho+\rho^2,\quad V_\epsilon=\sum_iP_i,\quad
X_\epsilon=\cK_\epsilon^\perp.
\]
Norms and orthogonality on $H^1_0(\Om)$ refer to the Dirichlet
inner product. We also use
$i^*=(-\Delta_D)^{-1}:L^{4/3}(\Om)\to H^1_0(\Om)$ and the
linearized bilinear form
\begin{equation}\label{eq:forms}
Q_\epsilon(f,h)=\int_\Om\nabla f\cdot\nabla h-D_\epsilon(f,h).
\end{equation}
We define the rescaled Jacobi fields associated with
$Z_\beta$ in \eqref{eq:Jacobi} and their Dirichlet projections by
\begin{equation}\label{eq:scaled-Jacobi}
\begin{gathered}
J_{i\beta}(x)=\delta_i^{-1}Z_\beta((x-\xi_i)/\delta_i),
\qquad 0\le\beta\le4,\\
-\Delta PJ_{i\beta}=3U_i^2J_{i\beta}\ \text{in }\Om,
\qquad PJ_{i\beta}=0\ \text{on }\partial\Om.
\end{gathered}
\end{equation}
Then $PJ_{i0}=\delta_i\partial_{\delta_i}P_i$ and
$PJ_{ik}=\delta_i\partial_{\xi_{i,k}}P_i$, so that
$\cK_\epsilon=\Span\{PJ_{i\beta}:1\le i\le m,\ 0\le\beta\le4\}$.
We abbreviate
\[
\alpha=2\sqrt2,\qquad A=\int_{\R^4}U^4=\frac{32\pi^2}{3},
\qquad d(x)=\min_i|x-\xi_i|.
\]
Using \eqref{eq:balance}, we obtain
\begin{equation}\label{eq:exponential}
\rho\epsilon^{-j}(1+s)^K\longrightarrow0
\quad\hbox{for every fixed }j,K,\qquad s\approx\epsilon^{-1}.
\end{equation}
Here we say that $a\approx b$ if $a\geq Cb$ and $a\leq Cb$, and $c$ and $C$ are indiscriminately used to denote various positive constants in the sequel.

\begin{lem}\label{lem:matrix}
The matrices in \eqref{eq:N} satisfy
\begin{align}
z^T\cN z&=\sum_{i<j}G(\xi_i^0,\xi_j^0)t_it_j
\left(\frac{z_i}{t_i}-\frac{z_j}{t_j}\right)^2,\label{eq:graphN}\\
b^T\cS b&=\sum_{i<j}G(\xi_i^0,\xi_j^0)t_it_j(b_i-b_j)^2.
\label{eq:graphS}
\end{align}
Thus $\mathcal N$ is positive semidefinite, and the nullspaces
of $\mathcal N$ and $\mathcal S$ are spanned by $t$ and
$\mathbf1$, respectively. Moreover,
\[
\mathcal N^\dagger t=0,
\qquad
\mathcal N\mathcal N^\dagger
=\mathcal N^\dagger\mathcal N=I_m-tt^T.
\]
Finally, $\mathcal A_1$ is positive definite.
\end{lem}
\begin{proof}
Since $Mt=\Lambda t$ and $t_i>0$, we have
\[
\tau(\xi_i^0)-\Lambda
=\sum_{j\ne i}G(\xi_i^0,\xi_j^0)\frac{t_j}{t_i}.
\]
For any $z\in\mathbb R^m$, expanding the quadratic form
as in \cite[Lemma~7.4]{YZmulti} gives
\begin{align*}
z^T\mathcal N z
&=\sum_{i=1}^m
  \left(\sum_{j\ne i}
  G(\xi_i^0,\xi_j^0)\frac{t_j}{t_i}\right)z_i^2
  -2\sum_{i<j}G(\xi_i^0,\xi_j^0)z_i z_j\\
&=\sum_{i<j}G(\xi_i^0,\xi_j^0)t_it_j
  \left(\frac{z_i}{t_i}-\frac{z_j}{t_j}\right)^2.
\end{align*}
This proves \eqref{eq:graphN}.
Taking $z=D_tb$ and using $\mathcal S=D_t\mathcal N D_t$
yields \eqref{eq:graphS}.

Since $\Om$ is connected, $G(\xi_i^0,\xi_j^0)>0$ for $i\ne j$.
Thus $\mathcal N$ and $\mathcal S$ are positive semidefinite.
For $m\ge2$, equality in \eqref{eq:graphN} holds precisely
when all ratios $z_i/t_i$ are equal, whereas equality in
\eqref{eq:graphS} holds precisely when all components of $b$
are equal. Hence
\[
\ker\mathcal N=\operatorname{span}\{t\},
\qquad
\ker\mathcal S=\operatorname{span}\{\mathbf1\}.
\]
These conclusions also hold for $m=1$, since
$\mathcal N=\mathcal S=0$ in that case.
As $|t|=1$, the matrix $\mathcal N+tt^T$ is positive definite
and satisfies $(\mathcal N+tt^T)t=t$.
Consequently, its inverse maps $t$ to itself, and
the definition \eqref{eq:Ndagger} gives
\[
\mathcal N^\dagger t=0,
\qquad
\mathcal N\mathcal N^\dagger
=\mathcal N^\dagger\mathcal N=I_m-tt^T.
\]
In particular, for every $f\perp t$,
$\mathcal N^\dagger f$ is the unique vector $z\perp t$
satisfying $\mathcal N z=f$.

Finally, by \eqref{eq:amplitude}, for every
$b\in\mathbb R^m\setminus\{0\}$,
\[
b^T\mathcal A_1b
=b^T\mathcal Sb+\frac{2\Lambda}{3}
  \sum_{i=1}^m t_i^2b_i^2
\ge\frac{2\Lambda}{3}\sum_{i=1}^m t_i^2b_i^2>0.
\]
Therefore $\mathcal A_1$ is positive definite.
\end{proof}

\begin{lem}\label{tp:profilelemma}
As $\epsilon\to0$, for each $1\le i\le m$,
\[
\delta_i u_\epsilon(\xi_i+\delta_i y)
\longrightarrow U(y)
\quad\text{in }C^1_{\mathrm{loc}}(\mathbb R^4).
\]
Moreover,
\begin{equation}\label{tp:outeru}
\begin{gathered}
\rho^{-1}u_\epsilon
\longrightarrow
c_3\sum_{i=1}^m t_iG(\cdot,\xi_i^0)
\quad
\text{in }C^2_{\mathrm{loc}}
\bigl(\Om\setminus\{\xi_1^0,\ldots,\xi_m^0\}\bigr).
\end{gathered}
\end{equation}
There exist constants $C>0$ and $\epsilon_0>0$ such that
\begin{equation}\label{tp:udecay}
u_\epsilon(x)
\le C\frac{\rho}{\rho^2+d(x)^2}
\qquad
\text{for all }x\in\Om,\quad 0<\epsilon<\epsilon_0.
\end{equation}
\end{lem}
\begin{proof}
By the projection estimates in
\cite[(2.4)--(2.5)]{P-R-V},
\[
P_i
=U_i-c_3\delta_iH(\cdot,\xi_i)+O(\delta_i^3).
\]
Together with \eqref{eq:branch} and the Sobolev inequality,
this gives
\[
\delta_i u_\epsilon(\xi_i+\delta_i y)
\longrightarrow U(y)
\quad\text{in }L^4_{\mathrm{loc}}(\mathbb R^4).
\]
Indeed, the rescaled contribution of $\phi_\epsilon$
tends to zero in $L^4_{\mathrm{loc}}$, while those of the
other bubbles tend to zero locally uniformly.
Since $U$ is locally bounded, this convergence ensures
uniform smallness of the local $L^4$ norms on sufficiently
small balls. The local iteration argument in
\cite[Lemmas~6--7]{HANZCHAO} therefore gives uniform bounds
on compact sets. This argument also applies to the
rescaled equation
\[
-\Delta z_{i,\epsilon}
=z_{i,\epsilon}^3+\epsilon\delta_i^2z_{i,\epsilon},
\qquad
z_{i,\epsilon}(y)
:=\delta_i u_\epsilon(\xi_i+\delta_i y),
\]
since the additional zeroth-order coefficient
$\epsilon\delta_i^2$ is uniformly bounded.
Interior elliptic estimates
\cite[Chapter~9]{GT} and the preceding
$L^4_{\mathrm{loc}}$ convergence then yield
\[
\delta_i u_\epsilon(\xi_i+\delta_i y)
\longrightarrow U(y)
\quad\text{in }C^1_{\mathrm{loc}}(\mathbb R^4).
\]

On compact subsets of
$\Om\setminus\{\xi_1^0,\ldots,\xi_m^0\}$, the projection
estimates also give
\[
P_i=c_3\delta_iG(\cdot,\xi_i)+O(\delta_i^3).
\]
Since $\delta_i/\rho\to t_i$, $\xi_i\to\xi_i^0$, and
\[
\left\|\frac{\phi_\epsilon}{\rho}\right\|_{L^4(\Om)}
\le C\frac{\gamma}{\rho}
=C(\epsilon+\rho)\longrightarrow0,
\]
we obtain
\[
\frac{u_\epsilon}{\rho}
\longrightarrow
c_3\sum_{i=1}^m t_iG(\cdot,\xi_i^0)
\quad\text{in }L^4_{\mathrm{loc}}
\bigl(\Om\setminus\{\xi_1^0,\ldots,\xi_m^0\}\bigr).
\]
In particular, the local $L^4$ norm of $u_\epsilon$
away from the concentration points is $O(\rho)$.
Applying the same local iteration argument to $u_\epsilon$
gives a local $L^\infty$ bound of order $\rho$.
The interior $W^{2,p}$ estimates
\cite[Chapter~9]{GT} consequently imply
\[
\|u_\epsilon\|_{C^{1,\sigma}(K)}
\le C_K\rho
\]
for every
$K\Subset\Om\setminus\{\xi_1^0,\ldots,\xi_m^0\}$
and some $\sigma\in(0,1)$.
Thus $v_\epsilon:=u_\epsilon/\rho$ is locally bounded
in $C^{1,\sigma}$ and satisfies
\[
-\Delta v_\epsilon
=\rho^2v_\epsilon^3+\epsilon v_\epsilon.
\]
The right-hand side tends to zero locally in
$C^{0,\sigma}$. Interior Schauder estimates
\cite[Chapter~6]{GT}, together with the
$L^4_{\mathrm{loc}}$ limit above, give the
$C^2_{\mathrm{loc}}$ convergence in \eqref{tp:outeru}.
It remains to prove \eqref{tp:udecay}.
Fix a sufficiently large constant $R>0$.
If $d(x)\le R\rho$, then
$|x-\xi_i|\le R\rho$ for some $i$.
Since $\delta_i/\rho\to t_i>0$, the rescaled points
$(x-\xi_i)/\delta_i$ remain in a fixed ball.
The local convergence proved above therefore yields
\(u_\epsilon(x)\le C_R/\rho,
\)
which implies \eqref{tp:udecay} in this region.

Next, suppose that $d:=d(x)\ge R\rho$ and that $x$
lies outside a fixed neighborhood of $\partial\Om$.
Choose a fixed $\theta\in(0,1/4)$ sufficiently small
that $B(x,2\theta d)\subset\Om$, and set $r=\theta d$.
For $y\in B(x,2r)$ and every $j$, we have
\[
|y-\xi_j|\ge(1-2\theta)d.
\]
Using $0\le P_j\le U_j$ and
$\|\phi_\epsilon\|_{L^4(\Om)}\le C\gamma$, we obtain
\[
\|u_\epsilon\|_{L^4(B(x,2r))}
\le C\left(\frac{\rho}{d}+\gamma\right)
\le \frac{C}{R}+C\gamma.
\]
The rescaled function
\[
z(y):=r\,u_\epsilon(x+ry),\qquad y\in B_2,
\]
satisfies
\[
-\Delta z=z^3+\epsilon r^2z,
\qquad
\|z\|_{L^4(B_2)}
=\|u_\epsilon\|_{L^4(B(x,2r))}.
\]
For $R$ sufficiently large and $\epsilon$ sufficiently
small, its $L^4$ norm is uniformly small and
$0\le\epsilon r^2\le1$.
The local iteration argument in
\cite[Lemmas~6--7]{HANZCHAO}, with this bounded
zeroth-order term, gives
\[
z(0)\le C\|z\|_{L^4(B_2)}.
\]
Consequently,
\[
\begin{aligned}
u_\epsilon(x)
&\le \frac{C}{d}
       \|u_\epsilon\|_{L^4(B(x,2r))}
\le C\left(\frac{\rho}{d^2}
            +\frac{\gamma}{d}\right)
\le C\frac{\rho}{d^2}.
\end{aligned}
\]
Here the last inequality follows from
$\gamma/\rho=\epsilon+\rho$ and the boundedness of $\Om$.
Since $d\ge R\rho$, this proves \eqref{tp:udecay}
in the second region. Finally, choose fixed boundary neighborhoods
$\mathcal C\subset\mathcal C'$ whose closures are
separated from the concentration points, with
$\mathcal C$ compactly contained in $\mathcal C'$
relative to $\overline{\Om}$.
The projection estimates and
$\|\phi_\epsilon\|_{L^4(\Om)}\le C\gamma$ imply
$
\|u_\epsilon\|_{L^4(\mathcal C')}\le C\rho.
$
Since $u_\epsilon=0$ on $\partial\Om$, the same
local iteration on intersections of balls with $\Om$
applies up to the boundary; see
\cite[Lemmas~6--7]{HANZCHAO}.
A finite covering of the boundary therefore gives
$\|u_\epsilon\|_{L^\infty(\mathcal C)}\le C\rho.
$
As $d(x)$ is bounded above and away from zero
on $\mathcal C$, this also implies \eqref{tp:udecay}
there. Completing the proof.
\end{proof}

\subsection{The spectral estimates}\label{sec:coarse}

We recall the limiting eigenvalues and establish a
quantitative decomposition of the eigenfunctions
whose eigenvalues converge to one.
By \cite[Theorem~1.2]{YZmulti},
\begin{equation}\label{tp:coarse}
\mu_{1,\epsilon},\ldots,\mu_{m,\epsilon}\to\frac13,
\qquad
\mu_{m+1,\epsilon},\ldots,\mu_{6m,\epsilon}\to1,
\qquad
\mu_{6m+1,\epsilon}\to2.
\end{equation}
Under the normalization \eqref{tp:Dnormal},
\cite[Lemma~3.1 and Proposition~3.2]{YZmulti}
give, along a common subsequence,
\begin{equation}\label{tp:coregeneral}
\begin{gathered}
\delta_iw_{\ell,\epsilon}(\xi_i+\delta_i y)
\rightharpoonup
\sum_{\beta=0}^4 b_{\ell,i,\beta}Z_\beta
\quad\text{in }\mathcal D^{1,2}(\mathbb R^4),
\qquad 1\le i\le m,\\
g_1\sum_{i=1}^m\sum_{\beta=0}^4
b_{\ell,i,\beta}b_{j,i,\beta}
=\delta_{\ell j},
\qquad m+1\le\ell,j\le6m.
\end{gathered}
\end{equation}
By Lemma~\ref{tp:profilelemma} and interior elliptic
estimates \cite[Chapter~9]{GT}, the convergence also
holds in $C^1_{\mathrm{loc}}(\mathbb R^4)$.
\begin{lem}\label{tp:quantprojection}
Let $i^*$ denote the inverse of $-\Delta$ with zero
Dirichlet boundary data, and define
\[
\mathcal A_u v
:=v-i^*\bigl((3u_\epsilon^2+\epsilon)v\bigr),
\qquad v\in H^1_0(\Om).
\]
Then
\begin{equation}\label{tp:Jacres}
\|\mathcal A_u PJ_{i\beta}\|_{H^1_0(\Om)}
\le C\gamma,
\qquad 1\le i\le m,\quad 0\le\beta\le4.
\end{equation}
For $m+1\le\ell\le6m$, every eigenfunction
$w_\ell=w_{\ell,\epsilon}$ normalized by
\eqref{tp:Dnormal} admits the decomposition
\begin{equation}\label{tp:thresholdprojection}
\begin{gathered}
w_\ell
=\sum_{i=1}^m\sum_{\beta=0}^4
b_{\ell,i,\beta,\epsilon}PJ_{i\beta}
+\psi_\ell,\quad\mbox{with}\quad
|b_{\ell,\epsilon}|\le C,
\quad
\|\psi_\ell\|_{H^1_0(\Om)}\le C\gamma,
\end{gathered}
\end{equation}
where
$b_{\ell,\epsilon}
=(b_{\ell,i,\beta,\epsilon})_{1\le i\le m,\;0\le\beta\le4}$.
Along the subsequence in \eqref{tp:coregeneral},
\begin{equation}\label{tp:outerw}
\begin{gathered}
\rho^{-1}w_\ell
\longrightarrow
c_3\sum_{i=1}^m t_i b_{\ell,i,0}G(\cdot,\xi_i^0)
\quad
\text{in}\hspace{3mm}C^1_{\mathrm{loc}}
\bigl(\Om\setminus\{\xi_1^0,\ldots,\xi_m^0\}\bigr).
\end{gathered}
\end{equation}
The constant $C$ is independent of $\epsilon$.
\end{lem}
\begin{proof}
Write $u_\epsilon=V+\phi_\epsilon$, where
$V=\sum_{j=1}^mP_j$. By the definition of $PJ_{i\beta}$,
\[
\begin{aligned}
\mathcal A_uPJ_{i\beta}
=i^*\bigl(&3U_i^2(J_{i\beta}-PJ_{i\beta})
+3(U_i^2-V^2)PJ_{i\beta}\\
&-3(u_\epsilon^2-V^2)PJ_{i\beta}
-\epsilon PJ_{i\beta}\bigr).
\end{aligned}
\]
Since $|J_{i\beta}|\le CU_i$, the maximum principle gives
\[
|PJ_{i\beta}|\le CP_i\le CU_i.
\]
Moreover, $J_{i\beta}-PJ_{i\beta}$ is harmonic in $\Om$,
and its boundary values are bounded by $C\delta_i$.
Hence
\[
\|J_{i\beta}-PJ_{i\beta}\|_\infty\le C\delta_i.
\]
Together with
$\|U_i^2\|_{4/3}\le C\delta_i$ and
$\|PJ_{i\beta}\|_{4/3}\le C\delta_i$, this yields
\[
\|U_i^2(J_{i\beta}-PJ_{i\beta})\|_{4/3}
\le C\rho^2,
\qquad
\epsilon\|PJ_{i\beta}\|_{4/3}\le C\epsilon\rho.
\]

To estimate the interaction term, choose disjoint balls
of fixed radius around the centers.
On the ball around $\xi_i$, the projection estimates
give $V-U_i=O(\rho)$.
On each ball around $\xi_j$ with $j\ne i$,
we have $PJ_{i\beta}=O(\rho)$.
Outside these balls, all the bubbles and projected
Jacobi fields are $O(\rho)$.
Using $\|U_j^2\|_{4/3}=O(\rho)$, we therefore obtain
\[
\|(U_i^2-V^2)PJ_{i\beta}\|_{4/3}\le C\rho^2.
\]
Finally, H\"older's inequality and \eqref{eq:branch} give
\[
\begin{aligned}
\|(u_\epsilon^2-V^2)PJ_{i\beta}\|_{4/3}
&\le
\bigl(2\|V\|_4+\|\phi_\epsilon\|_4\bigr)
\|\phi_\epsilon\|_4\|PJ_{i\beta}\|_4\\
&\le C\gamma.
\end{aligned}
\]
The boundedness of
$i^*:L^{4/3}(\Om)\to H^1_0(\Om)$ now proves
\eqref{tp:Jacres}.
The same estimate holds with $u_\epsilon$ replaced by $V$.

Consider the compact self-adjoint operator
\[
T_uv:=i^*\bigl((3u_\epsilon^2+\epsilon)v\bigr)
\]
on $H^1_0(\Om)$ with the Dirichlet inner product;
see \cite[Section~3, (3.1)--(3.6)]{YZmulti}.
Its eigenvalues are $\mu_{\ell,\epsilon}^{-1}$.
Let $E_\epsilon$ be the orthogonal projection onto
\[
\operatorname{span}\{w_{m+1},\ldots,w_{6m}\}.
\]
By \eqref{tp:coarse}, the remaining spectrum of $T_u$
stays a fixed positive distance from $1$.
The spectral theorem and \eqref{tp:Jacres} thus imply
\[
\|(I-E_\epsilon)PJ_{i\beta}\|_{H^1_0(\Om)}
\le C\|\mathcal A_uPJ_{i\beta}\|_{H^1_0(\Om)}
\le C\gamma.
\]

Rescaling around the centers shows that the Dirichlet
Gram matrix of the fields $PJ_{i\beta}$ converges to
$g_1I_{5m}$. The preceding estimate gives the same limit
for the Gram matrix of $E_\epsilon PJ_{i\beta}$.
These $5m$ functions therefore form a basis of
$\operatorname{Ran}E_\epsilon$ for small $\epsilon$.
Consequently, for $m+1\le\ell\le6m$,
\[
w_\ell
=\sum_{i=1}^m\sum_{\beta=0}^4
b_{\ell,i,\beta,\epsilon}E_\epsilon PJ_{i\beta},
\qquad
|b_{\ell,\epsilon}|\le C,
\]
where the coefficient bound follows from
$\|w_\ell\|_{H^1_0(\Om)}^2=\mu_{\ell,\epsilon}\le C$
and the Gram matrix limit.
Setting
\[
\psi_\ell
:=\sum_{i=1}^m\sum_{\beta=0}^4
b_{\ell,i,\beta,\epsilon}(E_\epsilon-I)PJ_{i\beta},
\]
we obtain
$\|\psi_\ell\|_{H^1_0(\Om)}\le C\gamma$.
This proves \eqref{tp:thresholdprojection}.
Rescaling this decomposition and comparing with
\eqref{tp:coregeneral} also gives
\(
b_{\ell,i,\beta,\epsilon}\longrightarrow b_{\ell,i,\beta}
\)
along the subsequence in that formula.

It remains to prove \eqref{tp:outerw}.
The explicit expressions for $J_{i\beta}$, together
with the boundary values and harmonicity of
$J_{i\beta}-PJ_{i\beta}$, give
\[
PJ_{i0}
=c_3\delta_iG(\cdot,\xi_i)+O(\delta_i^3),
\qquad
PJ_{ik}=O(\delta_i^2),\quad 1\le k\le4,
\]
in $C^1(K)$ for every
$K\Subset\Om\setminus\{\xi_1^0,\ldots,\xi_m^0\}$.
Since $\delta_i/\rho\to t_i$ and
\[
\left\|\frac{\psi_\ell}{\rho}\right\|_{H^1_0(\Om)}
\le C(\epsilon+\rho)\longrightarrow0,
\]
the decomposition \eqref{tp:thresholdprojection} yields
\[
\rho^{-1}w_\ell
\longrightarrow
c_3\sum_{i=1}^m t_i b_{\ell,i,0}G(\cdot,\xi_i^0)
\quad\text{in }L^4_{\mathrm{loc}}
\bigl(\Om\setminus\{\xi_1^0,\ldots,\xi_m^0\}\bigr).
\]
Moreover,
\[
-\Delta(\rho^{-1}w_\ell)
=\mu_{\ell,\epsilon}(3u_\epsilon^2+\epsilon)
\,\rho^{-1}w_\ell.
\]
By Lemma~\ref{tp:profilelemma}, the coefficients in this
equation are locally uniformly bounded away from the
concentration points. Interior elliptic estimates
\cite[Chapter~9]{GT} therefore upgrade the preceding
convergence to $C^1_{\mathrm{loc}}$ convergence,
proving \eqref{tp:outerw}.
\end{proof}

\subsection{Sharp decay of the eigenfunctions}\label{sec:decay}
We establish the following a sharp decay estimate of the eigenfunctions.
\begin{lem}\label{tp:endpoint}
For every fixed $L\in\mathbb N$, there exist constants
$C_L,C_L'>0$ such that, for all sufficiently small
$\epsilon>0$, every eigenfunction $v_{\ell,\epsilon}$
of \eqref{eq:spectral} with
$
1\le\ell\le L$,
$\|v_{\ell,\epsilon}\|_{L^\infty(\Om)}=1,
$
satisfies
\begin{equation}\label{tp:endpointbound}
|v_{\ell,\epsilon}(x)|
\le C_L\sum_{i=1}^m
\frac{\delta_i^2}{\delta_i^2+|x-\xi_i|^2}
\le C_L'\frac{\rho^2}{\rho^2+d(x)^2},
\qquad x\in\Om.
\end{equation}
\end{lem}
\begin{proof}
The fixed-index eigenvalue bound $\mu_{\ell,\epsilon}\le C_L$ follows
from \eqref{tp:coarse} and, for larger fixed $L$, from
\cite[Theorem 1.2]{YZmulti}. Fix a large $R$ and set
\[
E_R=\Om\setminus\bigcup_i\overline{B(\xi_i,R\rho)},\qquad
B_\rho(x)=\rho^2\sum_i|x-\xi_i|^{-2}.
\]
For small $\epsilon$ the closed balls are disjoint and interior.
Let $h$ be harmonic in $E_R$, with boundary values $v$ on the
inner spheres and zero on $\partial\Om$.
Since $B_\rho$ is harmonic in $E_R$ and
$R^2B_\rho\ge1$ on every inner sphere, comparison gives
$|h|\le R^2B_\rho$.
The Green representation on $E_R$ is
\[
v(x)=h(x)+\mu_{\ell,\epsilon}\int_{E_R}
G_{E_R}(x,y)(3u^2(y)+\epsilon)v(y)\,dy.
\]
The Green function satisfies
$0\le G_{E_R}(x,y)\le(4\pi^2|x-y|^2)^{-1}$.
Moreover, \eqref{tp:udecay} gives
\[
u^2(y)B_\rho(y)
\le C\rho^4\sum_{i=1}^m|y-\xi_i|^{-6},
\qquad y\in E_R.
\]
For $a=|x-\xi_i|\ge R\rho$, we divide the integration
domain into the regions
$|y-\xi_i|<a/2$, $|x-y|<a/2$, and their complement.
Estimating the integral on each region yields
\begin{equation}\label{tp:kernel6}
\rho^4\int_{|y-\xi_i|\ge R\rho}
\frac{dy}{|x-y|^2|y-\xi_i|^6}
\le C\rho^4\{a^{-2}(R\rho)^{-2}+a^{-4}\}
\le \frac{C}{R^2}\frac{\rho^2}{a^2}.
\end{equation}
The first region uses the radial integral
$\int_{R\rho}^\infty r^{-3}\,dr$; on the second region
$|y-\xi_i|\ge a/2$ and $\int_{|x-y|<a/2}|x-y|^{-2}\,dy\le Ca^2$.
The remaining integral is $O(a^{-4})$ after scaling by $a$.
The linear term is controlled without a logarithmic loss:
the triangle inequality gives
\[
\frac{|x-\xi_i|^2}{|x-y|^2|y-\xi_i|^2}
\le2\{|x-y|^{-2}+|y-\xi_i|^{-2}\}.
\]
Integration over the bounded domain therefore shows
\begin{equation}\label{tp:kernel2}
\int_{E_R}G_{E_R}(x,y)B_\rho(y)\,dy\le C_\Om B_\rho(x).
\end{equation}
For each fixed $\epsilon$, the quantity
$M_\epsilon=\sup_{E_R}|v|/B_\rho$ is finite. Taking absolute values
in the representation, dividing by $B_\rho$, and using
\eqref{tp:kernel6}--\eqref{tp:kernel2} yields
\[
M_\epsilon\le R^2+C C_L(R^{-2}+\epsilon)M_\epsilon.
\]
Choose $R$ first and $\epsilon$ second so that the last coefficient
is at most $1/2$. Then $|v|\le2R^2B_\rho$ on $E_R$.
Since $c\rho\le\delta_i\le\rho$, this is the first bound in
\eqref{tp:endpointbound} there, with a constant depending on the
fixed $R$. Inside the core balls it follows from $|v|\le1$.
The last inequality uses the same comparability and the
nearest-center distance.
\end{proof}

Lemma~\ref{tp:endpoint} establishes the decay estimate
corresponding to \cite[Lemma~2.8]{CKL-1} for $N=4$ and $q=1$.
The comparison argument in
\cite[Lemma~4.5 and Appendix~A]{YZmulti}
requires a decay exponent strictly greater than $2$,
a condition that fails in this case.
Here \eqref{tp:kernel2} controls the linear perturbation,
while \eqref{tp:kernel6} allows the remaining integral
term to be absorbed, yielding the decay exponent $2$.

\section{Proof of Theorem \ref{thm:amp}}\label{sec:amplitude}
In this section, we are devoted to proving Theorem \ref{thm:amp}.
The following several results will be used to complete the proof of Theorem \ref{thm:amp}. We note that
the local Green identity for \eqref{intro:BN} and \eqref{eq:spectral} is given by
\begin{equation}\label{am:identity}
\int_{\partial B_i}(\partial_\nu u_\epsilon\,w-u_\epsilon\partial_\nu w)
=(3\mu-1)\int_{B_i}u_\epsilon^3w
+(\mu-1)\epsilon\int_{B_i}u_\epsilon w,
\end{equation}
where $B_i=B(\xi_{i,\epsilon},a)$ and $\nu=\nu(x)$ denotes the unit outward normal to the boundary $\partial B_i$.
The logarithmic coefficient of the second volume term
will be evaluated below.
The following lemma gives the cutoff bound at order $\rho_\epsilon^2$.
\begin{lem}\label{am:minmax}
Every eigenvalue satisfies $\mu_{\ell,\epsilon}>1/3$, and
\begin{equation}\label{am:upper}
0<\mu_{\ell,\epsilon}-\frac13\le C\rho_\epsilon^2,
\quad 1\le\ell\le m.
\end{equation}
\end{lem}
\begin{proof}
For $f\in C_c^\infty(\Om)$, the ground-state identity gives
\[
\int_\Om|\nabla f|^2-\int_\Om(u_\epsilon^2+\epsilon)f^2
=\int_\Om u_\epsilon^2|\nabla(f/u_\epsilon)|^2\ge0
\]
By density, the resulting inequality extends to $H^1_0(\Om)$.
Therefore
\[
\int|\nabla f|^2-\frac13D_\epsilon(f,f)
\ge\frac{2\epsilon}{3}\int f^2.
\]
Applying this to a nonzero eigenfunction proves the strict bound.
For the upper bound, use the disjoint cutoff trial functions
$\chi_i u_\epsilon$ as in \cite[proof of Theorem 1.1, p.~504]{CKL-1}.
Take $\chi_i=1$ near $\xi_{i,\epsilon}$, with fixed separated supports.
Testing \eqref{intro:BN} with $\chi_i^2u_\epsilon$ shows that the Rayleigh
quotient minus $1/3$ equals
\[
\frac{(2\epsilon/3)\int\chi_i^2u_\epsilon^2+
\int|\nabla\chi_i|^2u_\epsilon^2}
{3\int\chi_i^2u_\epsilon^4+\epsilon\int\chi_i^2u_\epsilon^2}.
\]
By Lemma~\ref{tp:profilelemma}, the denominator is bounded below by a positive
constant, while \eqref{tp:udecay} gives the required upper bound
\[
\int_\Om u_\epsilon^2\le C\rho_\epsilon^2(1+|\log\rho_\epsilon|),
\qquad \int|\nabla\chi_i|^2u_\epsilon^2\le C\rho_\epsilon^2.
\]
 So \eqref{eq:balance} makes the quotient $O(\rho_\epsilon^2)$.
Disjoint supports make the quotient of any linear combination
a weighted average. The min--max principle proves \eqref{am:upper}.
\end{proof}

The next lemma gives the local profiles and uniform
estimates for the first $m$ eigenfunctions.

\begin{lem}\label{am:profiles}
Let $w_{\ell,\epsilon}$, $1\le\ell\le m$, be normalized
by \eqref{tp:Dnormal}. Up to the subsequence,
the limits \eqref{am:corew} hold, and the coefficient
vectors satisfy
\begin{equation}\label{am:Gram}
g_0\,b_\ell\cdot b_j=\delta_{\ell j},
\qquad 1\le\ell,j\le m.
\end{equation}
Moreover, for all sufficiently small $\epsilon>0$,
\begin{equation}\label{am:wbound}
\begin{gathered}
|w_{\ell,\epsilon}(x)|
\le C\frac{\rho_\epsilon}
{\rho_\epsilon^2+d(x)^2},
\qquad x\in\Om,\quad
\int_\Om w_{\ell,\epsilon}^2\,dx
\le C\rho_\epsilon^2
\bigl(1+|\log\rho_\epsilon|\bigr),
\end{gathered}
\end{equation}
uniformly for $1\le\ell\le m$.
\end{lem}
\begin{proof}
Write $w_\ell=w_{\ell,\epsilon}$.
By \eqref{tp:Dnormal} and \eqref{eq:spectral},
\[
\|w_\ell\|_{H^1_0(\Om)}^2
=\mu_{\ell,\epsilon}\le C,
\qquad 1\le\ell\le m.
\]
Since $\mu_{\ell,\epsilon}\to1/3$, applying
\cite[Lemma~3.1, Proposition~3.2 and (3.3)]{YZmulti}
to the energy-normalized functions
$w_\ell/\sqrt{\mu_{\ell,\epsilon}}$ gives
\eqref{am:corew} and \eqref{am:Gram} along a common
subsequence, with
$g_0=3\int_{\mathbb R^4}U^4$.
The preceding energy bound, Lemma~\ref{tp:profilelemma},
and interior elliptic estimates
\cite[Chapter~9]{GT} imply
\[
\bigl\|\delta_iw_\ell(\xi_i+\delta_i\,\cdot)\bigr\|
_{C^{1,\sigma}(B_R)}
\le C_R
\]
for every fixed $R>0$ and some $\sigma\in(0,1)$,
uniformly in $i$, $\ell$, and sufficiently small $\epsilon$.
This also gives the local $C^1$ convergence in
\eqref{am:corew}.

Set $v_\ell=w_\ell/\|w_\ell\|_\infty$.
Lemma~\ref{tp:endpoint} gives
\begin{equation}\label{am:vbound}
|v_\ell(x)|
\le C\frac{\rho^2}{\rho^2+d(x)^2},
\qquad x\in\Om.
\end{equation}
Choose $x_\epsilon\in\Om$ such that
$|v_\ell(x_\epsilon)|=1$.
Then \eqref{am:vbound} implies
$d(x_\epsilon)\le C\rho$.
Since $\delta_i/\rho\to t_i>0$,
$x_\epsilon$ lies in a fixed rescaled ball around
some center $\xi_i$.
The local bound above therefore yields
\[
\|w_\ell\|_\infty
=|w_\ell(x_\epsilon)|
\le \frac{C}{\delta_i}
\le \frac{C}{\rho}.
\]
Multiplying \eqref{am:vbound} by this bound proves
the pointwise estimate in \eqref{am:wbound}.
Finally, integration in dimension four gives
\[
\begin{aligned}
\int_\Om w_\ell^2\,dx
&\le C\rho^2
\int_0^{\operatorname{diam}(\Om)}
\frac{r^3}{(\rho^2+r^2)^2}\,dr\\
&\le C\rho^2\bigl(1+|\log\rho|\bigr),
\end{aligned}
\]
which completes the proof.
\end{proof}

We next determine the exterior limits and the logarithmic
contribution to the local integrals of
$u_\epsilon w_{\ell,\epsilon}$.
We omit the indices $\ell,\epsilon$ when no confusion arises.

\begin{lem}\label{am:kernel}
Let $G_\zeta$ be the Dirichlet Green function of
$-\Delta-\zeta$ in $\Om$.
For $\zeta_0>0$ sufficiently small, there exists $C>0$
such that, for all $0\le\zeta\le\zeta_0$ and
$x,y\in\Om$ with $x\ne y$,
\[
0<G_\zeta(x,y)\le C|x-y|^{-2},
\qquad
|\nabla_xG_\zeta(x,y)|\le C|x-y|^{-3},
\]
and
\begin{equation}\label{am:resolvent}
|G_\zeta(x,y)-G(x,y)|
\le C\zeta\bigl(1+|\log|x-y||\bigr).
\end{equation}
Moreover, as $\zeta\to0$,
\[
G_\zeta\longrightarrow G
\quad\text{in }\;
C^1_{\mathrm{loc}}
\bigl(\{(x,y)\in\Om\times\Om:x\ne y\}\bigr).
\]
\end{lem}

\begin{proof}
Choose $\zeta_0<\lambda_1(-\Delta,\Om)$.
Positivity follows from the maximum principle, while
the bounds for $G_\zeta$ and $\nabla_xG_\zeta$ follow
from the uniform Green function estimates in
\cite[Appendix~B, (B.17)]{K-L-2024}.

By the resolvent identity,
\[
G_\zeta(x,y)-G(x,y)
=\zeta\int_\Om G(x,z)G_\zeta(z,y)\,dz.
\]
The preceding bounds and the four-dimensional estimate
\[
\int_\Om
\frac{dz}{|x-z|^2|z-y|^2}
\le C\bigl(1+|\log|x-y||\bigr)
\]
give \eqref{am:resolvent}.
Finally, away from $x=y$,
\[
-\Delta_x(G_\zeta-G)=\zeta G_\zeta.
\]
Hence \eqref{am:resolvent}, interior elliptic estimates
\cite[Chapter~9]{GT}, and the symmetry of the Green
functions imply the stated $C^1_{\mathrm{loc}}$ convergence.
\end{proof}

\begin{lem}\label{am:monopoles}
Let $b=b_\ell$ be a coefficient limit from
Lemma~\ref{am:profiles}.
Up to the subsequence,
\begin{align}
&\rho^{-1}u_\epsilon\longrightarrow c_3F_t,\quad
F_t=\sum_{i=1}^m t_iG(\cdot,\xi_i^0),
\label{am:outeru}\\
&\rho^{-1}w\longrightarrow c_3F_b,\quad
F_b=\sum_{i=1}^m t_i b_iG(\cdot,\xi_i^0),
\label{am:outerwb}
\end{align}
in $C^1_{\mathrm{loc}}
(\Om\setminus\{\xi_1^0,\ldots,\xi_m^0\})$.
Moreover, for every fixed sufficiently small $a>0$,
\begin{equation}\label{am:criticalmass}
\int_{B(\xi_i,a)}u_\epsilon^3w\,dx
\longrightarrow A b_i,
\qquad 1\le i\le m.
\end{equation}
\end{lem}

\begin{proof}
The convergence \eqref{am:outeru} follows from
Lemma~\ref{tp:profilelemma}.
By \eqref{tp:udecay} and \eqref{am:wbound},
\begin{equation}\label{am:localbound}
u_\epsilon(x)+|w(x)|
\le C\frac{\delta_i}
{(\delta_i+|x-\xi_i|)^2},
\qquad x\in B(\xi_i,a).
\end{equation}
The local profile limits in
Lemmas~\ref{tp:profilelemma} and \ref{am:profiles},
together with \eqref{am:localbound}, give
\eqref{am:criticalmass} by rescaling and dominated
convergence, since $A=\int_{\mathbb R^4}U^4$.

To prove \eqref{am:outerwb}, fix
$K\Subset\Om\setminus\{\xi_1^0,\ldots,\xi_m^0\}$
and choose $r>0$ so small that the balls
$B(\xi_i,2r)$ are mutually disjoint and avoid $K$
for all sufficiently small $\epsilon$.
The equation for $w$ gives
\[
\frac{w(x)}{\rho}
=\frac{3\mu}{\rho}
\int_\Om G_{\mu\epsilon}(x,y)
u_\epsilon^2(y)w(y)\,dy.
\]
On each $B(\xi_i,r)$, rescaling and dominated
convergence, using \eqref{am:localbound} and
Lemma~\ref{am:kernel}, yield
\[
\frac{3\mu}{\rho}
\int_{B(\xi_i,r)}
G_{\mu\epsilon}(\cdot,y)u_\epsilon^2(y)w(y)\,dy
\longrightarrow
c_3t_i b_iG(\cdot,\xi_i^0)
\quad\text{in }C^1(K).
\]
Here we used $3\mu\to1$, $\delta_i/\rho\to t_i$,
and $\int_{\mathbb R^4}U^3=c_3$.
On $\Om\setminus\bigcup_i B(\xi_i,r)$,
\eqref{tp:udecay} and \eqref{am:wbound} imply
$\rho^{-1}u_\epsilon^2|w|\le C_r\rho^2$.
The Green function and gradient bounds in
Lemma~\ref{am:kernel} therefore show that the remaining
integral is $O(\rho^2)$ in $C^1(K)$.
Summing over $i$ proves \eqref{am:outerwb}.
\end{proof}

Since $U^2\notin L^1(\mathbb R^4)$, we next determine
the logarithmic contribution to the integral of
$u_\epsilon w$.

\begin{lem}\label{am:logmoment}
Along the subsequence of Lemma~\ref{am:profiles},
for every fixed sufficiently small $a>0$,
\begin{equation}\label{am:loglimit}
\frac{\epsilon}{\rho^2}
\int_{B(\xi_i,a)}u_\epsilon w\,dx
\longrightarrow c_3^2\Lambda t_i^2b_i.
\end{equation}
\end{lem}

\begin{proof}
Fix $a_0>0$ such that the balls $B(\xi_j,2a_0)$
are mutually disjoint and contained in $\Om$
for all sufficiently small $\epsilon$.
Write $B_i=B(\xi_i,a_0)$ and $r_x=|x-\xi_i|$.
Consider
\[
(v,h,\zeta)
=(u_\epsilon,u_\epsilon^3,\epsilon)
\quad\text{or}\quad
(w,3\mu u_\epsilon^2w,\mu\epsilon).
\]
By \eqref{am:localbound},
\begin{equation}\label{am:sourcebound}
|h(y)|\le
C\frac{\delta_i^3}{(\delta_i+|y-\xi_i|)^6},
\qquad y\in B_i.
\end{equation}
The local profile limits in
Lemmas~\ref{tp:profilelemma} and \ref{am:profiles}
therefore give
$\delta_i^{-1}\int_{B_i}h\to c_3$ and $c_3b_i$,
respectively.

We first reduce the Green representation of $v$
to its singular part. Indeed,
\eqref{tp:udecay} and \eqref{am:wbound} imply
$\|h\|_{L^1(\Om)}\le C\rho$ and
$\|h\|_{L^\infty(\Om)}\le C\rho^{-3}$.
Splitting the logarithmic integral at $|x-y|=\rho$,
\eqref{am:resolvent} and \eqref{eq:balance} yield
\[
\left|
\int_\Om
\bigl(G_\zeta(x,y)-G(x,y)\bigr)h(y)\,dy
\right|
\le C\epsilon\rho(1+|\log\rho|)
\le C\rho.
\]
The regular part of $G$ and the integral over
$\Om\setminus B_i$ are also $O(\rho)$ for
$r_x\le a_0/2$. Hence
\[
v(x)=\frac{1}{4\pi^2}
\int_{B_i}\frac{h(y)}{|x-y|^2}\,dy+O(\rho).
\]
For $\delta_i\le r_x\le a_0/2$, splitting $B_i$
at $|y-\xi_i|=r_x/2$ and using
\eqref{am:sourcebound}, we obtain
\[
\left|
\int_{B_i}
\left(\frac{1}{|x-y|^2}-\frac{1}{r_x^2}\right)
h(y)\,dy
\right|
\le C\left(
\frac{\delta_i^2}{r_x^3}
+\frac{\delta_i^3}{r_x^4}
\right)
\le C\frac{\delta_i^2}{r_x^3}.
\]
Here the inner region is estimated using
$\int_{B_i}|y-\xi_i||h(y)|\,dy\le C\delta_i^2$,
and the remaining region by rescaling with $r_x$.
Since $\delta_i\asymp\rho$ and
$c_3/(4\pi^2)=\alpha$, it follows that
\begin{equation}\label{am:annulus}
\frac{r_x^2}{\delta_i}u_\epsilon(x)
=\alpha+e_{i,\epsilon}(x),
\qquad
\frac{r_x^2}{\delta_i}w(x)
=\alpha b_i+f_{i,\epsilon}(x),
\end{equation}
where, for some $\eta_\epsilon\to0$,
\begin{equation}\label{am:annerror}
|e_{i,\epsilon}(x)|+|f_{i,\epsilon}(x)|
\le \eta_\epsilon+
C\left(\frac{\delta_i}{r_x}+r_x^2\right),
\qquad
\delta_i\le r_x\le a_0/2.
\end{equation}
Fix $0<a\le a_0/2$.
Multiplying \eqref{am:annulus} and integrating
\eqref{am:annerror}, we obtain
\begin{equation}\label{am:annintegral}
\int_{\delta_i<r_x<a}u_\epsilon w\,dx
=
16\pi^2b_i\delta_i^2\log\frac{a}{\delta_i}
+o\!\left(\delta_i^2\log\delta_i^{-1}\right)
+O(\delta_i^2),
\end{equation}
where $\alpha^2|\mathbb S^3|=16\pi^2$.
By \eqref{am:localbound}, the integral over
$B(\xi_i,\delta_i)$ is $O(\delta_i^2)$.
Finally, \eqref{eq:widths}--\eqref{eq:balance} give
\[
\epsilon\log\delta_i^{-1}\longrightarrow8\pi^2\Lambda,
\qquad
\frac{\delta_i^2}{\rho^2}\longrightarrow t_i^2.
\]
Multiplying \eqref{am:annintegral} by
$\epsilon/\rho^2$ and using
$16\pi^2\cdot8\pi^2=c_3^2$
proves \eqref{am:loglimit}.
\end{proof}

We next evaluate the boundary integral in \eqref{am:identity}
to complete the proof of Theorem~\ref{thm:amp}.
\begin{lem}\label{am:flux}
Let $a>0$ be fixed and sufficiently small, and set
$B_i=B(\xi_{i,\epsilon},a)$.
With $\nu$ denoting the outward unit normal on
$\partial B_i$, we have
\begin{equation}\label{am:fluxlimit}
\rho^{-2}\int_{\partial B_i}
\bigl(\partial_\nu u_\epsilon\,w
-u_\epsilon\partial_\nu w\bigr)\,dS
\longrightarrow c_3^2(\cS b)_i.
\end{equation}
\end{lem}

\begin{proof}
By Lemma~\ref{am:monopoles}, the left-hand side
of \eqref{am:fluxlimit} converges to $c_3^2 I_i$, where
\[
I_i:=
\int_{\partial B(\xi_i^0,a)}
\bigl(\partial_\nu F_t\,F_b
-F_t\partial_\nu F_b\bigr)\,dS.
\]
Writing
$\Phi_i(x)=(4\pi^2|x-\xi_i^0|^2)^{-1}$, we have
\[
F_t=t_i\Phi_i+F_t^{\rm reg},
\qquad
F_b=t_ib_i\Phi_i+F_b^{\rm reg},
\]
where the regular parts are harmonic in
$B(\xi_i^0,a)$ and satisfy
\[
F_t^{\rm reg}(\xi_i^0)=-(Mt)_i=-\Lambda t_i,
\qquad
F_b^{\rm reg}(\xi_i^0)=-(MD_tb)_i.
\]
Green's identity, as in \cite[Section~3]{CKL-1},
together with
$\int_{\partial B(\xi_i^0,a)}\partial_\nu\Phi_i\,dS=-1$,
therefore gives
\[
\begin{aligned}
I_i
&=-t_iF_b^{\rm reg}(\xi_i^0)
  +t_ib_iF_t^{\rm reg}(\xi_i^0)\\
&=t_i(MD_tb)_i-\Lambda t_i^2b_i
=(\cS b)_i.
\end{aligned}
\]
This proves \eqref{am:fluxlimit}.
\end{proof}

Now we are ready to prove Theorem \ref{thm:amp}.
\begin{proof}[Proof of Theorem~\ref{thm:amp}]
Let $\epsilon\to0$ along an arbitrary sequence.
By \eqref{am:upper} and Lemma~\ref{am:profiles},
we may pass to a common subsequence such that
the profile limits defining $b_\ell$ hold and
\[
\eta_\ell:=
\lim_{\epsilon\to0}
\frac{\mu_{\ell,\epsilon}-1/3}{\rho_\epsilon^2}
\]
exists for every $1\le\ell\le m$.
Applying \eqref{am:identity} to $w_{\ell,\epsilon}$,
dividing by $\rho_\epsilon^2$, and using
Lemmas~\ref{am:monopoles}, \ref{am:logmoment}
and \ref{am:flux}, we obtain
\[
c_3^2(\cS b_\ell)_i
=
3A\eta_\ell b_{\ell,i}
-\frac23c_3^2\Lambda t_i^2b_{\ell,i},
\qquad 1\le i\le m.
\]
By \eqref{eq:amplitude} and
$c_3^2/(3A)=4\pi^2$, this is equivalent to
\begin{equation}\label{am:finiteproblem}
\mathcal A_1b_\ell
=\frac{\eta_\ell}{4\pi^2}b_\ell,
\qquad 1\le\ell\le m.
\end{equation}

By \eqref{am:Gram}, the vectors $b_1,\ldots,b_m$
form an orthogonal basis of $\mathbb R^m$.
Consequently, \eqref{am:finiteproblem} identifies
$\eta_\ell/(4\pi^2)$, $1\le\ell\le m$,
with all eigenvalues of $\mathcal A_1$, counted
with multiplicity, as in
\cite[Section~3, proof of Theorem~1.1]{CKL-1}.
The ordering of $\mu_{\ell,\epsilon}$ implies
$\eta_1\le\cdots\le\eta_m$, and hence
\[
\eta_\ell=4\pi^2\sigma_\ell,
\qquad 1\le\ell\le m.
\]
Since the initial sequence was arbitrary,
\eqref{am:shift} follows for the full family.
Together with \eqref{am:Gram},
\eqref{am:finiteproblem} also proves
\eqref{am:bmatrix}.
The local and exterior profile limits follow from
Lemmas~\ref{am:profiles} and \ref{am:monopoles},
while the positive definiteness of $\mathcal A_1$
is established in Lemma~\ref{lem:matrix}.
\end{proof}

\section{Proof of Theorem~\ref{thm:pure}}\label{sec:pure}
In this section, we prove Theorem~\ref{thm:pure}.
We first give the upper bounds for the $\ell$-th eigenvalues.
\begin{Prop}\label{tp:upper}
For $m+1\le\ell\le5m$,
\begin{equation}\label{tp:upperbound}
\mu_{\ell,\epsilon}\le1+C\rho^4.
\end{equation}
\end{Prop}
\begin{proof}
We follow the cutoff min--max argument of
\cite[Proposition~4.1]{CKL-1}; see also the analogous
estimate in \cite[Lemma~6.2]{YZmulti}.
Choose cutoffs $\chi_i$ with mutually disjoint supports,
equal to one in $B(\xi_i,r)$ and supported in
$B(\xi_i,2r)$, where $r>0$ is fixed and sufficiently small.
Consider the space
\[
\mathcal V_\epsilon
=\operatorname{span}
\left\{
\chi_i u,\ \delta_i\chi_i\partial_k u:
1\le i\le m,\ 1\le k\le4
\right\}.
\]
Fix $i$, and write
\[
z=-\delta_i a\cdot\nabla u,
\qquad
f_i=\chi_i(a_0u+z),
\qquad
a_0\in\mathbb R,\quad a\in\mathbb R^4.
\]
Since $-\Delta z=(3u^2+\epsilon)z$, the cutoff identities
used in \cite[Proposition~4.1]{CKL-1} give
\[
\begin{aligned}
Q(\chi_i u,\chi_i u)&\le-c,\quad
|Q(\chi_i u,\chi_i z)|&\le C\rho^3|a|,\quad
Q(\chi_i z,\chi_i z)
&=\int_\Om|\nabla\chi_i|^2z^2\,dx
\le C\rho^4|a|^2.
\end{aligned}
\]
Here we used Lemma~\ref{tp:profilelemma}, which implies
$u,\nabla u=O(\rho)$ and
$z,\nabla z=O(\rho^2|a|)$ on
$\operatorname{supp}\nabla\chi_i$.
Young's inequality therefore yields
\[
Q(f_i,f_i)\le-ca_0^2+C\rho^4|a|^2.
\]

The local convergence in Lemma~\ref{tp:profilelemma}
and the linear independence of $U,Z_1,\ldots,Z_4$
give, by rescaling on $B(\xi_i,\delta_i)$,
\[
D(f_i,f_i)\ge c(a_0^2+|a|^2).
\]
In particular, $\dim\mathcal V_\epsilon=5m$.
Since the cutoff supports are disjoint, the preceding
estimates imply
\[
\sup_{f\in\mathcal V_\epsilon\setminus\{0\}}
\frac{\int_\Om|\nabla f|^2\,dx}{D(f,f)}
=
1+
\sup_{f\in\mathcal V_\epsilon\setminus\{0\}}
\frac{Q(f,f)}{D(f,f)}
\le1+C\rho^4.
\]
The min--max principle gives
$\mu_{5m,\epsilon}\le1+C\rho^4$.
The ordering of the eigenvalues then proves
\eqref{tp:upperbound}.
\end{proof}

\subsection{Local dilation identities}\label{sec:local}

For a fixed sufficiently small $r>0$, let
$B_i=B(\xi_i,r)$, and denote by $\nu$ the outward
unit normal on $\partial B_i$.
Set
\[
\omega_i=u+(x-\xi_i)\cdot\nabla u,
\qquad
L_u=-\Delta-3u^2-\epsilon,
\]
and
\[
F_i=
\int_{\partial B_i}
\bigl(\partial_\nu\omega_i\,w
-\omega_i\partial_\nu w\bigr)\,dS.
\]
Specializing \cite[(5.5), (5.7)]{YZmulti}
to $(p,q)=(3,1)$, we obtain
\begin{equation}\label{tp:scaleidentity}
\begin{gathered}
L_u\omega_i=2\epsilon u,\quad
F_i=(\mu-1)\int_{B_i}(3u^2+\epsilon)\omega_iw\,dx
-2\epsilon\int_{B_i}uw\,dx.
\end{gathered}
\end{equation}
\begin{lem}\label{tp:moments}
Let $w=w_{\ell,\epsilon}$, $m+1\le\ell\le6m$,
be normalized by \eqref{tp:Dnormal}, and set
\[
b^{(0)}=(b_i^{(0)})_{i=1}^m,
\qquad b_i^{(0)}=b_{\ell,i,0},
\]
where $b_{\ell,i,0}$ are the dilation coefficients
in \eqref{tp:coregeneral}.
Along the corresponding subsequence,
\begin{align}
\rho^{-2}F_i
&\longrightarrow
-c_3^2
\bigl[D_t(M+\Lambda\Id)D_tb^{(0)}\bigr]_i,
\label{tp:fluxlimit}\\
\int_{B_i}(3u^2+\epsilon)\omega_iw\,dx
&\longrightarrow -g_1b_i^{(0)},
\label{tp:criticalmoment}\\
\frac{\epsilon}{\rho^2}\int_{B_i}uw\,dx
&\longrightarrow c_3^2\Lambda t_i^2b_i^{(0)}.
\label{tp:logmoment}
\end{align}
\end{lem}

\begin{proof}
The limit \eqref{tp:fluxlimit} follows from
\eqref{tp:outeru} and \eqref{tp:outerw},
the boundary integral identity in
\cite[Lemma~4.5]{CKL-1}, and $Mt=\Lambda t$.
For \eqref{tp:criticalmoment}, we use the rescaling
and parity argument of \cite[Lemma~6.7]{YZmulti}.
By Lemma~\ref{tp:profilelemma},
\[
\delta_i\omega_i(\xi_i+\delta_i y)
\longrightarrow -Z_0(y)
\quad\text{in }C^0_{\mathrm{loc}}(\mathbb R^4).
\]
Together with \eqref{tp:coregeneral}, this gives
\[
3\int_{B_i}u^2\omega_iw\,dx
\longrightarrow
-3b_i^{(0)}\int_{\mathbb R^4}U^2Z_0^2\,dy
=-g_1b_i^{(0)}.
\]
Indeed, passage from fixed rescaled balls to $B_i$
is justified by \eqref{tp:udecay},
H\"older's inequality, and the bounded
$H^1_0(\Om)$ norms of $u$ and $w$, which yield
\[
\int_{B_i\setminus B(\xi_i,R\delta_i)}
u^2|\omega_iw|\,dx\le CR^{-2}.
\]
Moreover, $\|\omega_i\|_{L^2(B_i)}\le C$, so that
$\epsilon\int_{B_i}\omega_iw\,dx=O(\epsilon)=o(1)$.
This proves \eqref{tp:criticalmoment}.

For \eqref{tp:logmoment}, the decompositions
\eqref{eq:branch} and \eqref{tp:thresholdprojection}
give
\[
\int_{B_i}uw\,dx
=
b_{\ell,i,0,\epsilon}
\int_{B_i}P_iPJ_{i0}\,dx
+O(\rho^2+\rho\gamma+\gamma^2).
\]
Here the other bubble products are $O(\rho^2)$
by separation and oddness. The terms containing
$\phi_\epsilon$ or $\psi_{\ell,\epsilon}$ are
$O(\rho\gamma+\gamma^2)$ by H\"older's inequality,
using their $L^4$ bounds and the $O(\rho)$
bounds for the bubble factors in $L^{4/3}$.

Using the same radial computation as in
\cite[proof of Proposition~A.3]{P-R-V},
together with the projection estimates, we obtain
\begin{equation}\label{tp:logproduct}
\int_{B_i}P_iPJ_{i0}\,dx
=
16\pi^2\delta_i^2\log\delta_i^{-1}
+O(\rho^2).
\end{equation}
Since $\gamma/\rho=\epsilon+\rho\to0$,
multiplication by $\epsilon/\rho^2$ gives
\[
\frac{\epsilon}{\rho^2}\int_{B_i}uw\,dx
=
16\pi^2b_{\ell,i,0,\epsilon}
\frac{\delta_i^2}{\rho^2}
\epsilon\log\delta_i^{-1}+o(1).
\]
Now $b_{\ell,i,0,\epsilon}\to b_i^{(0)}$,
and \eqref{eq:widths}--\eqref{eq:balance}
yield \eqref{tp:logmoment}.
\end{proof}

\begin{Prop}\label{tp:leadingconstraint}
Let $b^{(0)}$ be as in Lemma~\ref{tp:moments}.
If $b^{(0)}\ne0$, then, along the same subsequence,
\begin{equation}\label{tp:graphconstraint}
\frac{\mu-1}{\rho^2}\longrightarrow\eta,
\qquad
\eta b^{(0)}=20\pi^2\cS b^{(0)}.
\end{equation}
In particular, for $m+1\le\ell\le5m$,
\[
b^{(0)}\in\operatorname{span}\{(1,\ldots,1)^T\}.
\]
\end{Prop}

\begin{proof}
Choose $i$ with $b_i^{(0)}\ne0$.
Dividing \eqref{tp:scaleidentity} by $\rho^2$ and
using Lemma~\ref{tp:moments}, we obtain the convergence
of $(\mu-1)/\rho^2$.
Applying the same identity for each $i$ gives
\eqref{tp:graphconstraint}, since $c_3^2/g_1=20\pi^2$.
For $m+1\le\ell\le5m$, \eqref{tp:upperbound} gives
$\eta\le0$. Lemma~\ref{lem:matrix} and
\eqref{tp:graphconstraint} therefore imply
$\eta=0$ and
$b^{(0)}\in\ker\cS
=\operatorname{span}\{(1,\ldots,1)^T\}$.
\end{proof}

\subsection{Scale reduction}\label{sec:scale}

We next analyze the remaining common dilation component
at the finer scale $\epsilon\rho^2$.
With $\epsilon$ and the centers fixed, write
\begin{equation}\label{tp:window}
\rho=e^{-s},\qquad
\delta_i=\rho t_i,\qquad
|t|=1,\quad t_i\ge c>0,\qquad
0<c_0\le\epsilon s\le C_0.
\end{equation}
Let $h$ denote fixed smooth local coordinates for $t$,
and let $\Pi_X$ be the Dirichlet orthogonal projection
onto $X=\cK^\perp$.
All estimates are uniform for the centers in compact
separated sets, and derivatives are taken only in $(s,h)$.
When $m=1$, $t_1=1$ and the $h$ variables are absent.

\begin{lem}\label{tp:scaledifferentiation}
Let $D_i=\delta_i\partial_{\delta_i}$ and
$H_i=H(\cdot,\xi_i)$.
Uniformly for the parameters in \eqref{tp:window},
and for $0\le k\le3$,
\begin{align}
|D_i^kU_i|+|D_i^kP_i|
&\le C_kU_i,
&
\|D_i^kP_i\|_{H^1_0(\Om)}
&\le C_k,
\label{tp:logderivative}\\
\big\|D_i^k(P_i-U_i+c_3\delta_iH_i)\big\|_
{C^2(\overline\Om)}
&\le C_k\delta_i^3.
\label{tp:projectionderivative}
\end{align}
The projected Jacobi fields and their first two
$(s,h)$ derivatives are uniformly bounded in $H^1_0(\Om)$.
Moreover,
\begin{equation}\label{tp:radialbounds}
\begin{gathered}
\|U_i\|_{4/3}\le C\delta_i,
\qquad
\|U_i^2\|_{4/3}\le C\delta_i,\\
\int_\Om U_i^2\,dx\le C\delta_i^2(1+s),
\qquad
\int_\Om U_i^3|x-\xi_i|\,dx\le C\delta_i^2.
\end{gathered}
\end{equation}
These bounds remain valid after up to two logarithmic
width derivatives of the integrands, with absolute
values taken.
\end{lem}

\begin{proof}
The explicit formula \eqref{eq:bubble} gives
$|D_i^kU_i|\le C_kU_i$ and
$|D_i^k(U_i^3)|\le C_kU_i^3$.
Since
$-\Delta D_i^kP_i=D_i^k(U_i^3)$ with zero boundary data,
the maximum principle and the boundedness of
$i^*:L^{4/3}(\Om)\to H^1_0(\Om)$ yield
\eqref{tp:logderivative}.

The undifferentiated projection expansion is
\cite[(2.4)]{P-R-V}.
To justify \eqref{tp:projectionderivative}, set
$R_i=P_i-U_i+c_3\delta_iH_i$.
Then $R_i$ is harmonic and
\[
R_i(x)=
\frac{\alpha\delta_i^3}
{|x-\xi_i|^2(\delta_i^2+|x-\xi_i|^2)}
\qquad\text{on }\partial\Om.
\]
Because the centers remain away from the boundary,
these boundary data and their first three $D_i$
derivatives have $C^{2,\sigma}(\partial\Om)$ norm
$O(\delta_i^3)$, for a fixed $\sigma\in(0,1)$.
Thus \eqref{tp:projectionderivative} follows from
the boundary Schauder estimates
\cite[Chapter~6]{GT}.
Similarly,
$|D_i^j(3U_i^2J_{i\beta})|\le CU_i^3$ for $j\le2$.
Applying the same $H^1_0$ estimate to
\eqref{eq:scaled-Jacobi} gives the asserted bounds
for the projected Jacobi fields; the change to
$(s,h)$ preserves these bounds.
Finally, \eqref{tp:radialbounds} follows by radial
integration, and its differentiated versions follow
from the same pointwise estimates.
\end{proof}

\begin{lem}\label{tp:scaleLS}
For sufficiently small $\epsilon$, uniformly for the
parameters in \eqref{tp:window}, the operator
$\mathcal A_{V,X}:X\to X$ defined by
\begin{equation}\label{tp:AV}
\mathcal A_{V,X}v
=\Pi_X\bigl[v-i^*((3V^2+\epsilon)v)\bigr]
\end{equation}
is invertible with uniformly bounded inverse.
The projected equation
\begin{equation}\label{tp:projected}
\Pi_X\bigl[\Psi-i^*(\Psi^3+\epsilon\Psi)\bigr]=0,
\qquad \Psi=V+\phi,
\end{equation}
admits a unique small solution
$\phi=\phi(s,h;\bm\xi)\in X$.
With $\epsilon$ and $\bm\xi$ fixed, this solution
depends $C^2$ on $(s,h)$ as an $H^1_0(\Om)$-valued map,
and
\begin{equation}\label{tp:C2phi}
\|\partial_{(s,h)}^\alpha\phi\|_{H^1_0(\Om)}
\le C_\alpha\gamma,
\qquad |\alpha|\le2.
\end{equation}
At the parameters corresponding to $u_\epsilon$,
the operator
$\mathcal A_{u,X}=\Pi_X\mathcal A_u|_X$
is also uniformly invertible.
\end{lem}
\begin{proof}
We first prove that $\mathcal A_{V,X}$ is uniformly invertible.
In this argument, norms and inner products refer to the
Dirichlet inner product on $H^1_0(\Om)$. Set
\[
T_Vv=i^*((3V^2+\epsilon)v),
\qquad
\mathcal A_V=\Id-T_V.
\]
The operator $T_V$ is compact and self-adjoint. Denote by
$\mu_{\ell,\epsilon}^V$ the ordered eigenvalues of
$-\Delta w=\mu(3V^2+\epsilon)w$ with zero Dirichlet boundary
data. Since $V=\sum_iP_i$, we may apply
\cite[Theorem~1.2]{YZmulti} with zero remainder; this theorem
does not require $V$ to solve the nonlinear equation. Hence
\[
\mu_{m,\epsilon}^V\longrightarrow\frac13,
\qquad
\mu_{\ell,\epsilon}^V\longrightarrow1
\quad(m+1\le\ell\le6m),
\qquad
\mu_{6m+1,\epsilon}^V\longrightarrow2.
\]

Let $E_V$ be the Dirichlet orthogonal projection onto the
span of the eigenfunctions corresponding to
$\mu_{m+1,\epsilon}^V,\ldots,\mu_{6m,\epsilon}^V$.
The eigenvalues of $\mathcal A_V$ are
$1-(\mu_{\ell,\epsilon}^V)^{-1}$, so
\[
\|\mathcal A_Vz\|\ge c_*\|z\|,
\qquad z\in\ker E_V,
\]
for some $c_*>0$.
This constant is uniform for the parameters in
\eqref{tp:window} and the centers in the fixed compact
separated set. Indeed, if not, one could choose a sequence
for which the bound fails and extract a subsequence of the
centers. Since $\rho\le e^{-c_0/\epsilon}\to0$ and the limiting
centers remain distinct and interior, the same spectral
theorem would give a contradiction.
The estimates in the proof of
Lemma~\ref{tp:quantprojection}, applied to $V$, and the
Gram matrix bounds for the fields $PJ_{i\beta}$ give
\[
\|\mathcal A_Vk\|\le C\gamma\|k\|,
\qquad k\in\cK.
\]
Since $E_V$ commutes with $\mathcal A_V$, the spectral gap
implies
\[
\|(\Id-E_V)k\|
\le c_*^{-1}\|(\Id-E_V)\mathcal A_Vk\|
\le C\gamma\|k\|,
\qquad k\in\cK.
\]
Thus $E_V|_{\cK}$ is injective for small $\epsilon$.
As $\dim\cK=\dim\operatorname{Ran}E_V=5m$, it is also onto,
and its inverse is uniformly bounded. For every
$e\in\operatorname{Ran}E_V$, we can therefore choose
$k\in\cK$ such that
\[
E_Vk=e,\qquad
\|k\|\le2\|e\|,\qquad
\|e-k\|\le C\gamma\|e\|.
\]
If $v\in X=\cK^\perp$, then
\[
|\langle v,e\rangle|
=|\langle v,e-k\rangle|
\le C\gamma\|v\|\|e\|.
\]
Taking the supremum over unit vectors
$e\in\operatorname{Ran}E_V$ yields
\[
\|E_Vv\|\le C\gamma\|v\|,
\qquad v\in X.
\]
Let $\Pi_{\cK}=\Id-\Pi_X$. By self-adjointness,
\[
\begin{aligned}
\|\Pi_{\cK}\mathcal A_Vv\|
&=\sup_{\substack{k\in\cK\\ \|k\|=1}}
  |\langle v,\mathcal A_Vk\rangle|
\le C\gamma\|v\|.
\end{aligned}
\]
Using again that $E_V$ commutes with $\mathcal A_V$, we obtain
\[
\begin{aligned}
\|\mathcal A_{V,X}v\|
&\ge \|\mathcal A_Vv\|
      -\|\Pi_{\cK}\mathcal A_Vv\|\\
&\ge c_*\|(\Id-E_V)v\|-C\gamma\|v\|\\
&\ge (c_*-C\gamma)\|v\|
\ge c\|v\|,
\qquad v\in X,
\end{aligned}
\]
where $c>0$ is independent of the parameters, since
$\gamma\to0$ uniformly under \eqref{tp:window}.
Finally,
$\mathcal A_{V,X}=\Id_X-\Pi_XT_V|_X$
is a compact perturbation of the identity on $X$.
The lower bound gives injectivity, and the Fredholm
alternative gives surjectivity. Consequently,
\[
\|\mathcal A_{V,X}^{-1}\|_{\mathcal L(X)}\le c^{-1}.
\]

Let $p=(s,h)$ and set
\[
\mathcal R_p=V-i^*(V^3+\epsilon V).
\]
By Lemma~\ref{tp:scaledifferentiation}, the estimates
in \cite[Lemma~A.1]{P-R-V} remain valid after
at most two $p$ derivatives. Hence
\begin{equation}\label{tp:residualC2}
\|\partial_p^\alpha\mathcal R_p\|_{H^1_0(\Om)}
\le C_\alpha\gamma,
\qquad |\alpha|\le2.
\end{equation}
The contraction argument of
\cite[Proposition~2.1]{P-R-V} applies also to the
polynomial cubic and yields a unique small solution
$\phi\in X$ with $\|\phi\|_{H^1_0(\Om)}\le C\gamma$.

It remains to prove the derivative estimates.
Write $\Pi_p=\Pi_X$.
The uniform Gram matrix bounds and
Lemma~\ref{tp:scaledifferentiation} imply that
$\Pi_p$ depends $C^2$ on $p$, with uniformly bounded
operator derivatives.
On the fixed space $H^1_0(\Om)$, define
\[
\begin{aligned}
\mathcal F(p,\varphi)
={}&\Pi_p\bigl[
V+\varphi-i^*((V+\varphi)^3+\epsilon(V+\varphi))
\bigr]\\
&+(\Id-\Pi_p)\varphi.
\end{aligned}
\]
Thus $\mathcal F(p,\varphi)=0$ is equivalent to
\eqref{tp:projected} and $\varphi\in X$.
At $\varphi=\phi$, its derivative is
\[
D_\varphi\mathcal F(p,\phi)
=\Pi_p\mathcal A_\Psi+(\Id-\Pi_p).
\]
Relative to $H^1_0(\Om)=X\oplus\cK$, this operator
is block triangular, with diagonal blocks
$\Pi_p\mathcal A_\Psi|_X$ and $\Id$.
It is uniformly invertible, since
\[
\|\mathcal A_\Psi-\mathcal A_V\|
\le C\|\Psi^2-V^2\|_{L^2(\Om)}
\le C\gamma.
\]
The implicit function theorem therefore gives
the $C^2$ dependence of $\phi$ on $p$.
By \eqref{tp:residualC2} and
$\|\phi\|_{H^1_0(\Om)}=O(\gamma)$, the first two
$p$ derivatives of $\mathcal F(p,\phi)$,
with $\phi$ held fixed, are $O(\gamma)$;
the required mixed derivatives are uniformly bounded.
Differentiating $\mathcal F(p,\phi(p))=0$
once and twice and applying the uniform inverse
gives \eqref{tp:C2phi}.

Finally, \eqref{eq:branch} yields
$\|u^2-V^2\|_{L^2(\Om)}=O(\gamma)$.
The same perturbation argument proves the uniform
invertibility of $\mathcal A_{u,X}$.
\end{proof}

By uniqueness in Lemma~\ref{tp:scaleLS},
$\Psi=u_\epsilon$ at the parameters in \eqref{eq:branch}.
Define
\[
E_\epsilon(f)
=\frac12\int_\Om|\nabla f|^2\,dx
-\frac14\int_\Om f^4\,dx
-\frac{\epsilon}{2}\int_\Om f^2\,dx,
\qquad
F_\epsilon(s,h;\bm\xi)=E_\epsilon(\Psi).
\]

\begin{Prop}\label{tp:energy}
Set
\[
a=\frac{c_3^2}{2}=64\pi^4,
\qquad b=8\pi^2,
\qquad C_m=\frac{mA}{4}.
\]
Uniformly for the parameters in \eqref{tp:window},
\begin{equation}\label{tp:F}
F_\epsilon
=C_m+\rho^2
\bigl[
a\,t^TM(\bm\xi)t-b\epsilon s
+\epsilon\beta(t;\bm\xi)
\bigr]+R_\epsilon,
\end{equation}
where $\beta$ is independent of $s$ and has uniformly
bounded $h$ derivatives up to order two.
For some fixed integer $K\ge0$, the remainder satisfies
\begin{equation}\label{tp:R}
\begin{gathered}
|\partial_{s,h}^{\alpha}R_\epsilon|
\le C\rho^2\eta_\epsilon,
\qquad |\alpha|\le2,\quad
\eta_\epsilon
:=\epsilon^2+\rho(1+s)^K=o(\epsilon).
\end{gathered}
\end{equation}
\end{Prop}
\begin{proof}
All remainder estimates below include derivatives
up to order two in $(s,h)$.
Integration by parts, the projection expansion
\cite[(2.4)--(2.5)]{P-R-V}, and
Lemma~\ref{tp:scaledifferentiation} give
\begin{equation}\label{tp:bareenergy}
E_0(V)
=C_m+\frac{c_3^2}{2}\bm\delta^TM(\bm\xi)\bm\delta
+O_{C^2(s,h)}\bigl(\rho^3(1+s)^K\bigr).
\end{equation}
For the $L^2$ term, we have
\begin{equation}\label{tp:mass}
\int_\Om V^2\,dx
=16\pi^2\sum_i\delta_i^2\log\delta_i^{-1}
+\bm\delta^TB(\bm\xi)\bm\delta
+O_{C^2(s,h)}\bigl(\rho^4(1+s)\bigr),
\end{equation}
where $B(\bm\xi)$ is a uniformly bounded symmetric
matrix independent of the widths.
Indeed, for fixed $r>0$, direct radial integration gives
\[
\int_{B_r(0)}U_{\delta,0}^2\,dx
=16\pi^2\delta^2\log(r/\delta)-8\pi^2\delta^2
+O_{C^2(\log\delta)}(\delta^4).
\]
After evaluating this singular square, the remaining
terms are estimated using
\[
q_\delta(r)
:=\frac{1}{\delta^2+r^2}-\frac{1}{r^2},
\qquad
|(\delta\partial_\delta)^kq_\delta(r)|
\le
\frac{C\delta^2}{r^2(\delta^2+r^2)},
\quad 0\le k\le2.
\]
The right-hand side has integral
$O(\delta^2(1+|\log\delta|))$ on a fixed ball.
Applied to products with factors smooth near the
corresponding center, this estimate and
\eqref{tp:projectionderivative} control the remainder
in \eqref{tp:mass}.
Products involving distinct centers are treated on
disjoint fixed balls.
The remaining coefficients are independent of the
widths and define $B(\bm\xi)$.
Since $\delta_i=\rho t_i$,
$\log\delta_i^{-1}=s-\log t_i$ and $\sum_i t_i^2=1$,
combining \eqref{tp:bareenergy} and \eqref{tp:mass}
gives the required expansion for $E_\epsilon(V)$, with
\[
\beta(t;\bm\xi)
=b\sum_i t_i^2\log t_i
-\frac12t^TB(\bm\xi)t.
\]
Its derivatives up to order two in $h$ are uniformly
bounded because $t_i\ge c>0$.

Finally, Taylor's formula, \eqref{tp:residualC2},
and Lemma~\ref{tp:scaleLS} yield
\[
E_\epsilon(V+\phi)-E_\epsilon(V)
=O_{C^2(s,h)}(\gamma^2).
\]
Together with
$\gamma^2\le C(\epsilon^2\rho^2+\rho^4)$,
this proves \eqref{tp:F} and \eqref{tp:R}.
The relation $\eta_\epsilon=o(\epsilon)$ follows
from \eqref{eq:exponential}.
\end{proof}

\begin{lem}\label{tp:curvature}
At the solution parameters, with the centers fixed,
\begin{equation}\label{tp:scaleblocks}
\begin{gathered}
F_{ss}=16\pi^2\epsilon\rho^2
       +O(\eta_\epsilon\rho^2),
\qquad
F_{sh}=O(\eta_\epsilon\rho^2),\\
F_{hh}\ge c\rho^2\Id.
\end{gathered}
\end{equation}
For $m\ge2$, choose the sphere coordinates so that
\[
S_\epsilon:=D_ht(h_\epsilon)\longrightarrow S,
\qquad S^TS=\Id,\qquad S^Tt=0.
\]
Then
\begin{equation}\label{tp:hhlimit}
\rho^{-2}F_{hh}
\longrightarrow c_3^2S^T\cN S.
\end{equation}
Moreover,
\begin{equation}\label{tp:scalepositive}
F''[(c_s,\mathbf h),(c_s,\mathbf h)]
\ge c\rho^2\bigl(\epsilon|c_s|^2+|\mathbf h|^2\bigr).
\end{equation}
When $m=1$, the $h$ terms are omitted.
\end{lem}

\begin{proof}
Since $\Psi=u_\epsilon$ at the solution parameters,
$F_s=F_h=0$.
Differentiating \eqref{tp:F} and using these identities gives
\begin{equation}\label{tp:exactcancel}
F_{ss}=2b\epsilon\rho^2+2R_s+R_{ss},
\qquad
F_{sh}=2R_h+R_{sh}.
\end{equation}
Thus the first two estimates in \eqref{tp:scaleblocks}
follow from \eqref{tp:R}.

Differentiating $t^TMt$ on the unit sphere and using
$Mt=\Lambda t$ gives \eqref{tp:hhlimit}.
By Lemma~\ref{lem:matrix}, its limiting matrix is
positive definite, which yields the estimate for $F_{hh}$.
Finally, absorbing the mixed term using
$\eta_\epsilon=o(\epsilon)$ proves
\eqref{tp:scalepositive}.
\end{proof}

\subsection{Translation directions}\label{sec:frame}
For each scale coordinate
$\alpha\in\{s,h_1,\ldots,h_{m-1}\}$, set
$W_\alpha=\partial_\alpha\Psi$.
At the solution parameters, Lemma~\ref{tp:scaleLS} gives
\begin{equation}\label{tp:scaletangents}
\begin{aligned}
W_s
&=-\sum_iPJ_{i0}+O_{H^1_0}(\gamma),\\
W_{h_a}
&=\sum_i\frac{\partial_{h_a}t_i}{t_i}PJ_{i0}
  +O_{H^1_0}(\gamma),
\qquad 1\le a\le m-1.
\end{aligned}
\end{equation}
Differentiating \eqref{tp:projected} against a smooth
local section of $X$, and using
$E_\epsilon'(u_\epsilon)=0$, yields the first identity below;
the second follows from the chain rule:
\begin{equation}\label{tp:Qorthscale}
Q(W_\alpha,z)=0\quad(z\in X),
\qquad
Q(W_\alpha,W_\beta)=F_{\alpha\beta},
\end{equation}
where $\alpha,\beta$ range over the scale coordinates.
\begin{lem}\label{tp:translations}
Let $T_{ik}=-\chi_i\delta_i\partial_k u$, where
$\chi_i$ are the cutoffs in Proposition~\ref{tp:upper},
and define
\begin{equation}\label{tp:correctedT}
\zeta_{ik}
=\mathcal A_{u,X}^{-1}\Pi_X\mathcal A_uT_{ik},
\qquad
\mathcal T_{ik}=T_{ik}-\zeta_{ik}.
\end{equation}
Then
\begin{gather}
\|\zeta_{ik}\|\le C\rho^2,
\qquad Q(\mathcal T_{ik},X)=0,
\qquad
\|\mathcal T_{ik}-PJ_{ik}\|\le C\gamma,
\label{tp:Tapprox}\\
Q(\mathcal T_{ik},\mathcal T_{jr})=O(\rho^4),
\qquad
Q(W_\alpha,\mathcal T_{ik})=O(\rho^3),
\label{tp:Tblocks}
\end{gather}
where $\alpha$ ranges over the scale coordinates.
\end{lem}

\begin{proof}
Since $L_u\partial_k u=0$, the function $L_uT_{ik}$
is supported in the cutoff annulus.
By \eqref{tp:outeru},
\[
\|L_uT_{ik}\|_{L^2(\Om)}\le C\rho^2,
\qquad
\|\mathcal A_uT_{ik}\|\le C\rho^2.
\]
Lemma~\ref{tp:scaleLS} and \eqref{tp:correctedT}
therefore give the bound for $\zeta_{ik}$ and
$Q(\mathcal T_{ik},X)=0$.

To compare $\mathcal T_{ik}$ with $PJ_{ik}$,
the projection estimates
\cite[(2.4), (2.9)]{P-R-V}, together with
\eqref{tp:projectionderivative}, give
\[
\|-\chi_i\delta_i\partial_kV-PJ_{ik}\|
\le C\rho^2.
\]
For the term containing $\phi$, integration by parts yields
\[
\begin{aligned}
\bigl|\langle-\chi_i\delta_i\partial_k\phi,
PJ_{j\beta}\rangle\bigr|
&=
\left|
\int_\Om
\delta_i\partial_k(3\chi_iU_j^2J_{j\beta})\phi\,dx
\right|\\
&\le C\|\phi\|_{L^4(\Om)}
\le C\gamma.
\end{aligned}
\]
Here the coefficient has uniformly bounded
$L^{4/3}$ norm by rescaling and comparability of the widths.
The Gram matrix estimate in
Lemma~\ref{tp:quantprojection} thus gives
\[
\|\Pi_{\cK}(T_{ik}-PJ_{ik})\|\le C\gamma.
\]
Since $\zeta_{ik}\in X$ and
$\Pi_X\mathcal A_u\mathcal T_{ik}=0$,
decomposing $\mathcal T_{ik}-PJ_{ik}$ in $X\oplus\cK$
and applying Lemma~\ref{tp:scaleLS} yields
\[
\|\mathcal T_{ik}-PJ_{ik}\|
\le C\bigl(
\|\Pi_{\cK}(T_{ik}-PJ_{ik})\|
+\|\mathcal A_uPJ_{ik}\|
\bigr)
\le C\gamma,
\]
where the last step uses \eqref{tp:Jacres}.
This proves \eqref{tp:Tapprox}.

The cutoff identity used in the proof of
Proposition~\ref{tp:upper} gives
$Q(T_{ik},T_{jr})=O(\rho^4)$.
Consequently, by the established orthogonality,
\[
Q(\mathcal T_{ik},\mathcal T_{jr})
=Q(T_{ik},T_{jr})-Q(\zeta_{ik},\zeta_{jr})
=O(\rho^4).
\]
Finally, \eqref{tp:scaletangents} implies
\[
\|W_\alpha\|_{L^2(\operatorname{supp}L_uT_{ik})}
\le C\rho.
\]
Using \eqref{tp:Qorthscale}, we obtain
\[
|Q(W_\alpha,\mathcal T_{ik})|
=|Q(W_\alpha,T_{ik})|
=\left|\int_\Om W_\alpha L_uT_{ik}\,dx\right|
\le C\rho^3,
\]
which completes the proof.
\end{proof}

\begin{lem}\label{tp:decomposition}
At the solution parameters, let
\[
Y=\operatorname{span}\{W_s,W_{h_a},\mathcal T_{ik}\}.
\]
Then $H^1_0(\Om)=Y\oplus X$, with uniformly bounded
associated projections, and $Q(Y,X)=0$.
For $m+1\le\ell\le6m$, let $w=w_{\ell,\epsilon}$
satisfy $D(w,w)=1$, and set
$\theta=\mu_{\ell,\epsilon}-1$.
The unique decomposition
\begin{equation}\label{tp:vdecomp}
w=c_sW_s+\sum_a h_aW_{h_a}
+\sum_{i,k}p_{ik}\mathcal T_{ik}+x,
\qquad x\in X,
\end{equation}
has uniformly bounded coefficients, and
\begin{equation}\label{tp:xsmall}
\|x\|_{H^1_0(\Om)}\le C|\theta|.
\end{equation}
\end{lem}

\begin{proof}
By \eqref{tp:scaletangents}, \eqref{tp:Tapprox},
and the Jacobi Gram matrix estimate,
$\Pi_{\cK}|_Y:Y\to\cK$ is uniformly invertible.
Indeed, after multiplication by $D_t$, the limiting
scale coefficient matrix is $(-t,S)$, whose columns
form a basis of $\mathbb R^m$; the tangent columns
are omitted when $m=1$.
This proves the direct sum and the coefficient bounds,
since $\|w\|^2=\mu_{\ell,\epsilon}\le C$.
The orthogonality follows from
\eqref{tp:Qorthscale} and \eqref{tp:Tapprox}.

Testing the eigenvalue equation on $X$ gives
\[
Q(x,z)=\theta D(w,z),
\qquad z\in X.
\]
The uniform inverse in Lemma~\ref{tp:scaleLS}
therefore yields
$\|x\|\le C|\theta|\|w\|\le C|\theta|$.
\end{proof}

\begin{lem}\label{tp:twosided}
For $m+1\le\ell\le5m$, with
$\theta_\ell=\mu_{\ell,\epsilon}-1$,
\begin{equation}\label{tp:thetabound}
|\theta_\ell|\le C\frac{\rho^4}{\epsilon}.
\end{equation}
\end{lem}

\begin{proof}
Write $w=y+x$ as in \eqref{tp:vdecomp}.
By \eqref{tp:scalepositive}, \eqref{tp:Tblocks},
and Young's inequality,
\begin{equation}\label{tp:formlower}
Q(y,y)\ge-C\frac{\rho^4}{\epsilon}|p|^2.
\end{equation}
Since $Q(Y,X)=0$, Lemma~\ref{tp:decomposition} gives
\[
\theta_\ell
=Q(w,w)=Q(y,y)+Q(x,x)
\ge-C\frac{\rho^4}{\epsilon}-C\theta_\ell^2.
\]
When $\theta_\ell<0$, the quadratic term is absorbed
using $\theta_\ell\to0$ from \eqref{tp:coarse}.
When $\theta_\ell\ge0$, the claim follows directly
from \eqref{tp:upperbound}.
\end{proof}

\begin{lem}\label{tp:coefficients}
For $m+1\le\ell\le5m$, the scale coefficients
in \eqref{tp:vdecomp} satisfy
\begin{equation}\label{tp:rates}
|h|\le C\rho,
\qquad
|c_s|\le C\frac{\rho}{\epsilon}.
\end{equation}
\end{lem}

\begin{proof}
Test $Q(w,z)=\theta_\ell D(w,z)$ with
$z=W_s,W_{h_a}$.
The $x$ terms vanish by $Q(Y,X)=0$.
Using \eqref{tp:scaleblocks}, \eqref{tp:Tblocks},
\eqref{tp:thetabound}, and
$\rho/\epsilon\to0$ from \eqref{eq:exponential},
we obtain
\[
|h|\le C\eta_\epsilon|c_s|+C\rho,
\qquad
|c_s|\le
C\frac{\eta_\epsilon}{\epsilon}|h|
+C\frac{\rho}{\epsilon}.
\]
Since $\eta_\epsilon=o(\epsilon)$,
substitution and absorption give \eqref{tp:rates}.
For $m=1$, only the second inequality is needed.
\end{proof}

\subsection{Proof of Theorem~\ref{thm:pure}}

We use the normalization \eqref{tp:Dnormal}.
For $m+1\le\ell\le5m$,
Lemmas~\ref{tp:translations}--\ref{tp:coefficients}
and \eqref{eq:exponential} give
\begin{equation}\label{tp:finalapprox}
\begin{gathered}
w_\ell
=\sum_{i,k}p_{\ell,i,k,\epsilon}PJ_{ik}
+\sum_i b_{\ell,i,0,\epsilon}PJ_{i0}
+r_{\ell,\epsilon},\\
|b_{\ell,0,\epsilon}|\le C\rho/\epsilon,
\qquad
\|r_{\ell,\epsilon}\|_{H^1_0(\Om)}\le C\gamma,
\end{gathered}
\end{equation}
where $b_{\ell,0,\epsilon}
=(b_{\ell,i,0,\epsilon})_{i=1}^m$.
After passing to a common subsequence,
$p_{\ell,i,k,\epsilon}\to b_{\ell,i,k}$.
Combining \eqref{tp:coregeneral} and
\eqref{tp:finalapprox}, we obtain
\begin{equation}\label{tp:Dprofiles}
\begin{gathered}
\delta_iw_\ell(\xi_i+\delta_i y)
\longrightarrow
\sum_{k=1}^4b_{\ell,i,k}Z_k(y)
\quad\text{in }C^1_{\mathrm{loc}}(\mathbb R^4),\\
g_1\sum_i b_{\ell,i}\cdot b_{j,i}
=\delta_{\ell j},
\qquad
b_{\ell,i}=(b_{\ell,i,k})_{k=1}^4.
\end{gathered}
\end{equation}
Since $Z_k\in L^2(\mathbb R^4)$ for $k\ge1$,
the projection estimates
\cite[(2.8)--(2.9)]{P-R-V} and
\eqref{tp:radialbounds} yield
\[
\|PJ_{ik}\|_2\le C\rho
\quad(1\le k\le4),
\qquad
\|PJ_{i0}\|_2\le C\rho\sqrt{1+s}.
\]
Thus, by \eqref{tp:finalapprox} and
\eqref{eq:exponential},
\[
\|w_\ell\|_2
\le C\rho\left(1+\frac{\rho}{\epsilon}\sqrt{1+s}\right)
+C\gamma
\le C\rho.
\]
Since $\delta_i\asymp\rho$ and
$\|w_\ell\|_{H^1_0(\Om)}^2=\mu_{\ell,\epsilon}\le C$,
the zero extensions of
\[
\widetilde w_{\ell,i,\epsilon}(y)
=\delta_iw_\ell(\xi_i+\delta_i y)
\]
are bounded in $H^1(\mathbb R^4)$.
Consequently, the convergence in \eqref{tp:Dprofiles}
also holds weakly in $H^1(\mathbb R^4)$.

Apply Lemma~\ref{tp:endpoint} to
$v_{\ell,\epsilon}=w_\ell/\|w_\ell\|_\infty$.
Every maximum point of $|w_\ell|$ lies in
$\bigcup_iB(\xi_i,R_0\rho)$ for some fixed $R_0$.
The local convergence in \eqref{tp:Dprofiles}
and the explicit translation fields therefore give
\begin{equation}\label{tp:supnorm}
\rho\|w_\ell\|_\infty
\longrightarrow
\mathfrak n_\ell
:=\frac{3\sqrt6}{4}
\max_i\frac{|b_{\ell,i}|}{t_i}>0.
\end{equation}
The Gram identity in \eqref{tp:Dprofiles}
bounds these limits uniformly away from zero.
Together with $\|w_\ell\|_2\le C\rho$, this implies
\begin{equation}\label{tp:translationL2}
\|v_{\ell,\epsilon}\|_2\le C\rho^2.
\end{equation}

Finally,
$\delta_i\|w_\ell\|_\infty\to t_i\mathfrak n_\ell$.
Dividing the weak $H^1$ and local $C^1$ limits by
this factor proves \eqref{tp:mainprofile}, with
\[
a_{\ell,i,k}
=\frac{b_{\ell,i,k}}{t_i\mathfrak n_\ell}.
\]
The Gram identity also shows that
\[
d_\ell=(t_i a_{\ell,i,k})_{i,k}
=\left(\frac{b_{\ell,i,k}}{\mathfrak n_\ell}\right)_{i,k}
\]
are nonzero and mutually orthogonal.
\qed

\section{Proof of Theorem~\ref{thm:exterior}}
\label{sec:exterior}

Following the approach of
\cite[Lemma~7.6, Corollary~7.8, and Lemma~7.9]{YZmulti},
we first estimate the remainder in the Green representation
and then determine its coefficients using the scale equations.
Fix $m+1\le\ell\le5m$ and work along a subsequence
given by Theorem~\ref{thm:pure}.
We suppress $\ell$ and $\epsilon$ when no confusion arises
and write $\theta=\mu_{\ell,\epsilon}-1$.
Let $p_{i,k,\epsilon}$ be the translation coefficients
in \eqref{tp:vdecomp}, with
\[
p_{i,k,\epsilon}\longrightarrow p_{i,k}=b_{\ell,i,k},
\qquad p_i=(p_{i,k})_{k=1}^4.
\]
By \eqref{tp:Dprofiles}--\eqref{tp:supnorm}, we have
\begin{align}
\mathfrak n_\epsilon
:=\rho\|w_\ell\|_\infty
&\longrightarrow
\mathfrak n
=\frac{3\sqrt6}{4}\max_i\frac{|p_i|}{t_i}>0,\qquad
a_{\ell,i,k}
=\frac{p_{i,k}}{t_i\mathfrak n},
\label{ep:normconvert}\\
g_1\sum_i|p_i|^2&=1,\qquad
\mathfrak n
=\left(g_1\sum_i t_i^2|a_{\ell,i}|^2\right)^{-1/2}.
\label{ep:nunique}
\end{align}

Using the cutoffs $\chi_i$ from Proposition~\ref{tp:upper},
and choose a fixed $r>0$ such that
\[
B_i=B(\xi_{i,\epsilon},r)\Subset\Om
\]
are pairwise disjoint and contain
$\operatorname{supp}\chi_i$ compactly.
Define
\begin{equation}\label{ep:moments}
\begin{aligned}
\kappa_{i0}
&=\int_{B_i}u^3\,dy,
\qquad
q_i=3\int_{B_i}u^2v\,dy,\\
p_i^{\rm mom}
&=3\int_{B_i}
(y-\xi_{i,\epsilon})u(y)^2v(y)\,dy.
\end{aligned}
\end{equation}
Following \cite[(7.14), (7.49)]{YZmulti},
introduce the auxiliary scalar
\begin{equation}\label{ep:temporaryscale}
\mathfrak h_\epsilon
:=\rho+\rho^{-2}\|q\|_\infty,
\qquad
q=(q_1,\ldots,q_m)^T.
\end{equation}
We also set
\begin{equation}\label{ep:weighteddefinitions}
\begin{gathered}
W_\rho(x)
:=\sum_{i=1}^m
\frac{\rho^2}{\rho^2+|x-\xi_{i,\epsilon}|^2},
\qquad
\|f\|_{*,\rho}
:=\sup_{x\in\Om}\frac{|f(x)|}{W_\rho(x)},\\
\mathscr T:=(-\Delta_D)^{-1}.
\end{gathered}
\end{equation}
\begin{lem}\label{ep:weightedlemma}
Let $m+1\le\ell\le5m$ and define
\begin{equation}\label{ep:weightedh}
\psi_{ik}=-\delta_i^2\chi_i\partial_k u_\epsilon.
\end{equation}
There exist uniformly bounded coefficients $c_{ik,\epsilon}$
such that
\begin{equation}\label{ep:weightedzeta}
\zeta_\epsilon
=v_{\ell,\epsilon}-\sum_{i,k}c_{ik,\epsilon}\psi_{ik}
\end{equation}
satisfies
\begin{equation}\label{ep:weightedconclusions}
\|\zeta_\epsilon\|_{*,\rho}\le C\mathfrak h_\epsilon,
\qquad
\|v_{\ell,\epsilon}\|_{C^1(K)}
\le C_K\rho^2\mathfrak h_\epsilon
\end{equation}
for every compact set
$K\Subset\Om\setminus\{\xi_1^0,\ldots,\xi_m^0\}$,
where $\mathfrak h_\epsilon$ is defined in
\eqref{ep:temporaryscale}.
\end{lem}

\begin{proof}
Following \cite[Definition~7.5]{YZmulti}, choose
$c_{ik,\epsilon}$ so that
\begin{equation}\label{ep:weightedorthogonality}
\int_{\R^4}\vartheta U^2Z_j\,
\zeta_\epsilon(\xi_i+\delta_i y)\,dy=0,
\qquad 1\le i\le m,\quad 1\le j\le4,
\end{equation}
where $\vartheta\in C_c^\infty(\R^4)$ is nonnegative,
radial, and nonzero.
Since
$\psi_{ik}(\xi_i+\delta_i y)\to Z_k(y)$
in $C^1_{\mathrm{loc}}(\R^4)$, the coefficient matrices
converge to $c_\vartheta\Id_4$, with $c_\vartheta>0$.
Thus the coefficients are uniformly bounded.
Moreover, \eqref{tp:endpointbound} and the interior
gradient estimates associated with \eqref{tp:udecay} give
\[
\|v_{\ell,\epsilon}\|_{*,\rho}
+\sum_{i,k}\|\psi_{ik}\|_{*,\rho}\le C.
\]

For $R\ge2$, set
$\mathcal O_{R,\rho}
=\Om\setminus\bigcup_i B(\xi_i,R\rho)$.
The argument for \cite[(7.17)--(7.18)]{YZmulti},
using \eqref{tp:udecay}, yields
\begin{align}
\mathscr T(3u^2W_\rho)
&\le CW_\rho,
\label{ep:weightedcritical}\\
\mathscr T(\mathbf1_{\mathcal O_{R,\rho}}3u^2W_\rho)
&\le CR^{-2}W_\rho.
\label{ep:weightedtail}
\end{align}
The argument for \eqref{tp:kernel2}, applied to the
regularized weight $W_\rho$, also gives
\begin{equation}\label{ep:weightedlinear}
\mathscr T(\epsilon W_\rho)
\le C_\Om\epsilon W_\rho.
\end{equation}
This replaces \cite[(7.19)]{YZmulti} at $N=4$, $q=1$.
Write $L_u=-\Delta-3u^2-\epsilon$ and
$\theta=\mu_{\ell,\epsilon}-1$.
Since $L_u(\partial_k u)=0$, we have
\[
L_u\zeta_\epsilon
=F_\epsilon
:=\theta(3u^2+\epsilon)v_{\ell,\epsilon}
-\sum_{i,k}c_{ik,\epsilon}E_{ik},
\qquad E_{ik}:=L_u\psi_{ik}.
\]
The functions $E_{ik}$ are supported on the fixed cutoff
annuli and satisfy $\|E_{ik}\|_\infty\le C\rho^3$,
by \eqref{tp:outeru}.
Since $W_\rho\ge c_\Om\rho^2$,
\eqref{ep:weightedcritical}--\eqref{ep:weightedlinear}
and \eqref{tp:thetabound} imply
\[
\|\mathscr T(|F_\epsilon|)\|_{*,\rho}
\le C(\rho+|\theta|)\le C\rho.
\]

Suppose, by contradiction, that
\[
n_\epsilon:=\|\zeta_\epsilon\|_{*,\rho},
\qquad
\frac{n_\epsilon}{\mathfrak h_\epsilon}\to\infty,
\qquad
Z_\epsilon:=\frac{\zeta_\epsilon}{n_\epsilon}.
\]
Then $\|Z_\epsilon\|_{*,\rho}=1$ and
\[
\left\|\mathscr T\left(
\frac{|F_\epsilon|}{n_\epsilon}\right)\right\|_{*,\rho}
\longrightarrow0.
\]
The compactness and kernel argument in
\cite[Lemma~7.6]{YZmulti}, using
\cite[Lemma~2.10]{CKL-1} and
\eqref{ep:weightedorthogonality}, gives, along a subsequence,
\[
Z_\epsilon(\xi_i+\delta_i y)
\longrightarrow b_iZ_0(y)
\quad\text{in }C^1_{\mathrm{loc}}(\R^4).
\]

To show that $b_i=0$, integration by parts gives
\[
\int_{B_i}3u^2\psi_{ik}
=\delta_i^2\int_{B_i}(\partial_k\chi_i)u^3
=O(\rho^5).
\]
Consequently,
\[
\int_{B(0,r/\delta_i)}
3[\delta_i u(\xi_i+\delta_i y)]^2
Z_\epsilon(\xi_i+\delta_i y)\,dy
=\frac{q_i+O(\rho^5)}{n_\epsilon\delta_i^2}
\longrightarrow0,
\]
where we used $\delta_i\asymp\rho$ and
$\|q\|_\infty\le\rho^2\mathfrak h_\epsilon$.
By \eqref{tp:udecay}, $|Z_\epsilon|\le W_\rho$,
and the separation of the centers, the integrand is
bounded by $C(1+|y|^2)^{-3}$ on the expanding ball.
Dominated convergence therefore yields
\[
0=b_i\int_{\R^4}3U^2Z_0=b_ic_3,
\]
and hence $b_i=0$.

For fixed $R$, it follows that
\[
\alpha_\epsilon(R)
:=\sup_{\bigcup_iB(\xi_i,R\rho)}
\frac{|Z_\epsilon|}{W_\rho}\longrightarrow0.
\]
Splitting the Green representation over these balls
and their complement, and using
\eqref{ep:weightedcritical}--\eqref{ep:weightedlinear},
we obtain
\[
1\le C\alpha_\epsilon(R)+CR^{-2}
+C\epsilon+o(1).
\]
Letting $\epsilon\to0$ and then $R\to\infty$
gives a contradiction. This proves the first estimate
in \eqref{ep:weightedconclusions}.

Finally, on a slightly larger compact neighborhood of $K$,
we have $W_\rho=O(\rho^2)$ and $\psi_{ik}=O(\rho^3)$.
Thus \eqref{ep:weightedzeta} and
$\mathfrak h_\epsilon\ge\rho$ give
$|v_{\ell,\epsilon}|\le C_K\rho^2\mathfrak h_\epsilon$.
Interior elliptic estimates for the eigenvalue equation
yield the second estimate in
\eqref{ep:weightedconclusions}.
\end{proof}

We next estimate the linear term in the Green representation.

\begin{lem}\label{ep:linearpotential}
For every compact set
$K\Subset\Om\setminus\{\xi_1^0,\ldots,\xi_m^0\}$,
\begin{equation}\label{ep:linearpotentialbound}
\bigl\|\mathscr T(\mu_{\ell,\epsilon}\epsilon
v_{\ell,\epsilon})\bigr\|_{C^1(K)}
\le C_K\epsilon\rho^2\mathfrak h_\epsilon
=o(\rho^2\mathfrak h_\epsilon).
\end{equation}
\end{lem}

\begin{proof}
Following \cite[Corollary~7.8]{YZmulti}, use the
decomposition \eqref{ep:weightedzeta}.
By \eqref{ep:weightedlinear}, Lemma~\ref{ep:weightedlemma},
and interior elliptic estimates,
\begin{equation}\label{ep:linearzetabound}
\|\mathscr T(\mu\epsilon\zeta_\epsilon)\|_{C^1(K)}
\le C_K\epsilon\rho^2\mathfrak h_\epsilon.
\end{equation}

For $P_{ik}:=\mathscr T(\mu\epsilon\psi_{ik})$,
integration by parts gives
\begin{equation}\label{ep:lineartranslationparts}
P_{ik}(x)
=\mu\epsilon\delta_i^2
\int_\Om
\bigl[\chi_i(y)\partial_{y_k}G(x,y)
+G(x,y)\partial_k\chi_i(y)\bigr]u(y)\,dy.
\end{equation}
By \eqref{tp:udecay}, $\int_\Om u\le C\rho$.
Choose a fixed neighborhood $K'$ of $K$ away from the
limiting centers. Splitting the integral into a small ball
about $\xi_i$, separated from $K'$, and its complement,
the local integrability of $G$ and $\nabla_yG$ gives
\[
\|P_{ik}\|_{L^\infty(K')}
\le C_K\epsilon\rho^3.
\]
Moreover,
$-\Delta P_{ik}=\mu\epsilon\psi_{ik}
=O_{K'}(\epsilon\rho^3)$ by \eqref{tp:outeru}.
Interior elliptic estimates therefore yield
\[
\|P_{ik}\|_{C^1(K)}
\le C_K\epsilon\rho^3.
\]
The conclusion follows from \eqref{ep:linearzetabound},
the boundedness of $c_{ik,\epsilon}$, and
$\mathfrak h_\epsilon\ge\rho$.
\end{proof}

We now derive the exterior expansions.

\begin{lem}\label{ep:multipole}
For every compact set
$K\Subset\Om\setminus\{\xi_1^0,\ldots,\xi_m^0\}$,
\begin{align}
u_\epsilon(x)
&=\sum_i\kappa_{i0}G(x,\xi_{i,\epsilon})
+o_{C^1(K)}(\rho),
\qquad
\kappa_{i0}=c_3\delta_i+o(\rho),
\label{ep:umulti}\\
\frac{v_{\ell,\epsilon}(x)}{\mu_{\ell,\epsilon}}
&=\sum_i
\left[
q_iG(x,\xi_{i,\epsilon})
+p_i^{\rm mom}\cdot\nabla_yG(x,\xi_{i,\epsilon})
\right]
+o_{C^1(K)}(\rho^2\mathfrak h_\epsilon),
\label{ep:vmulti}\\
q_i&=o(\rho^2),
\qquad
p_i^{\rm mom}=O(\rho^3).
\label{ep:momentbounds}
\end{align}
In particular, $\mathfrak h_\epsilon=o(1)$.
\end{lem}

\begin{proof}
By Lemma~\ref{tp:profilelemma},
Theorem~\ref{thm:pure}, and the decay estimates
\eqref{tp:udecay} and \eqref{tp:endpointbound},
rescaling and dominated convergence give
\[
\frac{\kappa_{i0}}{\delta_i}\longrightarrow c_3,
\qquad
\frac{q_i}{\delta_i^2}
\longrightarrow
3\sum_{k=1}^4a_{\ell,i,k}
\int_{\R^4}U^2Z_k=0.
\]
The last equality follows by oddness.
Together with \eqref{tp:outeru}, the first limit proves
\eqref{ep:umulti}; the second implies
$q_i=o(\rho^2)$ and $\mathfrak h_\epsilon=o(1)$.

The same decay estimates yield
\[
\int_{B_i}|y-\xi_i|\,3u^2|v|\,dy\le C\rho^3,
\qquad
\int_{B_i}|y-\xi_i|^2\,3u^2|v|\,dy
\le C\rho^4(1+s).
\]
The first bound proves the remaining assertion in
\eqref{ep:momentbounds}.

To prove \eqref{ep:vmulti}, use
\[
\frac{v}{\mu}
=\mathscr T(3u^2v)+\epsilon\mathscr T v
\]
and apply the Green-kernel expansion used in
\cite[Lemma~7.9]{YZmulti} and
\cite[Lemma~5.3]{CKL-1}.
If necessary, replace the integration balls by smaller
fixed balls separated from a neighborhood of $K$.
By Lemma~\ref{ep:weightedlemma} and
$\mathfrak h_\epsilon\ge\rho$, outside these smaller balls,
\[
u=O(\rho),
\qquad
v=O(\rho^2\mathfrak h_\epsilon),
\qquad
3u^2v=O(\rho^4\mathfrak h_\epsilon).
\]
Thus restoring the original coefficients $q_i$ and
$p_i^{\rm mom}$ changes the expansion by at most
$O_{C^1(K)}(\rho^4\mathfrak h_\epsilon)$.

On the smaller balls, Taylor expansion in the second
variable of $G$, together with the quadratic integral
bound above, gives an error
$O_{C^1(K)}(\rho^4(1+s))$.
The contribution from their complement is
$O_{C^1(K)}(\rho^4\mathfrak h_\epsilon)$, while
Lemma~\ref{ep:linearpotential} controls the linear term.
Hence the total error is bounded by
\[
C_K\bigl[
\epsilon\rho^2\mathfrak h_\epsilon
+\rho^4\mathfrak h_\epsilon
+\rho^4(1+s)
\bigr]
=o(\rho^2\mathfrak h_\epsilon),
\]
since $\mathfrak h_\epsilon\ge\rho$ and
$\rho(1+s)\to0$.
This proves \eqref{ep:vmulti}.
\end{proof}

Following \cite[(5.5)--(5.7)]{CKL-1}, set
$G_i=G(\cdot,\xi_i^0)$,
$\Gamma_i=(4\pi^2|\cdot-\xi_i^0|^2)^{-1}$, and
\begin{equation}\label{ep:harmonicparts}
f_t=\sum_jt_jG_j=t_i\Gamma_i+g_i,
\qquad
g_i=-t_iH(\cdot,\xi_i^0)+\sum_{j\ne i}t_jG_j.
\end{equation}
For $z\in\R^m$, similarly write
$f_z=\sum_jz_jG_j=z_i\Gamma_i+g_i^z$.
Both $g_i$ and $g_i^z$ are harmonic near $\xi_i^0$.

\begin{lem}\label{ep:pohozaev}
The limiting centers satisfy
\begin{equation}\label{ep:centerbalance}
-\frac{t_i}{2}\nabla\tau(\xi_i^0)
+\sum_{j\ne i}t_j\nabla_xG(\xi_i^0,\xi_j^0)=0,
\qquad
t^T\cP=0.
\end{equation}
\end{lem}

\begin{proof}
Multiplying the equation for $u$ by $\partial_k u$
and integrating over $B_i$ gives
\begin{equation}\label{ep:translationstress}
\int_{\partial B_i}
\left[
\partial_k u\,\partial_\nu u
-\frac12|\nabla u|^2\nu_k
+\left(\frac{u^4}{4}+\frac{\epsilon u^2}{2}\right)\nu_k
\right]\,dS=0.
\end{equation}
Divide by $c_3^2\rho^2$ and use \eqref{tp:outeru}.
By \eqref{ep:harmonicparts} and Green's identity,
the limiting boundary integral equals
$-t_i\partial_k g_i(\xi_i^0)$.
Since $t_i>0$, we obtain $\nabla g_i(\xi_i^0)=0$.
The identity
$\nabla_xH(\xi_i^0,\xi_i^0)=\nabla\tau(\xi_i^0)/2$
then proves the first assertion.
Finally, symmetry of $G$ and \eqref{eq:P} give
\[
(t^T\cP)_{(i,k)}
=t_i^2\partial_k g_i(\xi_i^0)=0.
\]
\end{proof}

Since $UZ_0\notin L^1(\R^4)$, the argument in
\cite[Corollary~7.7]{YZmulti} does not apply directly.
We instead relate $\int_{B_i}3u^2w$ to the dilation
coefficients.

\begin{lem}\label{ep:sourcemoments}
Suppose
\[
w=\sum_{j=1}^m\sum_{\beta=0}^4b_{j\beta}PJ_{j\beta}+e,
\qquad
|b|\le C,
\qquad
\|e\|_{H^1_0(\Om)}\le C\gamma.
\]
Then, for the fixed disjoint balls $B_i$,
\begin{equation}\label{ep:momentestimate}
\int_{B_i}3u^2w
=c_3\delta_i b_{i0}
+O\bigl(\rho\gamma+\rho^3(1+s)\bigr).
\end{equation}
\end{lem}

\begin{proof}
Set $V=\sum_jP_j$ and
$W=\sum_{j,\beta}b_{j\beta}PJ_{j\beta}$.
By \eqref{tp:udecay}, \eqref{tp:radialbounds},
and the projection estimates
\cite[(2.4), (2.8)--(2.9)]{P-R-V},
\[
\|u^2\|_{4/3}+\|(u+V)W\|_{4/3}\le C\rho.
\]
Since $u-V=\phi_\epsilon$ and $w-W=e$,
Hölder's inequality and \eqref{eq:branch} yield
\[
\int_{B_i}3u^2w
=\int_{B_i}3V^2W+O(\rho\gamma).
\]

On $B_i$, the same projection estimates give
\[
V=U_i+O(\rho),
\qquad
W=\sum_{\beta=0}^4b_{i\beta}J_{i\beta}+O(\rho).
\]
Using $|J_{i\beta}|\le CU_i$ and
\eqref{tp:radialbounds}, we therefore obtain
\[
\int_{B_i}3V^2W
=\sum_{\beta=0}^4b_{i\beta}
\int_{B_i}3U_i^2J_{i\beta}
+O\bigl(\rho^3(1+s)\bigr).
\]
Finally, scaling and oddness give
\[
\int_{B_i}3U_i^2J_{i0}
=c_3\delta_i+O(\rho^3),
\qquad
\int_{B_i}3U_i^2J_{ik}=0
\quad(1\le k\le4),
\]
where the first identity uses
$\int_{\R^4}3U^2Z_0=c_3$ and the omitted tail is
$O(\delta_i^3)$.
Combining these estimates proves
\eqref{ep:momentestimate}.
Its remainder is $o(\rho^2)$ by
\eqref{eq:exponential}.
\end{proof}

We next refine the estimate in the common scale direction
$s=\log\rho^{-1}$.
All $s$-derivatives are taken at fixed
$\epsilon,t,\bm\xi$ and evaluated at the solution parameters.

\begin{lem}\label{ep:commonlemma}
Let $\phi_\epsilon$ be the correction in
\eqref{tp:projected}, and set
\[
f_0=i^*(c_3f_t)\in H^1_0(\Om),
\]
where $f_t$ is defined in \eqref{ep:harmonicparts}.
Then
\begin{equation}\label{ep:philim}
\frac{\phi_\epsilon}{\epsilon\rho_\epsilon}\to f_0,
\qquad
\frac{\phi_\epsilon+\partial_s\phi_\epsilon}
{\epsilon\rho_\epsilon}\to0
\quad\text{strongly in }H^1_0(\Om).
\end{equation}
Moreover,
\begin{equation}\label{ep:commonmix}
Q(W_s,\mathcal T_{ik})=o(\epsilon\rho^3).
\end{equation}
\end{lem}

\begin{proof}
Write
\[
\mathcal A_V=\Id-i^*((3V^2+\epsilon)\,\cdot),
\qquad
\mathcal R=V-i^*(V^3+\epsilon V),
\qquad
\mathcal B(\phi)=i^*(3V\phi^2+\phi^3).
\]
The residual of the projected family is
\begin{equation}\label{ep:actualresidual}
\mathcal H
:=\Psi-i^*(\Psi^3+\epsilon\Psi)
=\mathcal R+\mathcal A_V\phi-\mathcal B(\phi),
\qquad
\Pi_X\mathcal H=0.
\end{equation}

By the projection expansions
\cite[(2.4)--(2.5)]{P-R-V},
\begin{equation}\label{ep:VLp}
\frac{V}{\rho}\longrightarrow c_3f_t
\quad\text{strongly in }L^{4/3}(\Om).
\end{equation}
The cubic estimates in the proof of
\eqref{tp:residualC2} give
\[
\left\|\sum_iU_i^3-V^3\right\|_{4/3}
+\left\|\partial_s\left(\sum_iU_i^3-V^3\right)
\right\|_{4/3}
\le C\rho^2.
\]
Furthermore, the explicit bubble formula and
\eqref{tp:projectionderivative} yield
\begin{equation}\label{ep:Vcancel}
\|V_s+V\|_{4/3}\le C\rho^2.
\end{equation}
Since
$\mathcal R=i^*(\sum_iU_i^3-V^3-\epsilon V)$
and $\rho/\epsilon\to0$, these estimates imply
\begin{equation}\label{ep:Rlimits}
\frac{\mathcal R}{\epsilon\rho}\to-f_0,
\qquad
\frac{\mathcal R_s+\mathcal R}{\epsilon\rho}\to0
\quad\text{strongly in }H^1_0(\Om).
\end{equation}

For every fixed $f\in H^1_0(\Om)$, we also have
\begin{equation}\label{ep:fixedoperators}
\begin{gathered}
(\mathcal A_V-\Id)f\to0,
\qquad
(\partial_s\mathcal A_V)f\to0,\\
\Pi_{\cK}f\to0,
\qquad
(\partial_s\Pi_{\cK})f\to0
\quad\text{strongly in }H^1_0(\Om).
\end{gathered}
\end{equation}
The first two limits follow by Hölder's inequality,
splitting into small fixed balls about the limiting centers
and their complement.
For the projection limits, let $Z$ be the column map of
the projected Jacobi fields. Its columns and those of $Z_s$
are bounded in $H^1_0$ and tend to zero in $L^2$,
by Lemma~\ref{tp:scaledifferentiation} and
\eqref{tp:radialbounds}; hence they converge weakly to zero
in $H^1_0$.
The formula
$\Pi_{\cK}=Z(Z^*Z)^{-1}Z^*$
and its derivative now give the last two limits.
The operators in \eqref{ep:fixedoperators} are uniformly
bounded.

By \eqref{tp:C2phi} and Sobolev's inequality,
\[
\|\mathcal B(\phi)\|_{H^1_0}
+\|\partial_s\mathcal B(\phi)\|_{H^1_0}
\le C\gamma^2=o(\epsilon\rho).
\]
As in the fixed-space argument in the proof of
Lemma~\ref{tp:scaleLS}, the operator
\[
\mathcal L_V:=\Pi_X\mathcal A_V+\Pi_{\cK}
\]
has a uniformly bounded inverse on $H^1_0(\Om)$.
Set $\Phi=\phi/(\epsilon\rho)$.
Equations \eqref{ep:actualresidual},
\eqref{ep:Rlimits}, and \eqref{ep:fixedoperators} give
\[
\mathcal L_V(\Phi-f_0)\longrightarrow0.
\]
Thus $\Phi\to f_0$ strongly in $H^1_0(\Om)$,
proving the first limit in \eqref{ep:philim}.

Differentiate $\Pi_X\mathcal H=0$ and
$\Pi_{\cK}\phi=0$ at the actual solution parameters.
Since $\mathcal H=0$ there, the term
$(\partial_s\Pi_X)\mathcal H$ vanishes.
Adding the original equations gives
\begin{equation}\label{ep:zetaprojected}
\begin{aligned}
\mathcal L_V(\phi_s+\phi)
={}&-\Pi_X\bigl[
\mathcal R_s+\mathcal R
+(\partial_s\mathcal A_V)\phi
-\partial_s\mathcal B(\phi)-\mathcal B(\phi)
\bigr]\\
&-(\partial_s\Pi_{\cK})\phi.
\end{aligned}
\end{equation}
After division by $\epsilon\rho$, the right-hand side
tends to zero by \eqref{ep:Rlimits},
\eqref{ep:fixedoperators}, and $\Phi\to f_0$.
The uniform inverse of $\mathcal L_V$ proves the second
limit in \eqref{ep:philim}.

Finally, on the cutoff annuli,
\eqref{tp:projectionderivative} gives $V_s+V=O(\rho^3)$.
Since $W_s+u=V_s+V+\phi_s+\phi$,
\eqref{ep:philim} and \eqref{eq:exponential} imply
\[
\|W_s+u\|_{L^2(\operatorname{supp}L_uT_{ik})}
=o(\epsilon\rho).
\]
The proof of Lemma~\ref{tp:translations} gives
$\|L_uT_{ik}\|_2=O(\rho^2)$, while integration by parts yields
\[
Q(u,T_{ik})
=-2\int_\Om u^3T_{ik}
=-\frac{\delta_i}{2}
\int_\Om(\partial_k\chi_i)u^4
=O(\rho^5).
\]
Using $Q(W_s,\zeta_{ik})=0$ from
\eqref{tp:Qorthscale}, we conclude that
\[
\begin{aligned}
Q(W_s,\mathcal T_{ik})
&=\int_\Om(W_s+u)L_uT_{ik}-Q(u,T_{ik})\\
&=o(\epsilon\rho^3)+O(\rho^5)
=o(\epsilon\rho^3),
\end{aligned}
\]
which proves \eqref{ep:commonmix}.
\end{proof}

Write $\mathbf h=(h_a)_{a=1}^{m-1}$ for the coefficients
in \eqref{tp:vdecomp}, and set
$z_\epsilon=S_\epsilon\mathbf h$.
When $m=1$, we set $z_\epsilon=0$.

\begin{lem}\label{ep:massclosure}
For $m+1\le\ell\le5m$, we have
\begin{equation}\label{ep:hclosure}
q=O(\rho^3),
\qquad
\mathfrak h_\epsilon=O(\rho),
\qquad
t^Tq=o(\rho^3).
\end{equation}
Moreover, $c_s=o(\rho)$.
\end{lem}

\begin{proof}
By Lemma~\ref{tp:coefficients},
$|\mathbf h|\le C\rho$.
Testing the eigenvalue equation with $W_s$ and using
\eqref{tp:Qorthscale}, \eqref{tp:scaleblocks},
\eqref{tp:thetabound}, and Lemma~\ref{ep:commonlemma},
we obtain
\[
\bigl(16\pi^2\epsilon\rho^2
+O(\eta_\epsilon\rho^2)\bigr)c_s
=O(\eta_\epsilon\rho^3)
+o(\epsilon\rho^3)
+O(\rho^4/\epsilon).
\]
Since $\eta_\epsilon=o(\epsilon)$ and
$\rho^2/\epsilon^2=o(\rho)$ by \eqref{eq:exponential},
division by $\epsilon\rho^2$ gives
\begin{equation}\label{ep:csrefine}
c_s=o(\rho).
\end{equation}

Next, \eqref{tp:scaletangents}, \eqref{tp:Tapprox},
and \eqref{tp:vdecomp} yield
\begin{equation}\label{ep:Jacexp}
\begin{gathered}
w=\sum_i
\left(-c_s+\frac{(z_\epsilon)_i}{t_{i,\epsilon}}\right)
PJ_{i0}
+\sum_{i,k}p_{i,k,\epsilon}PJ_{ik}+e_\epsilon,
\\
\|e_\epsilon\|_{H^1_0(\Om)}\le C\gamma.
\end{gathered}
\end{equation}
Applying Lemma~\ref{ep:sourcemoments} and dividing by
$\|w\|_\infty=\mathfrak n_\epsilon/\rho$, as in
\eqref{ep:normconvert}, we obtain
\begin{equation}\label{ep:premass}
\frac{q}{\rho^2}
=\frac{c_3}{\mathfrak n_\epsilon}
\bigl(-c_st_\epsilon+z_\epsilon\bigr)+o(\rho).
\end{equation}
Since $z_\epsilon=O(\rho)$ and
$\mathfrak n_\epsilon\to\mathfrak n>0$,
\eqref{ep:csrefine}--\eqref{ep:premass} imply
$q=O(\rho^3)$.
The definition \eqref{ep:temporaryscale} then gives
$\mathfrak h_\epsilon=O(\rho)$.

Finally, $t_\epsilon^Tz_\epsilon=0$ and
\eqref{ep:csrefine} imply
$t_\epsilon^Tq=o(\rho^3)$.
Since $t_\epsilon\to t$ and $q=O(\rho^3)$,
we also have $t^Tq=o(\rho^3)$, proving
\eqref{ep:hclosure}.
For later use, \eqref{ep:premass} simplifies to
\begin{equation}\label{ep:qcoordinate}
\frac{q}{\rho^3}
=\frac{c_3}{\mathfrak n_\epsilon}
\frac{z_\epsilon}{\rho}+o(1).
\end{equation}
\end{proof}
We now identify the limit of $\rho^{-3}q$.

\begin{lem}\label{ep:monopole}
For the vector $q$ defined in \eqref{ep:moments},
\begin{equation}\label{ep:qlimit}
\rho^{-3}q\longrightarrow
c_3\cN^\dagger\cP d_\ell.
\end{equation}
The limit is the unique vector $\beta\in t^\perp$
satisfying $\cN\beta=c_3\cP d_\ell$.
\end{lem}

\begin{proof}
If $m=1$, \eqref{ep:qcoordinate} gives
$q=o(\rho^3)$, while \eqref{ep:centerbalance}
implies $\cP=0$.
Thus it suffices to consider $m\ge2$.

Fix $z=Se_a\in t^\perp$ and use the harmonic
decompositions in \eqref{ep:harmonicparts}.
By \eqref{tp:Qorthscale},
\[
Q(W_{h_a},\mathcal T_{jk})
=Q(W_{h_a},T_{jk})
=\int_\Om W_{h_a}L_uT_{jk}.
\]
On the fixed cutoff annulus,
\eqref{tp:scaletangents},
\eqref{tp:projectionderivative}, and \eqref{tp:outeru}
give
\[
\frac{W_{h_a}}{\rho}\longrightarrow c_3f_z
\quad\text{in }L^2,
\qquad
\frac{u}{\rho}\longrightarrow c_3f_t
\quad\text{in }C^2.
\]
Choose a small fixed ball $D_j$ centered at $\xi_j^0$
such that the limiting cutoff equals one near
$\overline{D_j}$.
The cutoff identity in the proof of
Lemma~\ref{tp:translations}, followed by Green's identity,
therefore yields
\begin{equation}\label{ep:annularflux}
\lim_{\epsilon\to0}\rho^{-3}
Q(W_{h_a},\mathcal T_{jk})
=c_3^2t_j\int_{\partial D_j}
\left(
\partial_kf_t\,\partial_\nu f_z
-f_z\,\partial_\nu\partial_kf_t
\right)\,dS,
\end{equation}
where $\nu$ is the outward unit normal to $D_j$.

By the Green identities in
\cite[Lemma~6.2]{CKL-1} and
\eqref{ep:centerbalance}, the boundary integral equals
\[
-z_j\partial_kg_j(\xi_j^0)
-t_j\partial_kg_j^z(\xi_j^0)
=-t_j\partial_kg_j^z(\xi_j^0).
\]
Moreover, symmetry of $G$ and \eqref{eq:P} give
\[
t_j^2\partial_kg_j^z(\xi_j^0)
=(S^T\cP)_{a,(j,k)}.
\]
Consequently,
\begin{equation}\label{ep:mixedlimit}
\rho^{-3}Q(W_{h_a},\mathcal T_{jk})
\longrightarrow
-c_3^2(S^T\cP)_{a,(j,k)}.
\end{equation}

By Lemma~\ref{ep:massclosure},
$z_\epsilon=O(\rho)$ and
$t_\epsilon^Tz_\epsilon=0$.
Passing to a further subsequence, write
\[
\frac{z_\epsilon}{\rho}\longrightarrow z_0\in t^\perp.
\]
Test the eigenvalue equation with each $W_{h_a}$
and divide by $c_3^2\rho^3$.
The right-hand side tends to zero, since
$|\theta|/\rho^3=O(\rho/\epsilon)=o(1)$
by \eqref{tp:thetabound} and \eqref{eq:exponential}.
Using \eqref{tp:hhlimit}, \eqref{ep:csrefine},
\eqref{ep:mixedlimit}, and
$p_{j,k}=\mathfrak n t_ja_{\ell,j,k}$
from \eqref{ep:normconvert}, we obtain
\begin{equation}\label{ep:shaperow}
S^T\bigl(\cN z_0-\mathfrak n\cP d_\ell\bigr)=0.
\end{equation}
Both terms in parentheses belong to $t^\perp$:
$\cN t=0$ by Lemma~\ref{lem:matrix}, and
$t^T\cP=0$ by \eqref{ep:centerbalance}.
Since the columns of $S$ span $t^\perp$,
\begin{equation}\label{ep:shapeeq}
\cN z_0=\mathfrak n\cP d_\ell,
\qquad
t^Tz_0=0.
\end{equation}
The inverse on $t^\perp$ defined in \eqref{eq:Ndagger}
therefore gives
\[
z_0=\mathfrak n\cN^\dagger\cP d_\ell.
\]
This identifies the limit independently of the further
subsequence. Finally, \eqref{ep:qcoordinate} yields
\[
\frac{q}{\rho^3}
\longrightarrow
\frac{c_3}{\mathfrak n}z_0
=c_3\cN^\dagger\cP d_\ell,
\]
which proves \eqref{ep:qlimit}.
\end{proof}
\begin{lem}\label{ep:dipolelemma}
Along the subsequence of Theorem~\ref{thm:pure},
\begin{equation}\label{ep:dipolelimit}
\rho^{-3}p_{i,k}^{\rm mom}
\longrightarrow c_3t_i^3a_{\ell,i,k},
\qquad 1\le i\le m,\quad 1\le k\le4.
\end{equation}
\end{lem}

\begin{proof}
The rescaling argument in
\cite[(7.70)--(7.72)]{YZmulti} and
\cite[proof of Proposition~5.1]{CKL-1}
applies by Lemma~\ref{tp:profilelemma},
Theorem~\ref{thm:pure}, and \eqref{tp:endpointbound}.
Indeed, the rescaled integrands are bounded by
$C|y|(1+|y|^2)^{-3}\in L^1(\R^4)$.
Dominated convergence and the identity
\[
3\int_{\R^4}y_kU^2Z_r\,dy=c_3\delta_{kr},
\qquad 1\le k,r\le4,
\]
give \eqref{ep:dipolelimit}.
The identity follows by integration by parts from
$3U^2Z_r=-\partial_r(U^3)$.
\end{proof}

We can now prove Theorem~\ref{thm:exterior}.
\begin{proof}[Proof of Theorem~\ref{thm:exterior}]
Work along the chosen subsequence of
Theorem~\ref{thm:pure}, and fix
$K\Subset\Om\setminus\{\xi_1^0,\ldots,\xi_m^0\}$.
By Lemma~\ref{ep:massclosure},
$\mathfrak h_\epsilon=O(\rho)$, so the remainder in
\eqref{ep:vmulti} is $o_{C^1(K)}(\rho^3)$.
Dividing that expansion by $\rho^3$ and using
\eqref{ep:qlimit}, \eqref{ep:dipolelimit},
$\mu_{\ell,\epsilon}\to1$, and
$\xi_{i,\epsilon}\to\xi_i^0$, we obtain
\[
\rho^{-3}v_{\ell,\epsilon}
\longrightarrow
c_3\bigl(\cG\cN^\dagger\cP+\widetilde{\cG}\bigr)d_\ell
\quad\text{in }C^1(K),
\]
where the passage to the limiting centers follows from
the smoothness of $G$ away from the diagonal.
Since $K$ is arbitrary, this proves \eqref{ep:final}.
The assertion
$\cN^\dagger\cP d_\ell\in t^\perp$
follows from Lemma~\ref{ep:monopole}.
\end{proof}

\section{Proof of Theorem~\ref{thm:translation}}\label{sec:translation}
In this section we prove Theorem~\ref{thm:translation}.
Fix $m+1\le\ell\le5m$, retain $f_t$ from
\eqref{ep:harmonicparts}, and set
\[
V_\ell
=\bigl(\cG\cN^\dagger\cP+\widetilde{\cG}\bigr)d_\ell.
\]

\subsection{Boundary and integral estimates}

Fix $r>0$ such that the balls
$B(\xi_{i,\epsilon},2r)$ are pairwise disjoint and
compactly contained in $\Om$ for all small $\epsilon$.
Set
\[
B_{i,\epsilon}=B(\xi_{i,\epsilon},r),
\qquad
B_i^0=B(\xi_i^0,r).
\]
For functions smooth near $\partial B_i^0$, define
\[
\mathscr B_i(f,h)
=\int_{\partial B_i^0}
\bigl(\partial_\nu f\,h-f\,\partial_\nu h\bigr)\,dS,
\]
where $\nu$ denotes the outward unit normal.

\begin{lem}\label{tr:flux}
For $1\le i\le m$ and $1\le k\le4$,
\begin{equation}\label{tr:rows}
\begin{aligned}
t_i\mathscr B_i(\partial_kf_t,\cG)
&=(\cP^T)_{(i,k),\cdot},\\
t_i\mathscr B_i(\partial_kf_t,\widetilde{\cG})
&=(\cQ_4)_{(i,k),\cdot},
\end{aligned}
\end{equation}
with $\mathscr B_i$ acting componentwise on row vectors.
Moreover, if
\begin{equation}\label{tr:Fdef}
F_{\ell ik,\epsilon}
=\int_{\partial B_{i,\epsilon}}
\left(
\partial_\nu\partial_k u_\epsilon\,v_{\ell,\epsilon}
-\partial_k u_\epsilon\,\partial_\nu v_{\ell,\epsilon}
\right)\,dS,
\end{equation}
then, along the subsequence of Theorem~\ref{thm:exterior},
\begin{equation}\label{tr:Flimit}
t_{i,\epsilon}\rho_\epsilon^{-4}F_{\ell ik,\epsilon}
\longrightarrow
c_3^2(\cB_4d_\ell)_{(i,k)}.
\end{equation}
\end{lem}

\begin{proof}
Applying the Green identities in
\cite[Lemma~6.2]{CKL-1} to the entries of
$\cG$ and $\widetilde{\cG}$, and using
\eqref{ep:centerbalance}, gives \eqref{tr:rows}
by the definitions \eqref{eq:P} and \eqref{eq:Q4}.
This is the same computation as in
\cite[(8.1)]{YZmulti}.
Bilinearity and \eqref{eq:B4} therefore yield
\begin{equation}\label{tr:BV}
t_i\mathscr B_i(\partial_kf_t,V_\ell)
=(\cB_4d_\ell)_{(i,k)}.
\end{equation}

The convergences \eqref{tp:outeru} and \eqref{ep:final}
hold in $C^2$ and $C^1$, respectively, on a fixed
neighborhood of the limiting spheres.
Parameterizing $\partial B_{i,\epsilon}$ by
$x=\xi_{i,\epsilon}+r\omega$, we obtain
\[
\rho_\epsilon^{-4}F_{\ell ik,\epsilon}
\longrightarrow
c_3^2\mathscr B_i(\partial_kf_t,V_\ell).
\]
Together with $t_{i,\epsilon}\to t_i$ and
\eqref{tr:BV}, this proves \eqref{tr:Flimit}.
\end{proof}

We next determine the interior pairing which multiplies
$\mu_{\ell,\epsilon}-1$. This is the endpoint version of
\cite[Lemma 7.2]{YZmulti}. An energy estimate for the tail suffices,
so the calculation uses only the pointwise bounds \eqref{tp:udecay} and \eqref{tp:endpointbound} and local $C^1$ convergence.

\begin{lem}\label{tr:pairing}
Put
\begin{equation}\label{tr:Idef}
I_{\ell ik,\epsilon}=\int_{B_{i,\epsilon}}
(3u_\epsilon^2+\epsilon)\partial_ku_\epsilon\,v_{\ell,\epsilon}.
\end{equation}
Then
\begin{equation}\label{tr:L2}
\|\nabla u_\epsilon\|_2\le C,\qquad
\|v_{\ell,\epsilon}\|_2^2\le C\rho_\epsilon^4,
\end{equation}
and, along a core-profile subsequence,
\begin{equation}\label{tr:Ilimit}
I_{\ell ik,\epsilon}\longrightarrow-g_1 a_{\ell,i,k},\qquad
t_{i,\epsilon}I_{\ell ik,\epsilon}\longrightarrow-g_1(d_\ell)_{(i,k)},
\qquad g_1=3\int_{\R^4}U^2(\partial_1U)^2=\frac{32\pi^2}{5}.
\end{equation}
\end{lem}
\begin{proof}
The first estimate in \eqref{tr:L2} follows from
\eqref{eq:branch}; the second is \eqref{tp:translationL2}.
Consequently, the linear part of \eqref{tr:Idef} satisfies
\begin{equation}\label{tr:Ilinear}
\left|\epsilon\int_{B_{i,\epsilon}}\partial_ku_\epsilon v_{\ell,\epsilon}\right|
\le\epsilon\|\nabla u_\epsilon\|_2\|v_{\ell,\epsilon}\|_2
\le C\epsilon\rho^2=o(1).
\end{equation}

For fixed $R\ge1$, set
$A_{i,\epsilon}(R)=B_{i,\epsilon}\setminus B(\xi_{i,\epsilon},R\delta_{i,\epsilon})$.
Within $B_{i,\epsilon}$ the nearest center is $\xi_{i,\epsilon}$.
The width comparability in \eqref{eq:widths} and \eqref{tp:udecay} and \eqref{tp:endpointbound} imply
\[
\|u_\epsilon^2v_{\ell,\epsilon}\|_{L^2(A_{i,\epsilon}(R))}^2
\le C\rho^8\int_{R\delta_i}^{r}
\frac{r_1^3}{(\rho^2+r_1^2)^6}\,dr_1
\le CR^{-8}.
\]
Consequently,
\begin{equation}\label{tr:Itail}
\left|\int_{A_{i,\epsilon}(R)}3u_\epsilon^2\partial_ku_\epsilon v_{\ell,\epsilon}\right|
\le C R^{-4}.
\end{equation}
On the fixed rescaled core, let
$U_{i,\epsilon}(y)=\delta_i u_\epsilon(\xi_i+\delta_i y)$ and
$V_{\ell i,\epsilon}(y)=v_{\ell,\epsilon}(\xi_i+\delta_i y)$.
The change of variables contains no remaining power of $\delta_i$:
\[
\int_{B(\xi_i,R\delta_i)}3u_\epsilon^2\partial_ku_\epsilon v_{\ell,\epsilon}
=3\int_{B_R}U_{i,\epsilon}^2\partial_kU_{i,\epsilon}V_{\ell i,\epsilon}.
\]
First let $\epsilon\to0$ using \eqref{tp:mainprofile}, and then $R\to\infty$
using \eqref{tr:Itail}. Since $Z_h=-\partial_hU$ and
$U$ is radial, the limit is
\[
-3\sum_h a_{\ell,i,h}\int_{\R^4}U^2\partial_kU\partial_hU
=-g_1 a_{\ell,i,k}.
\]
Together with \eqref{tr:Ilinear}, this proves \eqref{tr:Ilimit}.
\end{proof}

\subsection{Proof of the eigenvalue expansion}
We now follow the order of \cite[Section 8]{YZmulti}. The preceding
lemmas supply all the estimates needed at the four-dimensional
endpoint.
\begin{proof}[Proof of the spectral assertions in
Theorem~\ref{thm:translation}]
Let $\epsilon_n\to0$ be arbitrary. By
Theorems~\ref{thm:pure} and \ref{thm:exterior}, we may pass
to a common subsequence along which the corresponding
eigenfunction limits hold for all $m+1\le\ell\le5m$.

Applying Green's identity to $v_{\ell,\epsilon}$ and
$\partial_k u_\epsilon$ on $B_{i,\epsilon}$, as in
\cite[(7.7)--(7.8)]{YZmulti}, gives
\begin{equation}\label{tr:Greenidentity}
F_{\ell ik,\epsilon}
=(\mu_{\ell,\epsilon}-1)I_{\ell ik,\epsilon}.
\end{equation}
Set
\[
\eta_{\ell,\epsilon}
:=\frac{\mu_{\ell,\epsilon}-1}{\rho_\epsilon^4}.
\]
Multiplying \eqref{tr:Greenidentity} by
$t_{i,\epsilon}\rho_\epsilon^{-4}$ and using
Lemmas~\ref{tr:flux} and \ref{tr:pairing}, we obtain
\begin{equation}\label{tr:component}
c_3^2(\cB_4d_\ell)_{(i,k)}+o(1)
=
\eta_{\ell,\epsilon}
\bigl[-g_1(d_\ell)_{(i,k)}+o(1)\bigr].
\end{equation}
Since $d_\ell\ne0$, choose $(i_0,k_0)$ such that
$(d_\ell)_{(i_0,k_0)}\ne0$. Then
\begin{equation}\label{tr:etaratio}
\eta_{\ell,\epsilon}
=
\frac{c_3^2(\cB_4d_\ell)_{(i_0,k_0)}+o(1)}
{-g_1(d_\ell)_{(i_0,k_0)}+o(1)}
\longrightarrow \eta_\ell\in\R.
\end{equation}
In particular, $\eta_{\ell,\epsilon}$ is bounded.
We may therefore pass to the limit in every component
of \eqref{tr:component}, obtaining
\begin{equation}\label{tr:limiteigen}
\cB_4d_\ell
=-\frac{g_1}{c_3^2}\eta_\ell d_\ell.
\end{equation}

By Theorem~\ref{thm:pure}, the vectors
$d_{m+1},\ldots,d_{5m}$ form an orthogonal basis of
$\R^{4m}$. Thus the numbers
$-g_1\eta_\ell/c_3^2$ exhaust the eigenvalues of $\cB_4$,
counted with multiplicity.
Since $\mu_{\ell,\epsilon}$ are ordered increasingly,
so are $\eta_\ell$. Consequently,
$-g_1\eta_\ell/c_3^2$ are ordered decreasingly and hence
\[
-\frac{g_1}{c_3^2}\eta_\ell=\nu_\ell,
\qquad
\cB_4d_\ell=\nu_\ell d_\ell.
\]
Using
\begin{equation}\label{tr:constant}
\frac{c_3^2}{g_1}=20\pi^2,
\end{equation}
we conclude that
$\eta_{\ell,\epsilon}\to-20\pi^2\nu_\ell$
along the chosen subsequence.

Since every sequence $\epsilon_n\to0$ admits a subsequence
along which
\[
\frac{\mu_{\ell,\epsilon_n}-1}{\rho_{\epsilon_n}^4}
\longrightarrow -20\pi^2\nu_\ell,
\]
and this limit is independent of the subsequence,
the subsequence criterion for convergence yields
\eqref{tr:mainshift} for the full family.
\end{proof}

To complete the proof of Theorem~\ref{thm:translation},
we establish \eqref{tr:Hessianidentity} by differentiating
$M(\bm\xi)$, following the Schur complement interpretation
in \cite[Remark~9.3, (9.8)]{YZmulti}, with $\cN^\dagger$
acting as the inverse of $\cN$ on $t^\perp$.
\begin{proof}[Proof of \eqref{tr:Hessianidentity}]
Since the least eigenvalue of $M(\bm\xi^0)$ is simple
by Lemma~\ref{lem:matrix}, its positive unit eigenvector
$t(\bm\xi)$ can be chosen smoothly near $\bm\xi^0$.
For $\alpha=(i,k)$, write
$M_\alpha=\partial_{\xi_{i,k}}M$ and
$M_{\alpha\beta}=\partial_{\xi_\beta}M_\alpha$.
All quantities below are evaluated at $\bm\xi^0$.

Differentiating $Mt=\Lambda t$ and $t^Tt=1$ gives
\[
\Lambda_\alpha=t^TM_\alpha t,
\qquad
t^Tt_\alpha=0.
\]
By \eqref{ep:centerbalance}, $t^TM_\alpha t=0$ for every
$\alpha$, so $\nabla\Lambda(\bm\xi^0)=0$ and
\[
\cN t_\alpha=-M_\alpha t,
\qquad
t_\alpha=-\cN^\dagger M_\alpha t.
\]
Define the matrices $\mathsf J$ and $\mathsf H_0$ by
\[
\mathsf J_\alpha=M_\alpha t,
\qquad
(\mathsf H_0)_{\alpha\beta}=t^TM_{\alpha\beta}t.
\]
Differentiating $\Lambda_\alpha=t^TM_\alpha t$ once more
therefore yields
\begin{equation}\label{tr:Hesscalc}
D^2\Lambda(\bm\xi^0)
=\mathsf H_0-2\mathsf J^T\cN^\dagger\mathsf J.
\end{equation}
The definitions \eqref{eq:M} and \eqref{eq:P},
together with \eqref{ep:centerbalance}, give the first
identity below, while the second follows from
\eqref{eq:M} and \eqref{eq:Q4}:
\begin{equation}\label{tr:PJ}
\cP=-\mathsf J\mathsf D,
\qquad
\cQ_4=-\frac12\mathsf D\mathsf H_0\mathsf D.
\end{equation}
Substituting \eqref{tr:PJ} into \eqref{eq:B4}
and using \eqref{tr:Hesscalc}, we obtain
\[
\cB_4
=
-\frac12\mathsf D
\bigl(\mathsf H_0-2\mathsf J^T\cN^\dagger\mathsf J\bigr)
\mathsf D
=
-\frac12\mathsf D\,D^2\Lambda(\bm\xi^0)\,\mathsf D.
\]
This proves \eqref{tr:Hessianidentity}.
\end{proof}

\begin{Rem}
Since $\mathsf D$ is invertible, \eqref{tr:Hessianidentity}
implies that $\cB_4$ and $-D^2\Lambda(\bm\xi^0)$ have
the same numbers of positive, negative and zero eigenvalues,
counted with multiplicity.
The coefficients in \eqref{tr:mainshift} depend on
the eigenvalues of
$\mathsf D D^2\Lambda(\bm\xi^0)\mathsf D$.
\end{Rem}

\section{Dilation eigenvalues and Morse indices}
\label{sec:dilation}

We now study the remaining $m$ eigenvalues converging
to $1$ and derive the Morse index conclusions.
Following \cite[Lemma~9.1]{YZmulti}, we first identify
their limiting profiles by orthogonality to the
translation modes.

\begin{lem}\label{ds:pure}
Along any subsequence for which \eqref{tp:coregeneral}
holds for all relevant indices, the coefficients satisfy
\[
b_{\ell,i,k}=0,
\qquad
5m+1\le\ell\le6m,\quad
1\le i\le m,\quad 1\le k\le4.
\]
Moreover, writing
$b_\ell=(b_{\ell,1,0},\ldots,b_{\ell,m,0})^T$, we have
\begin{equation}\label{ds:bGram}
g_1\,b_\ell\cdot b_j=\delta_{\ell j},
\qquad
5m+1\le\ell,j\le6m.
\end{equation}
In particular, $b_{5m+1},\ldots,b_{6m}$ form an
orthogonal basis of $\R^m$.
\end{lem}
\begin{proof}
Fix $5m+1\le\ell\le6m$.
By Theorem~\ref{thm:pure} and \eqref{tp:coregeneral},
$b_{j,i,0}=0$ for $m+1\le j\le5m$, and the vectors
$(b_{j,i,k})_{1\le i\le m,\,1\le k\le4}$
form an orthogonal basis of $\R^{4m}$.
Moreover, \eqref{tp:coregeneral} gives
\[
\sum_{i=1}^m\sum_{k=1}^4
b_{\ell,i,k}b_{j,i,k}=0,
\qquad m+1\le j\le5m.
\]
Thus $b_{\ell,i,k}=0$ for all $1\le i\le m$ and
$1\le k\le4$.
Formula \eqref{ds:bGram} then follows from
\eqref{tp:coregeneral}.
\end{proof}

We next use \eqref{tp:graphconstraint} to determine
the terms of order $\rho_\epsilon^2$, following
\cite[Lemma~9.2]{YZmulti}.

\begin{lem}\label{ds:leading}
As $\epsilon\to0$,
\begin{equation}\label{ds:zerolimit}
\frac{\mu_{5m+1,\epsilon}-1}{\rho_\epsilon^2}
\longrightarrow0,
\end{equation}
and \eqref{ds:relativeformula} holds.
For any subsequence such that \eqref{tp:coregeneral}
holds for all $m+1\le\ell\le6m$, the vectors $b_\ell$
satisfy \eqref{ds:relativevectors}.
In particular,
\[
\mu_{6m,\epsilon}\le1+C\rho_\epsilon^2
\]
for all small $\epsilon$.
\end{lem}

\begin{proof}
Let $\epsilon_n\to0$ be arbitrary. Pass to a subsequence
such that \eqref{tp:coregeneral} holds for all
$m+1\le\ell\le6m$.
By Lemma~\ref{ds:pure}, each vector $b_\ell$ with
$5m+1\le\ell\le6m$ is nonzero.
Proposition~\ref{tp:leadingconstraint} then gives,
after passing to a further subsequence,
\[
\eta_\ell
:=\lim_{\epsilon\to0}
\frac{\mu_{\ell,\epsilon}-1}{\rho_\epsilon^2}
\in\R,
\qquad
\eta_\ell b_\ell=20\pi^2\cS b_\ell.
\]
By \eqref{ds:bGram}, the vectors
$b_{5m+1},\ldots,b_{6m}$ form an orthogonal basis of
$\R^m$. Thus the numbers
$\eta_{5m+1},\ldots,\eta_{6m}$ are all the eigenvalues
of $20\pi^2\cS$, counted with multiplicity.
Since the eigenvalues $\mu_{\ell,\epsilon}$ are
ordered increasingly, Lemma~\ref{lem:matrix} gives
\[
\eta_{5m+1}=0,
\qquad
\eta_{5m+1+j}=20\pi^2\kappa_j,
\quad 1\le j\le m-1.
\]
Moreover, $\ker\cS=\Span\{\one\}$ and
\eqref{ds:bGram} imply
\[
b_{5m+1}=B\one,
\qquad
mg_1B^2=1.
\]
Together with the matrix equation above, this proves
\eqref{ds:relativevectors}.
These ordered limits do not depend on the chosen
subsequence. Since the initial sequence was arbitrary,
\eqref{ds:zerolimit} and \eqref{ds:relativeformula}
hold for the full family.
The last bound follows directly.
\end{proof}

We now determine the term of order
$\epsilon\rho_\epsilon^2$ in $\mu_{5m+1,\epsilon}$,
using the estimates in \eqref{tp:scaleblocks}.

\begin{lem}\label{ds:common}
The eigenvalue $\mu_{5m+1,\epsilon}$ satisfies
\eqref{ds:commonformula}.
\end{lem}

\begin{proof}
Set $w=w_{5m+1,\epsilon}$ and
$\theta=\mu_{5m+1,\epsilon}-1$.
Let $\epsilon_n\to0$ be arbitrary, and pass to a
subsequence such that Lemma~\ref{ds:leading} gives
$b_i=B$ for all $i$, where $B\ne0$.
By Lemma~\ref{tp:decomposition}, the coefficients in
\eqref{tp:vdecomp} are bounded and
$x\to0$ in $H^1_0(\Omega)$.
Comparing this decomposition with the limits of $w$,
using \eqref{tp:scaletangents} and \eqref{tp:Tapprox},
we obtain
\begin{equation}\label{ds:coefficientrelation}
-c_s+\frac{(S_\epsilon h)_i}{t_{i,\epsilon}}
=B+o(1),
\qquad 1\le i\le m.
\end{equation}
Multiplying by $t_{i,\epsilon}^2$ and summing,
we use $|t_\epsilon|=1$ and
$t_\epsilon^TS_\epsilon=0$ to get
\begin{equation}\label{ds:cs}
c_s\longrightarrow-B.
\end{equation}
Moreover, Lemma~\ref{tp:quantprojection},
\eqref{tp:coregeneral} and \eqref{tp:scaletangents}
give
\begin{equation}\label{ds:denominator}
D_\epsilon(w,W_s)\longrightarrow-mg_1B\ne0.
\end{equation}

Testing the eigenvalue equation with $W_s$ and using
\eqref{tp:Qorthscale}, we obtain
\begin{equation}\label{ds:commonrow}
\theta D_\epsilon(w,W_s)
=
c_sF_{ss}
+\sum_a h_aF_{h_as}
+\sum_{i,k}p_{ik}Q_\epsilon(\mathcal T_{ik},W_s).
\end{equation}
By \eqref{tp:scaleblocks} and \eqref{ds:cs},
\[
\frac{c_sF_{ss}}{\epsilon\rho_\epsilon^2}
\longrightarrow-16\pi^2B.
\]
Since $h$ and $p$ are bounded,
\eqref{tp:scaleblocks}, \eqref{tp:Tblocks} and
\eqref{eq:exponential} also give
\begin{align}
\frac{\left|\sum_a h_aF_{h_as}\right|}
{\epsilon\rho_\epsilon^2}
&\le C\frac{\eta_\epsilon}{\epsilon}
\longrightarrow0,
\label{ds:relativeerror}\\
\frac{\left|\sum_{i,k}p_{ik}
Q_\epsilon(\mathcal T_{ik},W_s)\right|}
{\epsilon\rho_\epsilon^2}
&\le C\frac{\rho_\epsilon}{\epsilon}
\longrightarrow0.
\label{ds:translationerror}
\end{align}
Thus \eqref{ds:commonrow} and
\eqref{ds:denominator} yield
\[
\frac{\theta}{\epsilon\rho_\epsilon^2}
\longrightarrow
\frac{-16\pi^2B}{-mg_1B}
=\frac5{2m}.
\]
This limit does not depend on the chosen subsequence.
Since the initial sequence was arbitrary,
\eqref{ds:commonformula} holds as $\epsilon\to0$.
\end{proof}
\begin{Rem}
The estimate for $F_{ss}$ in \eqref{tp:scaleblocks}
has an error of $o(\epsilon\rho_\epsilon^2)$.
This is needed to determine the coefficient in
\eqref{ds:commonformula}; an $o(\rho_\epsilon^2)$
error alone would not be enough.
\end{Rem}

We next give the proof of Theorem~\ref{thm:dilation}.
\begin{proof}[Proof of Theorem~\ref{thm:dilation}]
Lemmas~\ref{ds:pure} and \ref{ds:leading} give
\eqref{ds:Dprofile} and \eqref{ds:relativevectors}.
The expansions \eqref{ds:relativeformula} and
\eqref{ds:commonformula} follow from
Lemmas~\ref{ds:leading} and \ref{ds:common}, respectively.
Since $\kappa_j>0$ and $5/(2m)>0$, we have
$\mu_{\ell,\epsilon}>1$ for $5m+1\le\ell\le6m$
and all small $\epsilon$.
For $m=1$, only \eqref{ds:commonformula} is needed.
\end{proof}

We next write these limits for
$v_{\ell,\epsilon}
=w_{\ell,\epsilon}/\|w_{\ell,\epsilon}\|_\infty$,
as in \cite[Theorem~1.8]{YZmulti}.

\begin{cor}\label{ds:supnorm}
Fix a subsequence such that \eqref{ds:Dprofile}
holds for all $5m+1\le\ell\le6m$, and set
\begin{equation}\label{ds:n}
\mathfrak n_\ell
=2\sqrt2\max_{1\le i\le m}\frac{|b_{\ell,i}|}{t_i}>0,
\qquad
a_{\ell,i}=\frac{b_{\ell,i}}{t_i\mathfrak n_\ell},
\qquad
f_\ell=(t_i a_{\ell,i})_{i=1}^m.
\end{equation}
Then, for these indices and every $1\le i\le m$,
\begin{equation}\label{ds:vprofile}
v_{\ell,\epsilon}
(\xi_{i,\epsilon}+\delta_{i,\epsilon}y)
\rightharpoonup a_{\ell,i}Z_0(y)
\quad\hbox{in }\cD(\R^4),
\end{equation}
and the convergence is strong in
$C^1_{\mathrm{loc}}(\R^4)$.
The vectors $f_{5m+1},\ldots,f_{6m}$ form an
orthogonal basis of $\R^m$ and satisfy
\[
f_{5m+1}\in\Span\{\one\},
\qquad
\cS f_{5m+1+j}=\kappa_jf_{5m+1+j},
\quad 1\le j\le m-1.
\]
\end{cor}

\begin{proof}
By Lemma~\ref{tp:endpoint}, the argument used to prove
\eqref{tp:supnorm} applies here.
Since $\|Z_0\|_\infty=2\sqrt2$, it gives
\[
\rho_\epsilon\|w_{\ell,\epsilon}\|_\infty
\longrightarrow
2\sqrt2\max_{1\le i\le m}
\frac{|b_{\ell,i}|}{t_i}
=\mathfrak n_\ell.
\]
Thus
$\delta_{i,\epsilon}\|w_{\ell,\epsilon}\|_\infty
\to t_i\mathfrak n_\ell$.
Dividing \eqref{ds:Dprofile} by this factor proves
\eqref{ds:vprofile} and the convergence in
$C^1_{\mathrm{loc}}(\R^4)$.
Finally, $f_\ell=b_\ell/\mathfrak n_\ell$, so the
matrix equations and orthogonality follow from
\eqref{ds:relativevectors}.
\end{proof}

\begin{Rem}
The function $Z_0$ belongs to $\cD(\R^4)$ but not to
$L^2(\R^4)$.
This explains why the weak convergence in
\eqref{ds:Dprofile} and \eqref{ds:vprofile}
is stated in $\cD(\R^4)$.
\end{Rem}

We now prove Corollary~\ref{cor:morse}
using the eigenvalue estimates above.
\begin{proof}[Proof of Corollary~\ref{cor:morse}]
By Theorems~\ref{thm:amp} and \ref{thm:dilation},
the first $m$ eigenvalues are below $1$, while those
with $5m+1\le\ell\le6m$ are above $1$ for small
$\epsilon$.
Moreover, \eqref{tp:coarse} gives
$\mu_{6m+1,\epsilon}\to2$.
Since the eigenvalues are ordered, all eigenvalues
with $\ell\ge6m+1$ are also above $1$.

For $m+1\le\ell\le5m$, Theorem~\ref{thm:translation}
gives $\mu_{\ell,\epsilon}<1$ if $\nu_\ell>0$
and $\mu_{\ell,\epsilon}>1$ if $\nu_\ell<0$.
The remaining $n_0(\cB_4)$ eigenvalues therefore give
\[
m+n_+(\cB_4)
\le \ind(u_\epsilon)
\le \ind(u_\epsilon)+\nul(u_\epsilon)
\le m+n_+(\cB_4)+n_0(\cB_4).
\]
Since $\mathsf D$ is invertible,
\eqref{tr:Hessianidentity} implies
\[
n_+(\cB_4)=n_-\bigl(D^2\Lambda(\bm\xi^0)\bigr),
\qquad
n_0(\cB_4)=n_0\bigl(D^2\Lambda(\bm\xi^0)\bigr).
\]
This proves \eqref{eq:morse}.
The bound for the nullity and the conclusion for
an invertible Hessian follow from these inequalities.
\end{proof}

\begin{Rem}
If $\nu_\ell=0$, Theorem~\ref{thm:translation} gives
only $\mu_{\ell,\epsilon}-1=o(\rho_\epsilon^4)$.
This does not determine whether $\mu_{\ell,\epsilon}$
is below, equal to or above $1$.
\end{Rem}

When $m=1$, the above results give the following corollary.
\begin{cor}\label{cor:single}
Let $m=1$, write
$\delta_\epsilon=\delta_{1,\epsilon}$ and
$\xi^0=\xi_1^0$, and let
$\lambda_1\le\cdots\le\lambda_4$ be the eigenvalues
of $D^2\tau(\xi^0)$.
Then $\rho_\epsilon=\delta_\epsilon$,
$\nabla\tau(\xi^0)=0$, and
\begin{align*}
\mu_{1,\epsilon}
&=\frac13+\frac{8\pi^2}{3}\tau(\xi^0)\delta_\epsilon^2
+o(\delta_\epsilon^2),\\
\mu_{j+1,\epsilon}
&=1+10\pi^2\lambda_j\delta_\epsilon^4
+o(\delta_\epsilon^4),
\qquad 1\le j\le4,\\
\mu_{6,\epsilon}
&=1+\frac52\epsilon\delta_\epsilon^2
+o(\epsilon\delta_\epsilon^2).
\end{align*}
The term $\cG\cN^\dagger\cP$ in \eqref{ep:final}
vanishes.
If $D^2\tau(\xi^0)$ is invertible, then
\[
\ind(u_\epsilon)=1+n_-\bigl(D^2\tau(\xi^0)\bigr),
\qquad
\nul(u_\epsilon)=0
\]
for all small $\epsilon$.
\end{cor}

\begin{proof}
For $m=1$, we have $t_1=1$,
$M=\Lambda=\tau(\xi^0)$ and $\cN^\dagger=\cS=0$.
By \eqref{tr:Hessianidentity},
\[
\nabla\tau(\xi^0)=0,
\qquad
\cB_4=-\frac12D^2\tau(\xi^0).
\]
The three expansions follow from
Theorems~\ref{thm:amp}, \ref{thm:translation}
and \ref{thm:dilation}.
The remaining conclusions follow from
\eqref{ep:final} and Corollary~\ref{cor:morse}.
\end{proof}

	\vspace{3mm}
\noindent \textbf{Acknowledgments:}
Silvia Cingolani was supported by INdAM-GNAMPA. Minbo Yang was partially supported by the National Natural Science Foundation of China (No.~12471114) and the Natural Science Foundation of Zhejiang Province (No.~LZ26A010002). Shunneng Zhao was partially supported by the National Natural Science Foundation of China (No.~12401146) and the Natural Science Foundation of Zhejiang Province (No.~LMS25A010007). The authors partly used AI as an auxiliary tool for language editing, organization of arguments, and preliminary checking of mathematical derivations. All mathematical statements, proofs, references, and conclusions appearing in the manuscript were independently verified by the authors, who take full responsibility for the content. The authors sincerely thank Professor Angela Pistoia (Sapienza University of Rome) and Professor Giusi Vaira (University of Bari Aldo Moro) for their helpful comments and suggestions on this work.

	\vspace{3mm}
\noindent \textbf{Declarations of interest:} None.

	\vspace{3mm}
\noindent \textbf{Data availability:} No data was used for the research described in the article.


\small

\begin{thebibliography}{99}

\bibitem{Bahri-1988}
A. Bahri and J.-M. Coron,
 \emph{On a nonlinear elliptic equation involving the critical {S}obolev exponent: the effect of the topology of the domain}, Comm. Pure Appl. Math.
\textbf{41} (1988), 253--294.

\bibitem{B-L-R}
A. Bahri, Y. Y. Li and O. Rey,
\emph{On a variational problem with lack of compactness: the topological effect of the critical points at infinity}, Calc. Var. Partial Differ. Equ.
\textbf{3} (1995), 67--93.

\bibitem{Brezis-1983}
H. Brezis and L. Nirenberg,
 \emph{Positive solutions of nonlinear elliptic equations involving critical Sobolev exponents}, Comm. Pure Appl. Math.
\textbf{36} (1983), 437--477.

\bibitem{BP}
H. Brezis and L. A. Peletier,
\emph{Asymptotics for elliptic equations involving critical growth},
in: Partial Differential Equations and the Calculus of Variations, Vol.~I,
Birkh\"auser, Boston, 1989, pp.~149--192.

\bibitem{CLP}
D. Cao, P. Luo and S. Peng,
\emph{The number of positive solutions to the Brezis--Nirenberg problem},
Trans. Amer. Math. Soc. \textbf{374} (2021), 1947--1985.

\bibitem{CKL-1}
W. Choi, S. Kim and K.-A. Lee,
\emph{Qualitative properties of multi-bubble solutions for nonlinear elliptic equations involving critical exponents},
Adv. Math. \textbf{298} (2016), 484--533.

\bibitem{DDM}
M. del Pino, J. Dolbeault and M. Musso,
\emph{The Brezis--Nirenberg problem near criticality in dimension 3},
J. Math. Pures Appl. \textbf{83} (2004), 1405--1456.

\bibitem{GT}
D. Gilbarg and N. S. Trudinger,
\emph{Elliptic Partial Differential Equations of Second Order},
Classics in Mathematics, Springer, Berlin, 2001, reprint of the 1998 edition.

\bibitem{GP05}
M. Grossi and F. Pacella, \emph{On an eigenvalue problem related to the critical exponent}, Math. Z. \textbf{250} (2005), 225--256.

\bibitem{HANZCHAO}
Z.-C. Han,
\emph{Asymptotic approach to singular solutions for nonlinear elliptic equations involving critical Sobolev exponent},
Ann. Inst. H. Poincar\'e Anal. Non Lin\'eaire \textbf{8} (1991), 159--174.

\bibitem{HLWW}
R. He, X. Liu, J. Wei and Y. Wu,
\emph{A refined blow-up analysis of the Brezis--Nirenberg equation and its application: The one-bubble case for $N\geq4$},
preprint, \href{https://arxiv.org/abs/2607.11132}{arXiv:2607.11132}, 2026.

\bibitem{K-L-2024}
T. K\"{o}nig and P. Laurain, \emph{Fine multibubble analysis in the higher-dimensional Brezis--Nirenberg problem}, Ann. Inst. H. Poincar\'{e} C Anal. Non Lin\'{e}aire \textbf{41} (2024), no. 5, 1239--1287.

\bibitem{K-L-2026}
T. K\"onig and P. Laurain,
\emph{Corrigendum to ``Fine multibubble analysis in the higher-dimensional Brezis--Nirenberg problem''},
Ann. Inst. H. Poincar\'e C Anal. Non Lin\'eaire \textbf{43} (2026), no.~3, 749--758.

\bibitem{LVWW}
F. Li, G. Vaira, J. Wei and Y. Wu,
\emph{Construction of bubbling solutions of the Brezis--Nirenberg problem in general bounded domains (I): The dimensions $4$ and $5$},
J. Lond. Math. Soc. \textbf{112} (2025), no. 2, Paper No. e70246, 39 pp.

\bibitem{Musso-Pistoia-2002}
M. Musso and A. Pistoia,
\emph{Multispike solutions for a nonlinear elliptic problem involving the critical Sobolev exponent}, Indiana Univ. Math. J. \textbf{51} (2002), 541--579.

\bibitem{MS-2018}
M. Musso and D. Salazar,
\emph{Multispike solutions for the Brezis--Nirenberg problem in dimension three},
J. Differential Equations \textbf{264} (2018), no. 11, 6663--6709.

\bibitem{P-R-V}
A. Pistoia, G. M. Rago and G. Vaira,
\emph{Multi-bubble solutions for the Brezis--Nirenberg problem in four dimensions}, Calc. Var. Partial Differential Equations \textbf{65} (2026), no. 2, Paper No. 37, 25 pp.

\bibitem{PR}
A. Pistoia and S. Rocci,
\emph{The Brezis--Nirenberg problem in 4D},
Discrete Contin. Dyn. Syst. Ser. S \textbf{17} (2024), no. 4, 1562--1572.

\bibitem{Pohozaev-1965}
S. I. Poho\v{z}aev,
\emph{Eigenfunctions of the equation $\Delta u+\lambda f(u)=0$},
Soviet Math. Dokl. \textbf{6} (1965), 1408--1411.

\bibitem{Rey-1989}
O. Rey, \emph{Proof of two conjectures of H. Brezis and L. A. Peletier},
Manuscripta Math. \textbf{65} (1989), 19--37.

\bibitem{Rey-1990}
O. Rey, \emph{The role of the Green's function in a nonlinear elliptic equation involving the critical Sobolev exponent}, J. Funct. Anal. \textbf{89} (1990), 1--52.

\bibitem{T-1}
F. Takahashi,
\emph{An eigenvalue problem related to blowing-up solutions for a semilinear elliptic equation with the critical Sobolev exponent}, Discrete Contin. Dyn. Syst. Ser. S \textbf{4} (2011), 907--922.

\bibitem{YZ}
M. Yang and S. Zhao,
\emph{Qualitative properties of blowing-up solutions of nonlinear elliptic equations with critical Sobolev exponent},
preprint, \href{https://arxiv.org/abs/2512.14401}{arXiv:2512.14401}, 2025.

\bibitem{YZmulti}
M. Yang and S. Zhao,
\emph{Spectral analysis of multi-bubble solutions to the Brezis--Nirenberg problem},
preprint, \href{https://arxiv.org/abs/2610.27752v1}{arXiv:2610.27752v1}, 2026.

\end{thebibliography}
\end{document}